\documentclass[11pt,a4paper,reqno]{amsart}
\usepackage[T1]{fontenc}
\usepackage{amsmath}
\usepackage{amsfonts}
\usepackage{amssymb}
\usepackage{graphicx}
\usepackage{amsbsy}
\usepackage{mathrsfs}
\usepackage{bbm}
\newtheorem{theorem}{Theorem}[section]
\newtheorem{lemma}[theorem]{Lemma}

\newtheorem{proposition}[theorem]{Proposition}

\theoremstyle{remark}
\newtheorem{remark}[theorem]{\it \bf{Remark}\/}

\numberwithin{equation}{section}
\catcode`@=11
\def\section{\@startsection{section}{1}%
  \z@{1.5\linespacing\@plus\linespacing}{.5\linespacing}%
  {\normalfont\bfseries\large\centering}}
\catcode`@=12
\newcommand{\be}{\begin{equation}}
\newcommand{\ee}{\end{equation}}
\newcommand{\bea}{\begin{eqnarray}}
\newcommand{\eea}{\end{eqnarray}}
\newcommand{\bee}{\begin{eqnarray*}}
\newcommand{\eee}{\end{eqnarray*}}

\renewcommand{\l}{\lambda}
\newcommand{\e}{\epsilon}
\renewcommand{\d}{\delta}
\renewcommand{\c}{\cdot}
\renewcommand{\a}{\alpha}
\renewcommand{\b}{\beta}

\newcommand{\D}{\Delta}
\newcommand{\n}{\nabla}
\newcommand{\p}{\partial}

\newcommand{\tu}{\tilde{u}}
\newcommand{\ttu}{\tilde{\tilde{u}}}
\newcommand{\tw}{\tilde{\e}}
\newcommand{\mq}{M^{(4)}}
\newcommand{\tq}{\widetilde{Q}}
\newcommand{\dt}{\frac{d}{dt}}
\newcommand{\ds}{\displaystyle}
\renewcommand{\t}{\Theta}
\newcommand{\s}{\Sigma}
\newcommand{\ts}{\widetilde{\Sigma}}
\renewcommand{\tt}{\widetilde{\Theta}}
\newcommand{\half}{\frac{1}{2}}

\def\R{{\bf R}}

\def\CC{\mathbb{C}}

\def\RR{\mathbb{R}}
\def\CC{\mathbb{C}}

\def\P{\mathcal{P}}

\def\L{\Lambda}
\def\ds{\displaystyle}

\def\e{\varepsilon}

\def\fref#1{{\rm (\ref{#1})}}

\def\R2+{\RR ^2_+}

\def\lsl{\frac{\lambda_s}{\lambda}}

\def\pa{\partial}

\def\lim{\mathop{\rm lim}}
\def\goto{\rightarrow}

\def\sup{\mathop{\rm sup}}

\def\e{\varepsilon}
\def\l{\lambda}

\def\log{{\rm log}}

\def\et{\tilde{\e}}
\def\te{\tilde{\e}}
\def\lsl{\frac{\lambda_s}{\lambda}}
\def\xsl{\frac{\alpha_s}{\lambda}}
\def\tgamma{{\tilde{\gamma}}}

\def\ut{\tilde{u}}

\def\n{\nabla}

\def\qp{Q_{\mathcal{P}}}
\def\pp{P_{\mathcal{P}}}

\def\S{\Sigma}
\def\T{\Theta}

\def\tpsi{\tilde{\psi}}
\def\pa{\partial}

\def\psit{\tilde{\psi}}

\def\zt{\tilde{\zeta}}

\def\gts{\tilde{\gamma}_s}

\def\wt{\tilde{w}}

\title[]{Existence and uniqueness of minimal blow up solutions to an inhomogeneous mass critical NLS}
\author[P. Rapha\"el]{Pierre Rapha\"el}
\address{Institut de Math\'ematiques de Toulouse, Universit\'e Paul Sabatier, France}
\email{pierre.raphael@math.univ-toulouse.fr}
\author[J. Szeftel]{Jeremie Szeftel}
\address{CNRS and DMA, ENS, Rue d'Ulm, Paris}
\email{jeremie.szeftel@ens.fr}
\begin{document}

\maketitle

\begin{abstract} We consider the 2-dimensional  focusing mass critical  NLS with an inhomogeneous nonlinearity: $i\partial_tu+\Delta u+k(x)|u|^{2}u=0$. From standard argument, there exists  a threshold $M_k>0$ such that $H^1$ solutions with $\|u\|_{L^2}<M_k$ are global in time while a finite time blow up singularity formation may occur for $\|u\|_{L^2}>M_k$. In this paper, we consider the dynamics at threshold $\|u_0\|_{L^2}=M_k$ and give a necessary and sufficient condition on $k$ to ensure the existence of critical mass finite time blow up elements. Moreover, we give a complete classification in the energy class of the minimal finite time blow up elements at a non degenerate point, hence extending the pioneering work by Merle \cite{M1} who treated  the pseudo conformal invariant case $k\equiv 1$.

\end{abstract}

\section{Introduction}

We deal in this paper with a two dimensional focusing mass critical nonlinear Schr\"odinger equation with an inhomogeneous nonlinearity:
\be
\label{nlsk:0}
(NLS) \ \ \left   \{ \begin{array}{ll}
    i\partial_tu=-\Delta u-k(x)|u|^{2}u, \ \ (t,x)\in [0,T)\times \RR^2,\\
    u(0,x)=u_0(x), \ \ u_0:\RR^2\to \CC,
         \end{array}
\right .
\ee
for some smooth bounded inhomogenity $k:\RR^2\goto \RR^*_+$. This kind of problem arises naturally in nonlinear optics for the propagation of laser beams. From the mathematical point of view, it is a canonical model to break the large group of symmetries of the $k\equiv 1$ homogeneous case.\\
From standard Cauchy theory and variational argument, \cite{GV}, \cite{W1}, there exists a critical number $M_k>0$ such that for $H^1$ initial data with $$\|u_0\|_{L^2}<M_k,$$ the unique $H^1$ solution to \fref{nlsk:0} is global in time and bounded in $H^1$, while on the contrary finite time blow up may occur for $\|u_0\|_{L^2}>M_k.$ We address in this paper the question of the dynamics at threshold $$\|u\|_{L^2}=M_k.$$ 

In the pioneering breakthrough work \cite{M1}, Merle proves in the homogeneous case $k\equiv 1$ the {\it existence and uniqueness} of a finite time blow up solution with critical mass. In this case, existence and uniqueness both rely dramatically on the explicit pseudo conformal symmetry\footnote{see \fref{vbvbeojr}} which is {\it lost} in the inhomogeneous case. Merle initiated in \cite{M2} the study of the inhomogeneous problem and obtained in particular a sufficient condition to ensure the {\it nonexistence} of critical blow up elements for the inhomogeneous problem \fref{nlsk:0}\footnote{see Theorem \ref{propblowup} below}.\\

In this paper, we give a necessary and sufficient condition on the inhomogeneity $k$ to ensure the existence of critical blow up elements, and prove the strong rigidity theorem of uniqueness of such minimal blow up elements in the energy class despite the absence of pseudo conformal symmetry.


\subsection{The homogeneous case}


Let us start with recalling some well-known facts in the homogeneous case $k\equiv 1$:
\be
\label{nls}
\left   \{ \begin{array}{ll}
    i\partial_tu=-\Delta u-|u|^{2}u, \ \ (t,x)\in [0,T)\times \RR^2,\\
    u(0,x)=u_0(x), \ \ u_0:\RR^2\to \CC.
         \end{array}
\right .
\ee
 $H^1$ solutions satisfy  the conservation of energy, $L^2$ norm and momentum:
 $$
 \left . \begin{array}{lll} L^2 \mbox{norm}: \
\ \int|u(t,x)|^2=\int|u_0(x)|^2,\\ \mbox{Energy}:\ \
E(u(t,x))=\frac{1}{2}\int|\nabla u(t,x)|^2-\frac{1}{2+\frac{4}{N}}\int
|u(t,x)|^{2+\frac{4}{N}}=E(u_0),\\ \mbox{Momentum}:\ \  M(u)=Im\left(\int\nabla u \overline{u}\right)=M(u_0),
         \end{array}
\right .
$$
and a large group of $H^1$ symmetries leaves the flow invariant: if $u(t,x)$ solves \fref{nls}, then $\forall (\lambda_0,t_0,x_0,\beta_0,\gamma_0)\in \RR_*^+\times\RR\times\RR^N\times\RR^N\times \RR$, so does 
\be
\label{symmetrygroup}
v(t,x)=\lambda_0^{\frac{N}{2}}u(t+t_0,\lambda_0 x+x_0-\beta_0t)e^{i\frac{\beta_0}{2}\cdot(x-\frac{\beta_0}{2} t)}e^{i\gamma_0}.
\ee
A last symmetry is not in the energy space $H^1$ but in the virial space $\Sigma=\{xu\in L^2\}\cap H^1$, the pseudo conformal transformation: if $u(t,x)$ solves (\ref{nls}), then so does 
\be
\label{vbvbeojr}
v(t,x)=\frac{1}{|t|^{\frac{N}{2}}}\overline{u}\left(\frac{1}{t},\frac{x}{t}\right)e^{i\frac{|x|^2}{4t}}.
\ee
Following \cite{GNN}, \cite{Kw}, let $Q$ be the unique $H^1$ nonzero positive radial solution to
\be\label{eqQ}
\D Q-Q+Q^3=0,
\ee
then the variational characterization of $Q$ ensures that initial data $u_0\in H^1$ with $\|u_0\|_{L^2}<\|Q\|_{L^2}$ yield global and bounded solutions $T=+\infty$, \cite{W1}.  On the other hand, finite time blow up may occur for data $\|u_0\|_{L^2}\geq\|Q\|_{L^2(\RR^2)}$. At the critical mass threshold, the {\it pseudo conformal symmetry} \fref{vbvbeojr} applied to the periodic solitary wave solution $u(t,x)=Q(x)e^{it}$ yields a {\it minimal mass blow-up solution}:
\be\label{nls:1}
S(t,x)=\frac{1}{t}Q\left(\frac{x}{t}\right)e^{i\frac{|x|^2}{4t}-\frac{i}{t}}, \ \ \|S(t)\|_{L^2}=\|Q\|_{L^2}
\ee
which blows up at $t=0$. In \cite{M1}, Merle proves the {\it uniqueness} of the critical mass blow up solution: a solution $u\in H^1$ with $\|u_0\|_{L^2}=\|Q\|_{L^2}$ and blowing up at $t=0$ is equal to $S(t)$ up to the symmetries of the flow. Through the conformal invariance, this results yields a complete dynamical classification of the solitary wave as the only nondispersive solution in $\Sigma$ with critical mass.


\subsection{The inhomogeneous case}


We now come back to \eqref{nlsk:0}. The conservation of momentum {\it no longer holds} but $L^2$ norm and energy are still conserved:
$$
 \left . \begin{array}{ll} L^2 \mbox{norm}: \
\ \int|u(t,x)|^2=\int|u_0(x)|^2,\\ \mbox{Energy}:\ \
E(u(t,x))=\frac{1}{2}\int|\nabla u(t,x)|^2-\frac{1}{2+\frac{4}{N}}\int
k(x)|u(t,x)|^{2+\frac{4}{N}}=E(u_0),\\         \end{array}
\right .
$$
The canonical effect of the inhomogenity is to completely destroy the group of symmetry \fref{symmetrygroup}\footnote{up to phase invariance}, and in this sense \eqref{nlsk:0} is a toy model to analyze the properties of NLS systems in the absence of symmetries.\\
From standard variational techniques and a virial type argument, Merle has derived in \cite{M2} the criterion of global existence for \eqref{nlsk:0}: given $\kappa>0$, let $Q_\kappa$ be
\be\label{defQk}
Q_\kappa(x)=\frac{1}{\kappa^{1/2}}Q(x),
\ee
and let $$k_2=\max_{x\in \RR^2}k(x)<+\infty,$$ then initial data with 
\be
\label{defmk}
\|u_0\|_{L^2}<M_k=\|Q_{k_2}\|_{L^2}
\ee 
yield global and $H^1$ bounded solutions while finite time blow up may occur for $\|u_0\|_{L^2}>M_k$. Moreover, the localization of the concentration point and a sufficient condition for the {\it non existence} of critical blow up elements are obtained:

\begin{theorem}[Merle, \cite{M2}]
\label{propblowup}
Assume that $k$ belongs to $C^1(\RR^2)$.\\
{\em 1. Localization of the concentration point}: Let $u\in H^1$ with $\|u\|_{L^2}=M_k$ be a solution to \eqref{nlsk:0}  blowing up at $T=0$. Assume that $\{x\in \RR^2 \ \ \mbox{such that} \ \ k(x)=k_2\}$ is finite. Assume also that there is $\delta>0$ and $R>0$ such that $k(x)\leq k_2-\delta$ for $|x|\geq R$.  Then there is $x_0\in\RR^2$ such that $k(x_0)=k_2$ and 
$$
|u(t)|^2\rightharpoonup \|Q_{k(x_0)}\|_{L^2}^2\delta_{x=x_0} \ \ \mbox{as} \ \ t\to 0.
$$
{\em 2. Criterion of non existence}: Assume that $k(x_0)=k_2$ and 
\be
\label{ceihgoihge}
\nabla k(x)\cdot (x-x_0)\leq -|x-x_0|^{1+\alpha_0} \ \ \mbox{near $x_0$ and for some} \ \ \alpha_0>0,
\ee
then there are no critical mass blow up solution at $x_0$.
\end{theorem}

In other words, the concentration of a possible minimal blow up solution must occur at a point $x_0$ where $k$ reaches its maximum, and hence in particular:
$$\nabla k(x_0)=0,$$ and the repulsivity condition \fref{ceihgoihge} in fact implies the non existence of such an object. Note however that \fref{ceihgoihge} implies that $k$ is not $\mathcal C^2$ at $x_0$, and hence Theorem \ref{propblowup} leaves completely open the case of a smooth $k$ which we address in this paper.\\

In the presence of a very smooth and flat $k$ at $x_0$, the existence of critical elements may be derived using a brutal fixed point argument, \cite{M2}, \cite{Merlemulti}, \cite{BGT}, see also \cite{KLR} for related statements in the setting of nonlocal nonlinearities. This argument has been recently sharpened by Banica, Carles, Duyckaerts \cite{BCD} for \eqref{nlsk:0} by adapting the method designed by Bourgain and Wang \cite{BW} and further revisited by Krieger and Schlag \cite{KS}: after linearizing the problem close to the explicit $S(t)$ approximate solution, one uses modulation theory and energy estimates to treat perturbatively the unstable modes and integrate the system backwards from the singularity \footnote{similarly like in scattering theory}. This allows one to lower the flatness of $k$ and indeed the existence of critical mass blow up solution is obtained under the flatness assumption: $$\nabla k(x_0)=\nabla^2k(x_0)=0.$$  
Let us say that in this approach, the problem is treated perturbatively from the homogeneous case and $x_0$ can be taken to be any point where $k$ is flat enough. The case of smooth $k$ with nondegenerate Hessian at $x_0$ i.e.  $$\nabla^2k(x_0)<0$$ is out of reach with these techniques because it generates large deformations  {\it which live at the same scaling like $S(t)$} as we shall see later.\\

Let us also say that even when existence is known, uniqueness in the energy class as proved in \cite{M1} for the homogeneous case does not follow. See \cite{Ba} for a pioneering investigation of this delicate problem in the setting of a domain. Here the strength and the weakness of Merle's pioneering proof for $k\equiv 1$ is that it fundamentally relies on the {\it pseudoconformal symmetry} \fref{vbvbeojr} which is lost here. The challenge is hence to provide a more {\it dynamical proof} of the classification theorem in the absence of symmetry. Eventually, let us recall that uniqueness is one of the keys towards the derivation of chaotic continuation after blow up results, \cite{M2}, which is a canonical open problem for inhomogeneous problems like \eqref{nlsk:0}.


\subsection{Statement of the result}


Let us now fix our assumptions on $k$ where we normalize without loss of generality the supremum $k_2=1$:\\

{\bf Assumption (H)}: $k$ is $\mathcal C^5\cap W^{1,\infty}$ with $$0<k_1\leq k(x)\leq 1 \ \ \mbox{and} \ \ \max_{x\in \RR^2} k(x)=1 \ \ \mbox{is attained}.$$
From \fref{defQk}, \fref{defmk}, the critical mass is then $$M_k=\|Q\|_{L^2}.$$
We start with the derivation of necessary conditions for the existence of a critical mass blow up solution which is based on a slight refinement of the variational techniques introduced in \cite{M2}:

\begin{proposition}[Necessary condition for the existence of a critical blow up element]
\label{prop:th2}
Let $u$ with $\|u\|_{L^2}=\|Q\|_{L^2}$ be a solution to \eqref{nlsk:0} which blows up at time $T=0$, then there exists $x_0\in \RR^2$ such that $$k(x_0)=1$$ and $u$ blows up at $x_0$ in the sense:
\be
\label{weakconxnot}
|u(t)|^2\rightharpoonup \|Q\|_{L^2}^2\delta_{x=x_0} \ \ \mbox{as} \ \ t\to 0.
\ee
Moreover, the energy $E_0$ of $u$ satisfies:
\be\label{condblowup}
E_0+\frac{1}{8}\int \nabla ^2k(x_0)(y,y)Q^4>0.
\ee
\end{proposition}

In order to classify critical mass blow up solutions, we hence pick $x_0\in \RR^2$ such that $k(x_0)=1$ and focus in the whole paper onto the case of a nondegenerate Hessian $$\nabla^2k(x_0)< 0$$ which we expect to be the most delicate one. We claim that the lower bound \fref{condblowup} is sharp, and the following theorem is the main result of this paper:

\begin{theorem}[Existence and uniqueness of a critical element at a nondegenerate critical point]
\label{th2}
Let $x_0\in \RR^2$ with $$k(x_0)=1 \ \ \mbox{and} \ \ \nabla^2k(x_0)<0.$$ Then for all $E_0$ satisfying \fref{condblowup}, there exists a unique up to phase shift $H^1$ critical mass blow up solution to \eqref{nlsk:0} which blows up at time $T=0$ and at the point $x_0$ in the sense of \fref{weakconxnot}, with energy $E_0$. Moreover, 
\be
\label{mommentgotozero}
\lim_{t\to 0}Im\left(\int\nabla u\overline{u}\right)=0.
\ee
\end{theorem}

In other words, for $k$ smooth, there exists a critical mass finite time blow up solution if and only if the supremum of $k$ is attained, and the corresponding minimal blow up elements at a non degenerate blow up point are completely classified.\\

{\bf Comments on Theorem \ref{th2}}\\

{\it 1. Structure of the minimal blow up elements}: Let 
\be
\label{constantcnot}
\tilde{E}_0=E_0+\frac{1}{8}\int \nabla^2k(0)(y,y)Q^4>0,\ \ C_0=\frac{\|yQ\|_{L^2}}{\sqrt{8\tilde{E}_0}}
\ee
and 
\be
\label{djhfeioyeioco}
S_{C_0}(t,x)=\frac{C_0}{|t|}Q\left(\frac{C_0(x-x_0)}{|t|}\right)e^{i\frac{|x-x_0|^2}{4t}-i\frac{C_0^2}{t}}
\ee be the exact pseudo conformal blow up solution to the homogeneous problem (\ref{nls}) with {\it modified energy} \fref{constantcnot},
 then up to a fixed phase shift, the critical mass blow up solution of Theorem \ref{th2} decomposes into a regular and a singular part 
$$
u(t)=S_{C_0}(t)+\tilde{u}(t)
$$
with $$\tilde{u}(t)\to 0\ \ \mbox{in} \ \ H^{1} \ \ \mbox{as} \ \ t\to 0.$$ In particular, the blow up speed is given by the conformal law: 
\be
\label{blowupconformal}
|\nabla u(t)|_{L^2}\sim\frac{1}{T-t}.
\ee
There are two fundamental remarks. First, the manifold of critical mass blow up solutions at a nondegenerate point is now strictly smaller than in the homogeneous case and is submitted to the constraint of the non trivial lower bound on the energy \fref{condblowup} and an {\it asymptotic zero momentum condition} \fref{mommentgotozero}. This non trivial structure is very much due to fine non trivial effects induced by the inhomogeneity $k$. In view of \fref{djhfeioyeioco}, the leading order effect of $k$ is to induce a pseudo conformal type blow up but parametrized by the modified energy \fref{constantcnot}. This shows in particular that the case $\nabla^2k(x_0)\neq 0$ cannot be treated perturbatively from the homogeneous case as for a given initial data -and hence a given level of energy-, the blow up profiles of minimal mass blow up solutions for \fref{nlsk:0} and \fref{nls} are strictly different due to the energy shift \fref{constantcnot}. Second and most importantly, the fact that the blow up speed is still given by the pseudo conformal law \fref{blowupconformal} is very much {\it a surprise}. Indeed, the nonvanishing condition $\nabla^2k(x_0)\neq 0$ induces deformations which leave at the same scaling like the leading terms driving the homogeneous pseudo conformal blow up. We shall see that  the derivation of the law \fref{blowupconformal} relies on an {\it algebraic cancellation}. Both these issues are detailed in the strategy of the proof below.\\

{\it 2. The degenerate case}: We expect the methods presented in this paper to also be able to treat the degenerate case $\nabla k(x_0)=\nabla ^2k(x_0)=0$ which should in fact be easier to handle. In particular, this becomes a perturbative problem with respect to the homogeneous case and we expect the manifold of minimal mass blow up solutions to be of the same size like in the homogeneous case. More precisely, we conjecture the following: for all $E_0>0$ and $M_0\in \RR^2$, there exists a unique up to phase shift critical mass blow up solution to \eqref{nlsk:0} which blows up at time $T=0$ and at the point $x_0$ in the sense of \fref{weakconxnot}, with energy $E_0$ and asymptotic kinetic momentum $$\lim_{t\to 0}\left(Im \int\nabla u\overline{u}\right)=M_0.$$ Note that the existence part for $M_0=0$ is proved in  \cite{BCD}.\\

{\it 3. More blow up solutions}: The existence of stable log-log type blow up dynamics for super critical mass initial data $\|u_0\|_{L^2}>\|Q\|_{L^2}$ can be derived by adapting the argument in \cite{MR4}, see \cite{PR}, \cite{BGR} for related results. Such solutions would blow up with speed $$|\nabla u(t)|_{L^2}\sim \sqrt{\frac{\log|\log(T-t)|}{T-t}}$$ and at {\it any} point $x_0\in \RR^2$ in the sense that $$|u(t)|^2\rightharpoonup \|Q_{k(x_0)}\|_{L^2}^2\delta_{x=x_0}+|u^*|^2 \ \ \mbox{as} \ \ t\to T \ \ \mbox{with} \ \ u^*\in L^2.$$ This is a consequence of the fact that {\it in the log-log regime}, the inhomogeneity $k$ can very much be treated perturbatively and the problem be essentially reduced to the homogeneous one. This is a manifestation of the stability of the log-log blow up, \cite{R1}. Theorem \ref{th2} shows that the structure of critical mass blow up solutions is of course much more fragile due in fact to its intrinsic instability. In the homogeneous case, there also exist excited pseudo conformal blow up solutions corresponding to the infinite sequence of $H^1$ solutions to \fref{eqQ} in dimension $N\geq 2$, which display a nondispersive behavior in the sense that all the mass is put into the singularity formation:
$$|u(t)|^2\rightharpoonup \|u_0\|_{L^2}^2\delta_{x=x_0}. $$ The existence and possible uniqueness properties of such objects with higher $L^2$ norm is open, see \cite{MR3} for some classification results about such solutions in the homogeneous case.\\

{\it 4. Extension to higher dimensions}: We have decided to restrict the analysis to dimension $N=2$ which is the physically relevant case. The case $N=1$ could be treated completely similarly. For higher dimensions, the only restriction of our analysis is of technical nature and one would have to face the problem of the lack of differentiability of the $L^2$ critical nonlinearity $u|u|^{\frac{4}{N}}$ for large $N$.\\

{\it 5. Connection to previous works}: This work lies within the question of classification of critical elements. This kind of question is related to  Liouville type theorems, \cite{KM}, and is connected to the question of the dynamical classification of solitary wave solutions which is only at its beginning,  see \cite{Marteluniqueness}, \cite{visanetal}. In particular, our techniques borrow from \cite{Marteluniqueness} the energy estimates type of strategy, and are in spirit related to the recent work by Duyckaerts and Merle \cite{DM1}  where deep rigidity properties of critical dynamics of the energy critical homogeneous NLS are proved. One main difference with these works however is the fact that we are dealing here with a slowly polynomially converging problem, while \cite{DM1}, \cite{Marteluniqueness} deal with exponential decay which somehow allows a more clear decoupling of the interactions.\\

Eventually, let  us say that our strategy both for existence and uniqueness does not rely at all on the existence of the pseudo conformal symmetry for the homogeneous problem. On the contrary, the singular dynamics is extracted directly in the renormalized variables and the presence of the symmetry only allows for some -unnecessary- algebraic simplifications. Let us recall that the striking example of the $L^2$ critical KdV shows that critical mass blow up dynamics do not always exist, \cite{MMduke}, and non unique\footnote{infinite time blow up} critical mass blow up dynamics are exhibited in \cite{Mas} for the critical Keller-Segel model. We expect our strategy to be robust enough to provide a canonical framework to address the question of existence and uniqueness of minimal blow up dynamics for various nonlinear dispersive equations in the absence of any explicit pseudo conformal symmetry.\\

{\bf Acknowledgments}: Both authors are supported by the French Agence Nationale de la Recherche, ANR jeunes chercheurs SWAP. P.R is also supported by ANR OndeNonLin.


\subsection{Strategy of the proof}


Let us give a brief insight into the strategy of the proof of Theorem \ref{th2}.\\

{\bf step 1} Refined blow up profile and formal derivation of the laws of motion.\\

Let us start with assuming the existence of a critical mass blow up solution at $T=0$ and $x_0=0$ and obtaining dynamical informations on such a solution. From standard variational argument, $u$ must admit near blow up time a decomposition of the form: $$u(t,x)=\frac{1}{\lambda^{\frac{N}{2}}(t)}\left(Q+\e\right)\left(t,\frac{x-\alpha(t)}{\lambda(t)}\right)e^{i\gamma(t)} \ \ \mbox{with} \ \ \|\e(t)\|_{H^1}\to 0 \ \ \mbox{as} \ t \to 0,$$ see \cite{W1}, \cite{M2}. In fact, following \cite{MR1}, \cite{MR2}, \cite{R1}, this decomposition should be sharpened by deforming the main nonlinear part in the direction of pseudo conformal blow up. More precisely, in the homogeneous case, the explicit pseudo conformal blow up solution \fref{nls:1} should be thought of as being the following object: $$u(t,x)=\frac{1}{\lambda^{\frac{N}{2}}(t)}Q_{b(t),\beta(t)}\left(t,\frac{x-\alpha(t)}{\lambda(t)}\right)e^{i\gamma(t)}  \ \ \mbox{with} \ \ Q_{b,\beta}(y)=Q(y)e^{-i\frac{b|y|^2}{4}}e^{i\beta\cdot y}$$ and where the geometrical parameters $(\lambda(t),b(t),\alpha(t),\beta(t))$ are driven by the dynamical system: 
\be
\label{dynamicalhomogeneous}
b_s=-b^2, \ \ -\lsl=b, \ \ \frac{\alpha_s}{\lambda}=2\beta, \ \ \beta_s=-b\beta, \ \ \gamma_s=-1
\ee
where $s$ is the global rescaled time $$\frac{ds}{dt}=\frac{1}{\lambda^2}.$$ Integrating the system with boundary condition $\lambda(0)=0$ -corresponding to blow up at $t=0$- and $\alpha(0)=0$ -corresponding to blow-up at $x_0=0$- involves three integration constants which can de directly related to the energy  and momentum of the initial data, and the choice of a fixed  phase shift.\\
The first step in the proof of Theorem \ref{th2} is to similarly decompose the solution as 
\be
\label{conioheo}
u(t,x)=\frac{1}{\lambda^{\frac{N}{2}}(t)}(Q_{\mathcal P}+\e)\left(t,\frac{x-\alpha(t)}{\lambda(t)}\right)e^{i\gamma(t)} \ \ \mbox{with} \ \ \mathcal{P}=(b,\lambda,\beta, \alpha)
\ee
where the profile $Q_{\mathcal{P}}$ takes fully into account the deformation induced by the $k$ inhomogeneity, and is an approximation of high order -$O(\lambda^5)$- of the renormalized equations. This construction generalizes the one introduced in \cite{MR2}, see also \cite{KMR} for a related approach in a different setting. {\it Under the assumption that the excess of mass $\e$ is of lower order in some suitable sense}, one may formally extract the leading order ODE driving the geometrical parameters which roughly takes the form:
\be
\label{odemain}
b_s=-b^2, \ \ \beta_s=-b\beta+c_1(\alpha)\lambda, \ \ \lsl=-b, \ \ \frac{\alpha_s}{\lambda}=2\beta, \ \ s(t)\sim \frac{1}{|t|}
\ee
where $c_1(\alpha)$ is a linear form in $\alpha$. The last two equations are standard. The first equation for $b$ is surprising as inspection reveals that we should expect from $\nabla^2k(0)\neq 0$:$$b_s=-b^2+c_0\lambda^2,$$ and $\lambda\sim b$ in the conformal regime. It is a specific {\it algebraic cancellation}\footnote{see \fref{explicitcomp}} which leads to $c_0=0$. The second law for $\beta$ reflects the difference between the degenerate  case for which $c_1=0$ and the non degenerate case for which $c_1(\alpha)\sim \alpha$. In this last case, the reintegration of the scaling law yields $$b(s)\sim \frac{1}{s}, \ \ \lambda(s)\sim \frac{1}{s}.$$ On the other hand, from \fref{conioheo},
\be
\label{contorlsolution}
\frac{\beta}{\lambda}\sim Im\left(\int \nabla u\overline{u}\right)\sim M_0
\ee 
and hence if the asymptotic momentum $M_0$ is non zero asymptotically, the reintegration of the translation parameter yields $$\alpha_t\sim M_0 \ \ \mbox{and thus} \ \ \alpha (t)\sim M_0|t|\sim \frac{M_0}{s}.$$ Injecting this into the momentum law of $\beta$ yields: $$\left(\frac{\beta}{\lambda}\right)_s=\frac{1}{\lambda}[\beta_s+b\beta]\sim c_1(\alpha)\sim \frac{M_0}{s}.$$ This generates a logarithmic divergence in $s$ by integration which contradicts  \fref{contorlsolution} for $M_0\neq 0$. Such a phenomenon is absent is the degenerate case $\nabla^2k(x_0)=0$ for which $c_1=0$, and hence the asymptotic value $M_0$ becomes a free parameter when integrating \fref{odemain} as in the homogeneous case.\\

{\bf step 2} Construction of a critical mass blow up solution.\\

The generalization to a critical setting of the subcritical strategy of slow variable ansatz developed in \cite{KMR} allows us to derive an approximate critical mass solution to a sufficiently high order. The derivation of an exact minimal mass blow up element then follows from a standard compactness argument using energy estimates and integrating the flow backwards from blow up time as initiated in \cite{M1}, see also \cite{Marteluniqueness}, \cite{MMmulti}, \cite{KMR}. Indeed, one decomposes the flow as $$u(t,x)=u_{app}(t,x)+\tilde{u}(t,x)$$ where $\tilde{u}$ roughly solves: 
\be
\label{cnkoheoh}
i\pa_t\ut+L_{\lambda}\ut=\mbox{l.o.t.}
\ee 
Here the linearized Hamiltonian $L_{\lambda}$ is after renormalization the matrix linearized operator close to the ground state $L=(L_+,L_-)$ with:
\be
\label{deflpluslmoins}
L_+=-\Delta +1-3Q^2, \ \ L_-=-\Delta+1-Q^2.
\ee
Similarly like in \cite{RodRaph}, a direct energy estimate on \fref{cnkoheoh} however fails and generates quadratic errors due to the time dependence of the Hamiltonian which cannot be treated perturbatively. We hence modify the energy functional by adding a local Morawetz type term and derive a monotonicity formula of the form:
\be
\label{firstmoneot}
\frac{d}{dt}\left\{\|\nabla\tilde{u}\|^2_{L^2}+\frac{\|\tilde{u}\|^2_{L^2}}{\lambda^2}+\Im\left(\int_{r\leq \lambda}\pa_r\tilde{u}\overline{\ut}\right)\right\}\geq 0+O(\mbox{l.o.t.}),
\ee
which for a sufficiently small error in the constructed approximation is enough to close the existence part.\\

{\bf step 3} Derivation of the blow up speed assuming dispersion.\\

We now turn to the proof of uniqueness. Standard variational arguments ensure that a critical mass blow up solution must admit a decomposition of the form \fref{conioheo}, but there holds no {\it a priori upper bound} on the blow up rate in this regime. In fact, if we rewrite \fref{conioheo} in the form 
$$u(t,x)=\frac{1}{\lambda^{\frac{N}{2}}(t)}Q_{\mathcal P}\left(t,\frac{x-\alpha(t)}{\lambda(t)}\right)e^{i\gamma(t)} +\tilde{u}.
$$
then standard variational estimates yield: $$\limsup_{t\to 0}\|\ut(t)\|_{H^1}\lesssim C, \ \ \lim_{t\to 0}\|\tilde{u}(t)\|_{L^2}=0.$$ The key to the derivation of the sharp blow up speed is to improve this energy bound for a {\it dispersive} estimate:
\be
\label{vndoheoergu}
\ut(t)\to 0\ \ \mbox{in} \ \ H^1.
\ee 
The proof of \fref{vndoheoergu} is the heart of the paper and relies on the derivation of a local Morawetz monotonicity formula which generalizes in some sense \fref{firstmoneot}\footnote{which is useless until we know the blow up speed} and provides another Lyapounov type rigidity in the problem in the form  
\be
\label{cniofheioyeoy}
\left\{-\frac{b}{\lambda}+\Im\left(\int _{|x|\lesssim \lambda}\partial_r\tu\overline{\tu}\right)\right\}_t\gtrsim \frac{1}{\lambda^3}\left(\int_{|x|\lesssim \lambda}|\nabla \tilde{u}|^2+\int|\tilde{u}|^2\right).
\ee Note that remarkably enough, the quadratic form involved in the above coercivity estimate is the linearized quadratic form of the {\it energy} for which the coercivity is known from the variational characterization of $Q$, and no further coercivity property will be needed in the paper\footnote{as opposed to the works \cite{MM1}, \cite{MR1}, \cite{FMR}}.\\
The control of the local interactions with the ground state part of the solution provided by \fref{cniofheioyeoy} coupled with the critical mass conservation laws will imply all together the dispersion \fref{vndoheoergu} and the derivation of the blow up speed:
$$b\sim \lambda \ \ \mbox{ie} \ \ \lambda(t)\sim |t|.$$

{\bf step 4} Uniqueness.\\

Once the blow up speed is known, one can use the energy machinery \fref{firstmoneot} again to compare the flow with the constructed solution $u_c$ and derive a bound 
\be
\label{vnoihohoh}
\|u-u_c\|_{H^1}\lesssim |t|^A
\ee
for some constant $A$. Uniqueness then follow by integrating the flow close to $u_c$ backwards from infinity using a dynamical system type of argument and the a priori bound \fref{vnoihohoh}. The difficulty is to beat the polynomial growth of the propagator $e^{itL}$ generated by  the null space of the -non self adjoint- operator $L$ in a context where all possible gains are also only polynomial\footnote{and not exponential as in \cite{Marteluniqueness}, \cite{DM1} where a similar issue arises}. Recall indeed that the group of symmetries of the homogeneous case generates polynomially growing solutions to the linearized flow, explicitly:
\be
\label{structurekernelimpair}
L_+(\nabla Q)=0, \ \  L_-(yQ)=-2\nabla Q,
\ee
\be
\label{structurekernelpair}
L_-Q=0,\ \ L_+(\Lambda Q)=-2Q, \ \ L_-(|y|^2Q)=-4\Lambda Q, \ \ L_+\rho=|y|^2Q, 
\ee
where $\rho$ is the unique radial $H^1$ solution to $L_+\rho=|y|^2Q$. Beating the induced $t^3$ a priori growth of the propagator $e^{itL}$ requires a sufficiently large a priori constant $A=3$ in \fref{vnoihohoh} which in turn forces the construction of an approximate solution to a sufficiently high order.\\

The rest of the paper is as follows. In section \ref{approx}, we construct an approximate solution to the order 5 in renormalized variables, Proposition \ref{propconstution}. In section \ref{sectionexistence}, we derive a mass critical energy identity, Lemma \ref{lemma:timederivative}, which is enough to prove the existence of minimal mass blow-up solutions, Proposition \ref{existenceprop}. In section \ref{sectionspeed}, we start the proof of uniqueness. We derive the local virial identity, Proposition \ref{prop:localvirial}, which is the key to the control of the blow up speed and dispersion, Proposition \ref{prop:improvecontrol}. In section \ref{uniqueness}, we iterate the use of the energy identity of Lemma \ref{lemma:timederivative} to control the growth induced by the null space, Lemma \ref{scalarproductterms}, and conclude the proof of uniqueness.\\

{\bf Notations} We let $$(f,g)=\int_{\RR^2} f(x)g(x)dx$$ be the $L^2$ scalar product. We introduce the differential operator $$\Lambda f=f+y\cdot\nabla f$$ and recall that $$(f,\Lambda g)=-(\Lambda f,g).$$


\section{Approximate solutions and leading order dynamics}
\label{approx}

From now on and for the rest of this paper, we assume that $k$ satisfies {\bf Assumption (H)}. Moreover, without loss of generality, we assume that $k$ attains its maximum at $x_0=0$ which is nondegenerate:
 \be
\label{nondgeen}
 k(0)=1,  \ \ \nabla k(0)=0, \ \ \nabla^2k(0)<0.
\ee Our aim in this section is to construct a high order approximate critical mass blow up solution by following the slow variables method implemented in particular in \cite{KMR}. The outcome is the derivation of the leading order dynamical system \fref{odemain} driving the geometrical parameters attached to the refined modulation theory. In particular, we will see how an {\it algebraic cancellation} leads to the preservation of the pseudoconformal speed, and a new rigidity occurs for the motion of the center of mass.


\subsection{The slow modulated ansatz}


Let us consider the general modulated ansatz:
$$u(t,x)=\frac{1}{[k(\alpha(t))]^{\frac{1}{2}}}\frac{1}{\lambda(t)}v\left(s,\frac{x-\alpha(t)}{\lambda(t)}\right)e^{i\gamma(t)}, \ \ \frac{ds}{dt}=\frac{1}{\lambda^2}$$ which maps the finite time blow up problem \fref{nlsk:0} onto the global in time renormalized equation \fref{th1:eqrenormlaizedv}:
\bea
\label{th1:eqrenormlaizedv}
\nonumber & & i\partial_s v+\Delta v-v+\frac{k(\lambda(t)y+\alpha(t))}{k(\alpha(t))}v|v|^{2}\\
& = & i\lsl \Lambda v+i\xsl \cdot\left(\nabla v+\frac{\lambda}{2}\frac{\nabla k(\alpha(t))}{k(\alpha(t))}v\right)+\gts v
\eea
with $\gts=\gamma_s-1.$ We freeze the homogeneous laws $$\lsl=-b, \ \ \xsl=2\beta, \ \ \gts=|\beta|^2
$$
and let
\be
\label{beonceoeo}
w(s,y)=v(s,y)e^{i\frac{b|y|^2}{4}-i\beta\cdot y},
\ee
so that \fref{th1:eqrenormlaizedv} becomes:
\begin{equation}\label{eqw}
\begin{array}{r}
\ds i\partial_s w+\Delta w-w+(b_s+b^2)\frac{|y|^2w}{4}-\left\{(\beta_s+b\beta)\cdot y+i\lambda\beta\cdot \frac{\nabla k(\alpha(t))}{k(\alpha(t))}\right\}w\\
\ds +\frac{k(\lambda(t)y+\alpha(t))}{k(\alpha(t))}w|w|^{2}=0.
\end{array}
\end{equation}
Note that $Q$ itself provides an approximate solution of order $b^2$ but this is not enough to carry out the analysis and to exhibit fine effects induced by the homogeneity $k$. In order to construct a higher order approximation, we adapt the strategy designed in \cite{KMR} in a subcritical setting to the critical case. Given small parameters $$\P=(b,\l,\b,\a),$$ we look for a solution to \fref{eqw} of the form $$w(s,y)=P_{\P(s)}(y), \ \  \lsl=-b, \ \ \xsl=2\beta, \ \ \gts=|\beta|^2, \ \ b_s+b^2=0, \ \ \beta_s+b\beta=\mathcal{B}(\l,\a).$$ This maps  \fref{eqw} onto:
\bea
\label{eqp}
\nonumber & &  -ib^2\partial_b\pp-i\l b\p_\l\pp+2i\b\l\p_\a\pp+i(-b\b +\mathcal{B})\p_\b\pp\\
\nonumber &  -& \left\{\mathcal{B}\cdot y+i\lambda\beta\cdot \frac{\nabla k(\alpha(t))}{k(\alpha(t))}\right\}\pp+  \Delta \pp-\pp+ \frac{k(\lambda(t)y+\alpha(t))}{k(\alpha(t))}\pp|\pp|^{2} =0.
\eea
The profile $P_{\P}$ is computed from an asymptotic expansion near $Q$ and at each step the laws for the modulation parameters $\mathcal{B}$ which contain the deformation induced by the inhomogeneity $k$  are adjusted to ensure the solvability of the obtained equations. Indeed, recall that the kernel of the linearized operator close to $Q$ is explicit: 
\be
\label{kerl}
Ker\{L_+\}=\mbox{span}\{\nabla Q\}, \ \ Ker\{L_-\}=\mbox{span}\{Q\},
\ee 
see \cite{W1}, \cite{nakspectral}. The 2 instability directions generated by the kernel of $L_+$ are precisely adjusted by computing the two parameters $\mathcal{B}\in \RR^2$. The instability direction generated by the kernel of $L_-$ should in principle be adjusted by modulating on the law of $b$. However, a spectacular cancelation (see \eqref{explicitcomp}) will take care of the kernel of $L_-$, which in turn will allow us to keep the same law for $b$ as in the pseudoconformal regime, i.e. $b_s+b^2=0$.
 
\begin{proposition}[Approximate solution]
\label{propconstution}
Let $C_0>0$ be a given constant. There exists a universal constant $c>0$ and a small constant $\eta^*(C_0)>0$ such that for all $|\P=(b,\l,\b,\a)|\leq \eta^*$, the following holds true. There exist  $\mathcal{B}$ of the form: 
\be
\label{loibone}
\mathcal{B}=\lambda c_0(\alpha)+\b_3\l^3+\b_4\l^4,
\ee
where $c_0$ is a linear map on $\RR^2$ and $\b_3, \b_4$ are vectors in $\RR^2$, and smooth well localized profiles $(T_j,S_j)_{1\leq j\leq 4}$ homogeneous of degree $j$ in $\P$, such that
\be
\label{solapp}
\pp=Q+\sum_{j=1}^4(T_j+iS_j)
\ee
is an approximate solution to \fref{eqp} in the sense:
\bea
\label{eqpbis}
\nonumber & &  -ib^2\partial_b\pp-i\l b\p_\l\pp+2i\b\l\p_\a\pp+i(-b\b +\mathcal{B})\p_\b\pp \\
\nonumber & - & \left\{\mathcal{B}\cdot y+i\lambda\beta\cdot \frac{\nabla k(\alpha(t))}{k(\alpha(t))}\right\}\pp+ \Delta \pp-\pp\\
 & + &\frac{k(\lambda(t)y+\alpha(t))}{k(\alpha(t))}\pp|\pp|^{2} =-\tilde{\Psi}_{\mathcal P}(y)
\eea
with $\tilde{\Psi}_{\mathcal P}(y)$ polynomial in $\mathcal P$ and smooth and well localized in $y$: $\forall \alpha\in \Bbb N^2$,
\be
\label{estpsitilde}
\left|\pa^{\alpha}\tilde{\Psi}_{\mathcal P}(y)\right|\lesssim \left(\P^5+\P(\a^2+\b^2)+\left(b-\frac{\l}{C_0}\right)\P^3\right) e^{-C_{\alpha}|y|}.
\ee
\end{proposition}

\begin{remark}\label{rmk:solapp}
Remark that the remainder \fref{estpsitilde} is not of order five with respect to all parameters. This simplifies the construction and is not a problem since we will obtain better decay properties for $\a$, $\b$ and $b-\l/C_0$. In particular, we will show $\a=O(\P^2)$, $\b=O(\P^2)$ and $b-\l/C_0=O(\P^3)$ where $C_0$ will be explicitly chosen according to \fref{constantcnot}\footnote{see Remark \ref{ie:rmk}}.
\end{remark}

{\bf Proof of Proposition \ref{propconstution}}\\

The proof proceeds by injecting the expansion \fref{solapp} in \eqref{eqpbis}, identifying the terms with the same homogeneity, and inverting the corresponding operator. Let us recall that if $L=(L_+,L_-)$ is the matrix linearized operator close to $Q$ given by \fref{deflpluslmoins}, then the explicit description of the kernel \fref{kerl} ensures the standard uniform elliptic estimates:
\be
\label{ellitpicestinamtes}
\forall f\in (\nabla Q)^{\perp}, \ \ \|e^{c|y|}L_+^{-1}f\|_{H^2}\lesssim \|e^{2c|y|}f\|_{L^2},
\ee
\be
\label{ellitpicestinamtesbis}
\forall f\in (Q)^{\perp}, \ \ \|e^{c|y|}L_-^{-1}f\|_{H^2}\lesssim \|e^{2c|y|}f\|_{L^2}.
\ee

{\bf step 1} Derivation of the law for the center of mass.\\

We inject \fref{solapp} in \eqref{eqpbis} and sort the terms of same homogeneity.\\

{\bf order 1}: From $\nabla k(0)=0$, we obtain
$$L_+(T_1)=0,\, L_-(S_1)=0\ \ \mbox{and hence} \ \ T_1\equiv 0,\,S_1\equiv 0.
$$

{\bf order 2}: We identify the terms homogeneous of order 2 in \eqref{eqp} and get:
\be\label{eqT2S2}
\left\{\begin{array}{l}
\ds L_+(T_2)=\nabla^2 k(0)(\a,y)\l Q^3+\frac{\l^2}{2}\nabla^2k(0)(y,y)Q^3 -\l c_0(\a)\c y Q,\\
\ds L_-(S_2)=0.
\end{array}\right.
\ee
For $S_2$, we may choose:
$$S_2\equiv 0.
$$
In view of the first equation in \eqref{eqT2S2} and \fref{ellitpicestinamtes}, we may solve for $T_2$ if we adjust:
$$
\left(\nabla^2k(0)(\a,y)\l Q^3+\frac{\l^2}{2}\nabla^2k(0)(y,y)Q^3-\l c_0(\a)\c y Q,\n Q\right)=0,
$$
which computes the value of $c_0$:
\be\label{choicec0}
(c_0(\a))_j:=\frac{\left(\int Q^4\right)}{2\left(\int Q^2\right)}\nabla^2k(0)(e_j,\a)\textrm{ for }j=1, 2.
\ee

{\bf step 2} Derivation of pseudo-conformal speed: an algebraic cancellation.\\

Recall that in the pseudoconformal regime $\lambda(t)\sim t$, there holds $$b=-\frac{\lambda_s}{\lambda}=-\lambda\lambda_t\sim \lambda.$$ Equivalently, the law for b is in this case $$b_s=-b^2.$$ Hence the term $\frac{\l^2}{2}\nabla^2k(0)(y,y)Q^3$ in the RHS of the $T_2$ equation \fref{eqT2S2} may potentially induce a drastic modification of the b law of the form: $$b_s=-b^2+C\lambda^2+\dots$$ This explains why the non degenerate case $\nabla^2k$ is a much larger deformation than the degenerate case $\nabla^2k(0)=0$. We now claim that in fact the choice $C=0$ is dictated by the algebraic cancellation \fref{explicitcomp}. To see this, we compute the system at the next order:\\

{\bf order 3}: We identify the terms homogeneous of order 3 in \eqref{eqp}. We get rid of the terms $\P\a^2$, $\P\a\b$ and $\P\b^2$ which are remainder terms. We obtain:
\be\label{eqT3S3}
\left\{\begin{array}{ll}
\ds L_+(T_3)= & \ds \nabla^3k(0)(y,y,y)\frac{\l^3}{6} Q^3+\nabla^3k(0)(y,y,\a)\frac{\l^2}{2}Q^3-\b_3\l^3\cdot yQ,\\
\ds L_-(S_3)= & \ds -\l b\p_\l T_2+2\b\l\p_\a T_2.
\end{array}\right.
\ee
In view of the first equation in \eqref{eqT3S3}, we may solve for $T_3$ if and only if:
$$\ds\bigg(\nabla^3k(0)(y,y,y)\frac{\l^3}{6} Q^3+\nabla^3k(0)(y,y,\a)\frac{\l^2}{2}Q^3 -\b_3\l^3\cdot yQ,\n Q\bigg)\ds =0,
$$
which is true provided we choose $\b_3$ to be:
$$(\b_3)_j=\frac{\left(\int \nabla^3k(0)(y,y,e_j)Q^4\right)}{4\left(\int Q^2\right)}\textrm{ for }j=1, 2.
$$

In view of the second equation in \eqref{eqT3S3}, we may solve for $S_3$ if:
\be\label{S3:1}
\left(-\l b\p_\l T_2+2\b\l\p_\a T_2,Q\right)=0.
\ee
We now integrate by parts using $L_+(\L Q)=-2Q$ and \fref{eqT2S2} to compute:
\be
\label{cbcoceoicoio}
\nonumber -2(T_2,Q)  =  (L_+T_2,\Lambda Q)=\frac{\lambda^2}{2}(\nabla^2k(0)(y,y)Q^3,\Lambda Q)= 0,
\ee
which shows that \eqref{S3:1} does actually hold. Here we used the {\it spectacular algebraic cancellation}: 
\be\label{explicitcomp}
\left(y_jy_lQ^3,\L Q\right)=0\textrm{ for }j, l=1,2.
\ee
{\it This cancellation is the reason why we did not need to adjust the law for $b$ at the $\lambda^2$ level and hence why the critical blow up will still be of conformal type}.\\

{\bf step 3} Higher order corrections.\\

At this stage, we have computed the leading order terms in the derivation of the modified modulation equations. We however keep expanding the profile to get higher order corrections, that will be required for the final step of the proof of uniqueness, see section \ref{uniqueness}.\\

{\bf order 4}: We identify the terms homogeneous of order 4 in \eqref{eqp}. We get rid of the terms $\P^3\a$ and $\P^3\b$ which are remainder terms. Since terms of type $(b-\l/C_0)\P^3$ are also remainder terms, we may replace everywhere $b$ by $\l/C_0$. We obtain:
\be\label{eqT4S4}
\left\{\begin{array}{ll}
\ds L_+(T_4)= & \ds - \b_4\l^4\cdot yQ+f_4\l^4,\\
\ds L_-(S_4)= & \ds -\frac{\l^2}{C_0}\p_\l T_3,
\end{array}\right.
\ee
where $f_4$ depends on $T_2, T_3, S_3$ and is a function from $\RR^2$ to $\RR$. In view of the first equation in \eqref{eqT4S4}, we may solve for $T_4$ if and only if we have:
$$\bigg(- \b_4\cdot yQ+f_4,\n Q\bigg)  =0,
$$
which is true provided we choose $\b_4$ to be:
$$(\b_4)_j=-\frac{2\left(\int f_4\p_j Q\right)}{\left(\int Q^2\right)}\textrm{ for }j=1, 2.
$$
In view of the second equation in \eqref{eqT4S4}, we may solve for $S_4$ if and only if we have:
\be\label{S4:1}
\left(-\frac{\l^2}{C_0}\p_\l T_3,Q\right)=0.
\ee
Since $L_+(\L Q)=-2Q$, this is equivalent to:
\be\label{S4:2}
\frac{\l^2}{2C_0}\left(\p_\l L_+(T_3),\L Q\right)=0.
\ee
Now, in view of the first equation of \eqref{eqT3S3} and the cancellation \fref{explicitcomp}, we have:
\be\label{S4:2bis}
\left(L_+(T_3),\L Q\right)=0.
\ee
In view of \eqref{S4:2} and \eqref{S4:2bis}, \eqref{S4:1} holds. This concludes the construction of $S_4$.\\

{\bf step 4} Proof of \fref{estpsitilde}.\\

The estimate \fref{estpsitilde} now follows from the above construction, the uniform exponential bounds \fref{ellitpicestinamtes}, \fref{ellitpicestinamtesbis} and the explicit polynomial development of the nonlinear term $u|u|^2$. The details are left to the reader.\\
This concludes the proof of Proposition \ref{propconstution}.


\subsection{Properties of the approximate profile}


Let $P_{\P}$ be the refined pseudo-conformal profile given by Proposition \ref{propconstution} which has been computed for simplicity in the conformal variables \fref{beonceoeo}. Going back to the self similar variables, we let: 
\be\label{solapp2}
\qp=\pp e^{-ib\frac{|y|^2}{4}+i\b\c y},
\ee
so that  from direct check:
\bea
\label{eqq}
 \nonumber & - & ib^2\partial_b\qp-i\l b\p_\l\qp+i(-b\b +\mathcal{B})\p_\b\qp+2i\b\l\p_\a\qp-i\lambda\beta\cdot \frac{\nabla k(\alpha(t))}{k(\alpha(t))}\qp\\
\nonumber & +& \Delta \qp-\qp+\frac{k(\lambda(t)y+\alpha(t))}{k(\alpha(t))}\qp|\qp|^{2}+ib\L\qp-2i\b\n\qp-|\b|^2\qp\\
& = & - \Psi_{\P}
\eea
where 
\be
\label{noiveauerror}
\Psi_{\P}=\tilde{\Psi}_{\P}e^{-ib\frac{|y|^2}{4}+i\b\c y}
\ee
 satisfies from \fref{estpsitilde}:
\be
\label{estpsitildebis}
\left\|e^{c|y|}\Psi_{\mathcal P}\right\|_{H^4}\lesssim \left(\P^5+\P(\a^2+\b^2)+\left(b-\frac{\l}{C_0}\right)\P^3\right).
\ee

Let us compute the $L^2$ norm and energy of $\qp$ which will appear as important quantities in the analysis.

\begin{lemma}[Invariants of $\qp$]
\label{lemmainvariants}
There holds:
\be\label{compm4}
\int |\qp|^2=\int Q^2+O(\P^4),
\ee
\bea
\label{compe4}
  \tilde{E}(\qp) &  =& \frac 12 \int|\nabla \qp|^2-\frac{1}{4}\int\frac{k(\lambda y+\alpha)}{k(\alpha)}|\qp|^4\\
\nonumber& = &   \frac{b^2}{8}\int |y|^2Q^2+\frac{|\b|^2}{2}\int Q^2 -\frac{\l^2}{8}\int \nabla^2k(0)(y,y)Q^4+O(\P^4).
\eea
\end{lemma}

{\bf Proof of Lemma \ref{lemmainvariants}}\\

We compute the mass of the approximate solution up to fourth order. Using \eqref{solapp2}, we have:
$$\int |\qp|^2 = \int (Q+T_2+T_3)^2+O(\P^4)=\int Q^2 +2(T_2,Q)+2(T_3,Q)+O(\P^4).
$$
From \fref{cbcoceoicoio}: 
\be
\label{calcultdeuxq}
(T_2,Q)=0,
\ee 
while from \fref{S4:2bis} and $L_+(\Lambda Q)=-2Q$, 
\be
\label{cajcojfor}
(T_3,Q)=0,
\ee and \fref{compm4} follows.\\
For the energy, we have from \eqref{solapp2}:
\bea
 \label{compe1}
 \tilde{E}(\qp) & = & \ds\frac{1}{2}\int |\n\qp|^2 -\frac{1}{4}\int \frac{k(\l y+\a)}{k(\a)}|\qp|^4\\
\nonumber &  = & \ds\frac{1}{2}\int \left|\n Q+\n T_2+\n T_3+\n T_4+\left(-\b+\frac{by}{2}\right)(S_3+S_4)\right|^2\\
\nonumber& + & \frac{1}{2}\int \left|\left(\b-\frac{by}{2}\right)(Q+T_2+T_3+T_4)+\n S_3+\n S_4\right|^2\\
\nonumber& - & \frac{1}{4}\int \frac{k(\l y+\a)}{k(\a)}((Q+T_2+T_3+T_4)^2+(S_3+S_4)^2)^2\\
\nonumber& = & \frac{1}{2}\int |\n Q|^2 -\frac{1}{4}\int Q^4 + \frac{b^2}{8}\int |y|^2Q^2+\frac{|\b|^2}{2}\int Q^2
+(T_2,-\D Q-Q^3)\\
\nonumber& + &(T_3,-\D Q-Q^3)-\frac{1}{4}\int \left(\frac{k(\l y+\a)}{k(\a)}-1\right)Q^4+O(\P^4).
\eea
Injecting the soliton equation $-\Delta Q-Q^3=-Q$ and \fref{calcultdeuxq}, \fref{cajcojfor} into \fref{compe1}, and finally expanding the inhomogenity $k(\l y+\a)$ near $0$ using the radiallity of $Q$ yields \fref{compe4}.\\
This concludes the proof of Lemma \ref{lemmainvariants}.


\section{Energy estimates and existence of critical elements}
\label{sectionexistence}

Our aim in this section is to provide a robust framework to prove the existence of critical blow up elements in the absence of explicit symmetries. Once the leading order behavior is exhibited through the construction of the approximate profile with a small enough correction, our strategy is to integrate the flow backwards from the singularity. Similar strategies are well known in scattering theory and it was used in the (NLS) framework in \cite{Merlemulti}, \cite{BW}, \cite{Marteluniqueness}, \cite{MMmulti}, \cite{KMR}. The key is energy estimates as recalled in the strategy of the proof, with this additional subtlety in the critical setting\footnote{which was in particular absent in \cite{KMR}} that one needs to add a Morawetz type information to derive a suitable Lyapounov type monotonicity information, see \fref{defI}.


\subsection{Nonlinear decomposition of the wave and modulation equations}
\label{refinedenergybounds}


Let $u(t)\in H^1$ be a solution to \fref{nlsk:0} on a time interval $[t_0,t_1]$, $t_1<0$. We assume that on this time interval, the solution is in what will be proved to be the asymptotic regime for critical mass blow up solutions, namely we assume that  $u(t)$ admits a geometrical decomposition 
\be
\label{newdecomp}
u(t,x)=\frac{1}{[k(\alpha(t))]^{\frac{1}{2}}}\frac{1}{\lambda(t)}(Q_{\P(t)}+\e)\left(t,\frac{x-\alpha(t)}{\lambda(t)}\right)e^{i\gamma(t)}
\ee
with a uniform smallness bound on $[t_0,t_1]$:
\be
\label{uniformmallinit}
|\mathcal \P(t)|+\|\e(t)\|_{H^1}\lesssim \lambda(t)\ll 1.
\ee
Moreover, we assume that $u(t)$ has almost critical mass in the sense: $\forall t\in [t_0,t_1]$,
\be
\label{almostcriticalmass}
\left|\|u\|^2_{L^2}-\|Q\|^2_{L^2}\right|\lesssim \lambda^4(t).
\ee
From standard modulation argument, see e.g. \cite{MR1}, the uniqueness of the nonlinear decomposition \fref{newdecomp} may be ensured by imposing a suitable set of orthogonality conditions on $\e$, namely: 
\be\label{ortho1}
(\e_2,\n\s)-(\e_1,\n\t)=0,
\ee
\be\label{ortho2}
(\e_1,y\s)+(\e_2,y\t)=0,
\ee
\be\label{ortho3}
-(\e_1,\L\t)+(\e_2,\L\s)=0,
\ee
\be\label{ortho4}
(\e_1,|y|^2\s)+(\e_2,|y|^2\t)=0,
\ee
\be\label{ortho5}
-(\e_1,\rho_2)+(\e_2,\rho_1)=0,   
\ee
where $\rho$ is the unique even $H^1$ solution to $L_+\rho=|y|^2Q$ and $$\rho_1+i\rho_2=\rho(y)e^{-ib\frac{|y|^2}{4}+i\b\cdot y}.$$ Recall the non degeneracy:
\be
\label{nondegenrho}
(\rho,Q)=-\frac{1}{2}(L_+\rho,\Lambda Q)=-\frac{1}{2}(|y|^2Q,\Lambda Q)=\frac{1}{2}|yQ|_{L^2}^2.
\ee 
These orthogonality conditions correspond exactly in the case $\P=(0,0,0,0)$ to the null space of the linearized operator close to $Q$, see \fref{structurekernelimpair}, \fref{structurekernelpair}, and the directions involved in \fref{ortho1}-\fref{ortho5} provide a first approximation of the null space close to $\qp$\footnote{A slight refinement will be needed, see \fref{diff8}, \fref{diff9}}. From standard argument, the obtained modulation parameters are $\mathcal C^1$ functions of time, see \cite{MR1} for related statements. Let 
\be
\label{rescaledtime}
s(t)=\int_{t_0}^{t_1}\frac{d\tau}{\lambda^2(\tau)}
\ee be the rescaled time\footnote{which will be a global time, $s([t_0,0))=[s(t_0),+\infty)$}, then $\e$ satisfies for $s\in [s_0,s_1]$ the equation:
\bea
\label{eqe1}
& &  (b_s+b^2)\p_b\s+\l\left(\frac{\l_s}{\l}+b\right)\p_\l\s\\
\nonumber & + & (\b_s+b\b-c_0(\a)\l-\b_3\l^3-\b_4\l^4)\p_\b\s+\l\left(\frac{\a_s}{\l}-2\b\right)\p_\a\s\\
\nonumber &+ &\p_s\e_1-M_2(\e)+b\L\e_1-2\b\left(\n\e_1+\frac{\l}{2}\frac{\n k(\a)}{k(\a)}\e_1\right)-|\b|^2\e_2\\
\nonumber &= &\left(\frac{\l_s}{\l}+b\right)(\L\s+\L\e_1)+(\tgamma_s-|\b|^2)(\t+\e_2)\\
\nonumber &+ &\left(\frac{\a_s}{\l}-2\b\right)\left(\n\s+\n\e_1+\frac{\l}{2}\frac{\n k(\a)}{k(\a)}\s+\frac{\l}{2}\frac{\n k(\a)}{k(\a)}\e_1\right)+\Im(\Psi_{\P})-R_2(\e),
\eea
and 
\bea
\label{eqe2}
& &  (b_s+b^2)\p_b\t+\l\left(\frac{\l_s}{\l}+b\right)\p_\l\t\\
\nonumber & + &(\b_s+b\b-c_0(\a)\l-\b_3\l^3-\b_4\l^4)\p_\b\t+\l\left(\frac{\a_s}{\l}-2\b\right)\p_\a\t\\
\nonumber &+ &\p_s\e_2+M_1(\e)+b\L\e_2-2\b\left(\n\e_2+\frac{\l}{2}\frac{\n k(\a)}{k(\a)}\e_2\right)+|\b|^2\e_1\\
\nonumber & = &\left(\frac{\l_s}{\l}+b\right)(\L\t+\L\e_2)-(\tgamma_s-|\b|^2)(\s+\e_1)\\
\nonumber & + &\left(\frac{\a_s}{\l}-2\b\right)\left(\n\s+\n\e_2+\frac{\l}{2}\frac{\n k(\a)}{k(\a)}\s+\frac{\l}{2}\frac{\n k(\a)}{k(\a)}\e_2\right) -\Re(\Psi_{\P})+R_1(\e),
\eea
where $M_1, M_2$ are small deformations of the linearized  operator $(L_+,L_-)$ close to $Q$:
\be\label{defm1}
\ds M_1(\e)= -\D\e_1+\e_1-\frac{k(\l y+\a)}{k(\a)}((|\qp|^2+2\s^2)\e_1+2\s\t\e_2),
\ee
\be\label{defm2}
\ds M_2(\e)= -\D\e_2+\e_2-\frac{k(\l y+\a)}{k(\a)}((|\qp|^2+2\t^2)\e_2+2\s\t\e_1),
\ee
where the nonlinear terms are given by 
\be\label{defr1}
\ds R_1(\e)=3\s\e_1^2+2\t\e_1\e_2+\s\e_2^2+|\e|^2\e_1,
\ee
\be\label{defr2}
\ds R_2(\e) =3\t\e_2^2+2\s\e_1\e_2+\t\e_1^2+|\e|^2\e_2,
\ee
and where  $\Psi_{\P}$ given by \fref{noiveauerror} denotes the remainder term in the equation of $\qp$ and satisfies \fref{estpsitildebis}.\\
Let us collect the standard preliminary estimates on this decomposition which rely on the conservation laws and the explicit choice of orthogonality conditions.

\begin{lemma}[Preliminary estimates on the decomposition]
\label{estprelimestdecomp}
There holds the bounds for $s\in [s_0,s_1]$:\\
{\em 1. Energy bound}: 
\be
\label{cme7}
b^2+|\b|^2+|\a|^2+\|\e\|_{H^1(\RR^2)}^2\lesssim \l^2\left(E_0+\frac{1}{8}\int \nabla^2k(0)(y,y)Q^4\right)+O(\P^4).
\ee
{\em 2. Control of the geometrical parameters}: Let the quadratic forms:
\be
\label{defdzeroalpha}
d_0(\alpha,\alpha)=\frac{2\|Q\|_{L^2}^2}{\|yQ\|_{L^2}^2}\nabla^2k(0)(\alpha,\alpha),\ \ d_1(\alpha,\alpha)=\frac{(|y|^2Q,\rho)}{4(\rho,Q)}d_0(\alpha,\alpha),
\ee
and the vector of modulation equations
\bea
\label{defmod}
\nonumber  Mod(t) & = & \Big(b_s+b^2-d_0(\alpha,\alpha), \tgamma_s-|\beta|^2+d_1(\alpha,\alpha),\xsl-2\beta,\lsl+b,\\
& & \beta_s+b\beta-c_0(\alpha)\lambda-\b_3\l^3\Big)
\eea
then the modulation equations are to leading order:
\be
\label{law2}
|Mod(t)| \lesssim  \P^4+\P^2\|\e\|_{L^2}+\|\e\|^2_{L^2}+\|\e\|_{H^1}^3+\left(b-\frac{\l}{C_0}\right)\P^3+\P(\a^2+\b^2),
\ee
with the improvement:
\be
\label{law2:bis}
\left|\frac{\l_s}{\l}+b\right| \lesssim  \P^5+\P^2\|\e\|_{L^2}+\|\e\|^2_{L^2}+\|\e\|_{H^1}^3+\left(b-\frac{\l}{C_0}\right)\P^3+\P(\a^2+\b^2).
\ee
\end{lemma}

{\bf Proof of Lemma \ref{estprelimestdecomp}}\\

{\bf step 1} Conservation of $L^2$ norm and energy.\\

Let us write down the conservation of $L^2$ norm using \fref{newdecomp}, \fref{almostcriticalmass}:
\be
\nonumber \int |\qp+\e|^2dy  =  k(\a)\int|u|^2=k(\alpha)\int Q^2+k(\alpha)\left[\int|u|^2-\int Q^2\right]
\ee
from which we obtain:
\be
\label{degeneracyeq}
2\Re(\e,\overline{\qp})+\int |\e|^2=k(\alpha)\left[\int|u|^2-\int Q^2\right]-\left[\int|\qp|^2-k(\alpha)\int Q^2\right].
\ee
We now use \eqref{compm4}, \fref{almostcriticalmass} to conclude:
\be\label{cm5}
-\frac{\nabla^2k(0)(\a,\a)}{2}\left(\int Q^2\right)+2\Re\left(\int \e\overline{\qp}\right)+\int |\e|^2=O(\P^4+\P|\alpha|^2).
\ee
Consider now the conservation of energy:
\be\label{ce3}
\frac{1}{2}\int |\n\qp+\n\e|^2 -\frac{1}{4}\int \frac{k(\l y+\a)}{k(\a)}|\qp+\e|^4=k(\a)\l^2E_0.
\ee
We expand the nonlinear term:
\bee
 |\qp+\e|^4 & = & \ds |\qp|^4+4\Re(\e\overline{|\qp|^2\qp})\\
&  + &2|\qp|^2\left((1+\frac{2\s^2}{|\qp|^2})\e_1^2+4\frac{\s\t}{|\qp|^2}\e_1\e_2+(1+\frac{2\t^2}{|\qp|^2})\e_2^2\right)\\
&  + &4\Re(|\e|^2\e\overline{\qp})+|\e|^4.
\eee
We now inject the value of $\tilde{E}(\qp)$ given by \eqref{compe4} and estimate the cubic and higher nonlinear terms using standard Gagliardo-Nirenberg estimates and the a priori smallness \fref{uniformmallinit} to derive:
\bee
& &\frac{1}{2}\int |\n\qp+\n\e|^2 -\frac{1}{4}\int \frac{k(\l y+\a)}{k(\a)}|\qp+\e|^4\\
\nonumber &= &\tilde{E}(\qp)+\Re\left(\e,-\overline{\D\qp-\frac{k(\l y+\a)}{k(\a)}|\qp|^2\qp}\right)\\
\nonumber& + & \frac{1}{2}\int |\n\e|^2-\frac{1}{2}\int |\qp|^2\left(\left(1+\frac{2\s^2}{|\qp|^2}\right)\e_1^2+4\frac{\s\t}{|\qp|^2}\e_1\e_2+\left(1+\frac{2\t^2}{|\qp|^2}\right)\e_2^2\right)\\
\nonumber& + & O(\|\e\|_{H^1(\RR^2)}^3+\|\e\|_{H^1(\RR^2)}^2\P^2)\\
\nonumber& = &\ds\frac{b^2}{8}\int |y|^2Q^2+\frac{|\b|^2}{2}\int Q^2-\frac{\l^2}{8}\int \nabla^2k(0)(y,y)Q^4\\
\nonumber&+ & \Re\left(\e,\overline{-\D\qp-\frac{k(\l y+\a)}{k(\a)}|\qp|^2\qp}\right)+O(\|\e\|_{H^1(\RR^2)}^3+\P^4)\\
\nonumber& + & \frac{1}{2}\int |\n\e|^2-\frac{1}{2}\int |\qp|^2\left(\left(1+\frac{2\s^2}{|\qp|^2}\right)\e_1^2+4\frac{\s\t}{|\qp|^2}\e_1\e_2+\left(1+\frac{2\t^2}{|\qp|^2}\right)\e_2^2\right),
\eee
and hence using \eqref{ce3}:
\bea
\label{ce6}
& & \l^2E_0= \ds\frac{b^2}{8}\int |y|^2Q^2+\frac{|\b|^2}{2}\int Q^2-\frac{\l^2}{8}\int \nabla^2k(0)(y,y)Q^4\\
\nonumber & + &\Re\left(\e,\overline{-\D\qp-\frac{k(\l y+\a)}{k(\a)}|\qp|^2\qp}\right)+O(\|\e\|_{H^1(\RR^2)}^3+\|\e\|_{H^1(\RR^2)}^2\P^2+\P^4)\\
\nonumber & + & +\frac{1}{2}\int |\n\e|^2-\frac{1}{2}\int |\qp|^2\left(\left(1+\frac{2\s^2}{|\qp|^2}\right)\e_1^2+4\frac{\s\t}{|\qp|^2}\e_1\e_2+\left(1+\frac{2\t^2}{|\qp|^2}\right)\e_2^2\right).
\eea

{\bf step 2} Coercivity of the linearized energy and proof of \fref{cme7}.\\

We now sum the conservation of mass \eqref{cm5} and the conservation of energy \eqref{ce6}. We obtain:
\be\label{cme1}
\begin{array}{l}
 \ds\frac{b^2}{8}\int |y|^2Q^2+\frac{|\b|^2}{2}\int Q^2-\frac{\nabla^2k(0)(\a,\a)}{4}\int Q^2\\
\ds +\Re\left(\e,\overline{\qp-\D\qp-\frac{k(\l y+\a)}{k(\a)}|\qp|^2\qp}\right)+\frac{1}{2}\int |\e|^2+\frac{1}{2}\int |\n\e|^2\\\ds -\frac{1}{2}\int |\qp|^2\left(\left(1+\frac{2\s^2}{|\qp|^2}\right)\e_1^2+4\frac{\s\t}{|\qp|^2}\e_1\e_2+\left(1+\frac{2\t^2}{|\qp|^2}\right)\e_2^2\right)\\
\ds =\ds \l^2\left(E_0+\frac{1}{8}\int \nabla^2k(0)(y,y)Q^4\right)+O(\|\e\|_{H^1}^3+|\alpha|^3+\P^4).
\end{array}
\ee
The remaining linear term is degenerate from \eqref{eqq}:
$$ \qp-\Delta \qp-\frac{k(\lambda(t)y+\alpha(t))}{k(\alpha(t))}\qp|\qp|^{2}=ib\L\qp-2i\b\n\qp+O(\P^2),
$$and thus using also the orthogonality conditions \eqref{ortho1} and \eqref{ortho3}:
\bea
\label{cme4}
\nonumber \Re\left(\e,\overline{\qp-\D\qp-\frac{k(\l y+\a)}{k(\a)}|\qp|^2\qp}\right)& = & b\Im(\e,\overline{\L\qp})-2\b\Im(\e,\overline{\n\qp})+O(\P^2\|\e\|_{L^2})\\
& = & O(\P^2\|\e\|_{L^2}).
\eea
\eqref{cme1} and \eqref{cme4} imply:
\bea
\label{cme5}
\nonumber & & \frac{b^2}{8}\int |y|^2Q^2+\frac{|\b|^2}{2}\int Q^2-\frac{\nabla^2k(0)(\a,\a)}{4}\int Q^2+\frac{1}{2}\int |\e|^2+\frac{1}{2}\int |\n\e|^2\\
& - & \frac{1}{2}\int |\qp|^2\left(\left(1+\frac{2\s^2}{|\qp|^2}\right)\e_1^2+4\frac{\s\t}{|\qp|^2}\e_1\e_2+\left(1+\frac{2\t^2}{|\qp|^2}\right)\e_2^2\right)\\
\nonumber &= & \l^2\left(E_0+\frac{1}{8}\int \nabla^2k(0)(y,y)Q^4\right)+O(\|\e\|_{H^1}^3+|\alpha|^3+\|\e\|_{H^1}\P^2+\P^4).
\eea
We now observe from the proximity of $\qp$ to $Q$ ensured by the a priori smallness \fref{uniformmallinit} that the quadratic form in the LHS of \fref{cme5} is a small deformation of the linearized energy close to $Q$. In other words, we rewrite \fref{cme5}:
\bea
\label{cme5bis}
\nonumber & & \frac{b^2}{8}\int |y|^2Q^2+\frac{|\b|^2}{2}\int Q^2-\frac{\nabla^2k(0)(\a,\a)}{4}\int Q^2+\frac{1}{2}\left[(L_+\e_1,\e_1)+(L_-\e_2,\e_2)\right]\\
& =&  \l^2\left(E_0+\frac{1}{8}\int \nabla^2k(0)(y,y)Q^4\right)+o(\|\e\|_{H^1}^2+|\alpha|^2)+O(\|\e\|_{H^1}\P^2+\P^4).
\eea
We now recall the following coercivity property of the linearized energy which is a well known consequence of the variational characterization of $Q$:

\begin{lemma}[Coercivity of the linearized energy, \cite{W2}, \cite{MR1}, \cite{MR4}]
\label{lemmacoerc}
There holds for some universal constant $c_0>0$ :  $\forall \e\in H^1$, 
\bea
\label{coerclinearenergy}
& & (L_+\e_1,\e_1)+(L_-\e_2,\e_2) \geq c_0  \|\e\|_{H^1}^2\\
\nonumber & - & \frac{1}{c_0}\left\{(\e_1,Q)^2+(\e_1,|y|^2Q)^2+(\e_1,yQ)^2+(\e_2,\rho)^2\right\}. \eea
\end{lemma}

The choice of orthogonality conditions \footnote{The standard orthogonality condition in \fref{coerclinearenergy} on $\e_2$ is $(\e_2,Q)=0$ because $Q$ is the bound state $L_-Q=0$. In \cite{MR1}, it is shown how to replace this by $(\e_2,\Lambda^2Q)=0$, see proof of Lemma 3, and the same proof applies verbatim with the choice $(\e_2,\rho)=0$ using the key non degeneracy \fref{nondegenrho}.} \fref{ortho1}-\eqref{ortho5} together with the degeneracy inherited from \fref{cm5}: 
$$ |(\e_1,Q)|^2=o(\|\e\|_{H^1}^2+|\alpha|^2)+O(\mathcal P^4)$$ yield: $$(L_+\e_1,\e_1)+(L_-\e_2,\e_2)\geq \frac{c_0}{2}\|\e\|_{H^1}^2+o(|\alpha|^2)+O(\mathcal P^4).$$ Injecting this into \fref{cme5bis} and using the nondegeneracy \fref{nondgeen} yields \fref{cme7}.\\

{\bf step 3} Computation of the modulation equations.\\

We now compute the modulation parameters using the sets of orthogonality conditions and the $\e$ equation \fref{eqe1}, \fref{eqe2}. This is done in detail in Appendix A, and the bounds \fref{law2}, \eqref{law2:bis} follow.\\
This concludes the proof of Lemma \ref{estprelimestdecomp}.


\subsection{Refined energy identity}


Our aim in this subsection is to derive a general refined mixed energy/Morawetz type estimate  which will allow us to derive a Lyapounov function for critical mass blow up solutions. We shall work within the following general framework. We let $u$ be a solution to (\ref{nlsk:0}) on $[t_0,0)$ and $w$ be an approximate solution to (\ref{nlsk:0}): 
\be\label{un2}
i\p_tw+\Delta w+k(x)|w|^2w=\psi,
\ee
with the a priori bounds 
\be
\label{aprioriboundwsw}
\|w\|_{L^2}\lesssim 1, \ \ \|\nabla w\|_{L^2}\lesssim \frac{1}{\lambda}, \ \ \|w\|_{H^{\frac{3}{2}}}\lesssim \frac{1}{\lambda^{\frac{3}{2}}}.
\ee
We then decompose $u=w+\tilde{u}$ so that $\tilde{u}$ satisfies: 
\be\label{un3}
i\p_t\tu+\Delta\tu+k(x)(|u|^2u-|w|^2w)=-\psi
\ee
and assume the a priori bounds on $\tu$\footnote{corresponding to the asymptotic regime near the singularity, see \fref{bis:ex17}}:
\be
\label{aprioritu}
\|\nabla \tu\|_{L^2}\lesssim \lambda, \ \ \|\tu\|_{L^2}\lesssim \lambda^2
\ee
and on the geometrical parameters:
\be
\label{aprioirigeometrical}
\ \ |\lambda\lambda_t+b|\lesssim \lambda^4, \ \ b\sim \lambda, \ \ |\lambda\alpha_t|\lesssim \lambda, \ \ |b_t|\lesssim 1
\ee
for some nonnegative parameters $0<\lambda,b\ll 1$ and $\alpha\in \RR^2$.\\ 
We let $A>0$ be a large enough constant which will be chosen later and let $\phi:\RR^2\goto \RR$ be a smooth radially symmetric cut off function with 
\be
\label{defphi}
\phi'(r)=\left\{\begin{array}{ll} r \ \ \mbox{for} \ \ r\leq 1,\\ 3-e^{-r} \ \ \mbox{for}\ \ r\geq 2.\end{array}\right .
\ee
Let $$F(u)=\frac{1}{4}|u|^4, \ \ f(u)=u|u|^2 \ \ \mbox{so that} \ \ F'(u)\cdot h=Re(f(u)\overline{h}).$$ We claim the following generalized energy estimate on the linearized flow \fref{un3}:

\begin{lemma}[Generalized energy estimate]
\label{lemma:timederivative}
Let 
\bea
\label{defI}
\nonumber \mathcal I & = & \frac{1}{2}\int |\n\tu|^2 +\half\int \frac{|\tu|^2}{\l^2}-\int k(x)\left[F(w+\tu)-F(w)-F'(w)\cdot\tilde{u}\right]\\
 & + &  \half\frac{b}{\l}\Im\left(\int A\nabla\phi\left(\frac{x-\a}{A\l}\right)\cdot\nabla\tu\overline{\tu}\right)
\eea
then there holds:
\bea
\label{crc6}
 \frac{d\mathcal I}{dt} & = &-\frac{1}{\l^2}\Im\left(\int k(x)w^2\overline{\tu}^2\right)-\Re\left(\int k(x)w_t\overline{(2|\tu|^2w+\tu^2\overline{w})}\right)\\
\nonumber &+ &\frac{b}{\l^2}\Bigg(\int \frac{|\tu|^2}{\l^2}+\Re\left(\int\nabla^2\phi\left(\frac{x-\a}{A\l}\right)(\n\tu,\overline{\n\tu})\right)-\frac{1}{4A^2}\left(\int\Delta^2\phi\left(\frac{x-\a}{A\l}\right)\frac{|\tu|^2}{\l^2}\right)\\
\nonumber & + & \frac{b}{\l}\Re\left(\int A\nabla\phi\left(\frac{x-\a}{A\l}\right)k(x)(2|\tu|^2w+\tu^2\overline{w})\cdot\overline{\nabla w}\right)\Bigg)\\
\nonumber & + & \Im\left(\int\left[\Delta\psi-\frac{\psi}{\l^2}+k(x)(2|w|^2\psi-w^2\overline{\psi})+i\frac{b}{\l}A\nabla\phi\left(\frac{x-\a}{A\l}\right)\cdot\nabla\psi\right .\right .\\
\nonumber & + & \left .\left . i\frac{b}{2\l^2}\Delta\phi\left(\frac{x-\a}{A\l}\right)\psi\right]\overline{\tu}\right)+ O\left(\l^2\|\psi\|^2_{L^2(\RR^2)}+\frac{1}{\l^2}\|\tu\|^{2}_{L^2(\RR^2)}+\|\tu\|^{2}_{H^1(\RR^2)}\right).
\eea
\end{lemma}

\begin{remark} The virtue of \fref{crc6} is to keep track of the {\it quadratic} terms in $\tu$. In the various situations we will encounter, the key will be to prove that the corresponding quadratic forms both in the boundary term in time and in the RHS of  \fref{crc6} are small deformations of the linearized energy and hence are definite positive. This will generate a control of the form: 
\be
\label{roughenegryindeitty}
\frac{d}{dt}\left\{\|\nabla \tu\|_{L^2}^2+\frac{\|\ut\|_{L^2}^2}{\lambda^2}\right\}\geq \frac{b}{\lambda^2}\left(\int_{|x-\alpha|\lesssim \lambda}|\nabla \ut|^2+\frac{\|\ut\|_{L^2}^2}{\lambda^2}+l.o.t.\right).
\ee
The two keys are first the sign of the above RHS and second the fact that this is a useful information once the blow up speed is known, that is $b\sim \lambda$. The Lyapounov property \fref{crc6} will be the key to estimate the solution backwards from the singularity, which will be used both for the proof of existence and uniqueness.
\end{remark}

{\bf Proof of Lemma \ref{lemma:timederivative}}\\

{\bf step 1} Algebraic derivation of the energetic part.\\

We compute from \fref{un3}:
\bea
\label{cpepjpvudpov}
& & \dt\bigg\{\frac{1}{2}\int |\n\tu|^2 +\half\int \frac{|\tu|^2}{\l^2}-\int k(x)\left[(F(u)-F(w)-F'(w)\cdot\tu)\right]\bigg\}\\
\nonumber & = &- \Re\left(\partial_t\tu,\overline{\Delta \tu-\frac{1}{\lambda^2}\tu+k(x)(f(u)-f(w))}\right)-\frac{\lambda_t}{\lambda^3}\int|\tu|^2\\
\nonumber & - & \Re\left(\partial_tw,k(x)(\overline{f(\tu+w)-f(w)-f'(w)\cdot\tu)}\right)\\
\nonumber & = & \Im\left(\psi,\overline{\Delta \tu-\frac{1}{\lambda^2}\tu+k(x)(f(u)-f(w))}\right)-\frac{1}{\lambda^2}\Im\left(k(x)(f(u)-f(w)),\overline{\ut}\right)\\
\nonumber& - & \frac{\lambda_t}{\lambda^3}\int|\tu|^2-\Re\left(\partial_tw,k(x)\overline{(f(\tu+w)-f(w)-f'(w)\cdot\tu)}\right)\\
\nonumber& = & \Im\left(\psi,\overline{\Delta \tu-\frac{1}{\lambda^2}\tu+k(x)(2|w|^2\ut+\overline{\ut}w^2)}\right)-\frac{1}{\lambda^2}\Im\int k(x)\overline{\ut}^2w^2\\
\nonumber& - & \frac{\lambda_t}{\lambda^3}\int|\tu|^2-\Re\left(\partial_tw,k(x)\overline{(\overline{w}\tu^2+2w|\tu|^2)}\right)\\
\nonumber& + & \Im\left(\psi-\frac{1}{\lambda^2}\ut,\overline{k(x)(f(w+\ut)-f(w)-f'(w)\cdot\ut)}\right)-\Re\left(\partial_tw,k(x)\overline{\ut|\ut|^2}\right)
\eea
where we used that $f'(w)\cdot\ut=2|w|^2\ut+w^2\overline{\ut}.$
We first estimate from \fref{aprioirigeometrical}:
\be
\label{neiohoeghe}
-\frac{\lambda_t}{\lambda^3}\int|\tu|^2=\frac{b}{\lambda^4}\int|\ut|^2-\frac{1}{\lambda^4}(\lambda\lambda_t+b)|\ut|_{L^2}^2=\frac{b}{\lambda^4}\int|\ut|^2+O(\|\ut\|_{H^1}^2).
\ee
It remains to estimate the last line in the RHS \fref{cpepjpvudpov}. For the quadratic and higher terms, we estimate using the a priori bounds \fref{aprioriboundwsw}, \fref{aprioritu}:
\bea
\label{crc3}
\nonumber & & \left|\Im\left(\psi-\frac{1}{\lambda^2}\ut,\overline{k(x)(f(w+\ut)-f(w)-f'(w)\cdot\ut)}\right)\right|\\
\nonumber& = & \left|\Im\left(\psi-\frac{1}{\lambda^2}\ut,\overline{k(x)(\ut^2\overline{w}+2|\tu|^2w+|\tu|^2\tu)}\right)\right|\\
\nonumber& \lesssim & \|\psi\|_{L^2(\RR^2)}\|\tu\|_{L^6}^2(\|w\|_{L^6}+\|\tu\|_{L^6})+(1+\|w\|_{L^2})\|\ut\|_{L^6}^3\\
\nonumber & \lesssim& \|\psi\|_{L^2}\frac{1}{\l^{2/3}}\|\tu\|_{L^2}^{2/3}\|\tu\|^{4/3}_{H^1}+\|\ut\|_{H^1}^2\\
& \lesssim&  \l^2\|\psi\|^2_{L^2}+\|\tu\|^{2}_{H^1}.
\eea
For the cubic terms hitting $w_t$, we replace $w_t$ using \eqref{un2}, integrate by parts and then rely on \fref{aprioriboundwsw} to estimate:
\bea
\label{cnkeneofho}
\nonumber \left|\int k(x)w_t\overline{|\tu|^2\tu}\right|& \lesssim & \|w\|_{H^{3/2}}\||\tu|^2\tu\|_{H^{1/2}(\RR^2)}+\|w\|^3_{L^6}\|\tu\|^3_{L^6}
+\|\psi\|_{L^2(\RR^2)}\|\tu\|^3_{L^6}\\
\nonumber & \lesssim & \frac{1}{\l^{3/2}}\|\tu\|^{1/2}_{L^2}\|\tu\|^{5/2}_{H^1}+\frac{1}{\lambda^2}\|\ut\|_{H^1}^2\|\tu\|_{L^2}+\|\psi\|_{L^2}\|\tu\|^{2}_{H^1}\|\tu\|_{L^2}\\
& \lesssim & \l^2\|\psi\|^2_{L^2}+\|\tu\|^2_{H^1}.
\eea
Injecting \fref{neiohoeghe}, \fref{crc3}, \fref{cnkeneofho} into \fref{cpepjpvudpov} yields the preliminary computation:
\bea
\label{cnkoheofh}
\nonumber & &  \dt\bigg\{\frac{1}{2}\int |\n\tu|^2+\half\int \frac{|\tu|^2}{\l^2}- \int k(x)(F(u)-F(w)-F'(w)\cdot\tu)\bigg\}\\
\nonumber & = & -\frac{1}{\l^2}\Im\left(\int k(x)w^2\overline{\tu}^2\right)-\Re\left(\int k(x)w_t\overline{(2|\tu|^2w+\tu^2\overline{w})}\right)+\frac{b}{\l^2}\int \frac{|\tu|^2}{\l^2}\\
\nonumber & + & \Im\left(\int\left[\Delta\psi-\frac{\psi}{\l^2}+k(x)(2|w|^2\psi-w^2\overline{\psi})\right]\overline{\tu}\right)\\
& + & O\left(\l^2\|\psi\|^2_{L^2}+\frac{1}{\l^2}\|\tu\|^{2}_{L^2}+\|\tu\|^{2}_{H^1}\right).
\eea

{\bf step 2} Algebraic derivation of the localized virial part.\\

Let $$\nabla\tilde{\phi}(t,x)=\frac{b}{\l}A\nabla\phi\left(\frac{x-\a}{A\l}\right),$$ then
\bea
\label{firsttermvirieloc}
 & & \half\dt\left(\frac{b}{\l}\Im\left(\int A\nabla\phi\left(\frac{x-\a}{A\l}\right)\nabla\tu\overline{\tu}\right)\right) \\
 \nonumber & =&  \half \Im\left(\int\partial_t\nabla \tilde{\phi}\cdot \nabla\ut\overline{\ut}\right)+  \Re\left(\int i\partial_t\ut\left[\overline{\half \Delta\tilde{\phi}\tu+\nabla\tilde{\phi}\cdot\nabla\tilde{u}}\right]\right).
\eea
Using \eqref{aprioirigeometrical}, we estimate in brute force:
$$
\left|\partial_t\nabla \tilde{\phi}\right|\lesssim \frac{1}{\lambda^3}(|\lambda^2b_t+b^2|+|\l\l_t+b|)+\frac{b}{\lambda^3}(|\l\alpha_t|+b)\lesssim \frac{1}{\lambda}
$$
from which: 
\be
\label{calculprline}
\left|\Im\left(\int\partial_t\nabla \tilde{\phi}\cdot \nabla\ut\overline{\ut}\right)\right|\lesssim \frac{1}{\lambda}\|\ut\|_{L^2}\|\nabla\ut\|_{L^2}=O\left(\frac{1}{\l^2}\|\tu\|^{2}_{L^2}+\|\tu\|^{2}_{H^1}\right).
\ee
The second term in \fref{firsttermvirieloc} corresponds to the localized Morawetz multiplier, and we get from \fref{un3} and integration by parts:
\bea
\label{nceohoeoud}
& & \Re\left(\int i\partial_t\ut\left[\overline{\half \Delta\tilde{\phi}\tu+\nabla\tilde{\phi}\cdot\nabla\tilde{u}}\right]\right)=\\
\nonumber & = &   \frac{b}{\l^2}\Re\left(\int\nabla^2\phi\left(\frac{x-\a}{A\l}\right)(\n\tu,\overline{\n\tu})\right)-\frac{1}{4}\frac{b}{A^2\l^4}\left(\int\Delta^2\phi\left(\frac{x-\a}{A\l}\right)|\tu|^2\right)\\
\nonumber &  - & \frac{b}{\l}\Re\left(\int A\nabla\phi\left(\frac{x-\a}{A\l}\right)k(x)(|u|^2u-|w|^2w)\cdot\overline{\nabla\tu}\right)\\
\nonumber &  - & \half\frac{b}{\l^2}\Re\left(\int\Delta\phi\left(\frac{x-\a}{A\l}\right)k(x)(|u|^2u-|w|^2w)\overline{\tu}\right)\\
\nonumber & - &\frac{b}{\l}\Re\left(\int A\nabla\phi\left(\frac{x-\a}{A\l}\right)\psi\cdot\overline{\nabla\tu}\right)\ds -\half\frac{b}{\l^2}\Re\left(\int\Delta\phi\left(\frac{x-\a}{A\l}\right)\psi\overline{\tu}\right).
\eea
We now expand the nonlinear terms and estimate the cubic and higher terms:
\bee
& &  \bigg|-\frac{b}{\l}\Re\left(\int A\nabla\phi\left(\frac{x-\a}{A\l}\right)k(x)(2|\tu|^2w+\tu^2\overline{w}+|\tu|^2\tu)\cdot\overline{\nabla\tu}\right)\\
& - & \half\frac{b}{\l^2}\Re\left(\int\Delta\phi\left(\frac{x-\a}{A\l}\right)k(x)(2|\tu|^2w+\tu^2\overline{w}+|\tu|^2\tu)\overline{\tu}\right)\bigg|\\
& \lesssim &(\|\tu\|^3_{L^6(\RR^2)}+\|\tu\|^2_{L^6(\RR^2)}\|w\|_{L^6(\RR^2)})\|\n\tu\|_{L^2(\RR^2)}+\frac{1}{\l}(\|\tu\|^4_{L^4(\RR^2)}+\|\tu\|^3_{L^4(\RR^2)}\|w\|_{L^4(\RR^2)})\\
& \lesssim &\|\tu\|^4_{H^1(\RR^2)}+\|\tu\|^{7/3}_{H^1(\RR^2)}\|\tu\|^{2/3}_{L^2(\RR^2)}\|w\|_{L^2(\RR^2)}^{1/3}\|w\|^{2/3}_{H^1(\RR^2)}\\
&+& \frac{1}{\l}\left(\|\tu\|^2_{L^2(\RR^2)}\|\tu\|^2_{H^1(\RR^2)}+\|\tu\|^{3/2}_{H^1(\RR^2)}\|\tu\|^{3/2}_{L^2(\RR^2)}\|w\|_{L^2(\RR^2)}^{1/2}\|w\|^{1/2}_{H^1(\RR^2)}\right)\\
& \lesssim &\ds \|\tu\|^2_{H^1(\RR^2)}
\eee
where we have used \eqref{aprioritu} and \eqref{aprioriboundwsw}.

The remaining quadratic terms in \fref{nceohoeoud} are integrated by parts:
\bea
\label{lpusfpfwu}
\nonumber &- &\frac{b}{\l}\Re\left(\int A\nabla\phi\left(\frac{x-\a}{A\l}\right)\psi\cdot\overline{\nabla\tu}\right)\ds -\half\frac{b}{\l^2}\Re\left(\int\Delta\phi\left(\frac{x-\a}{A\l}\right)\psi\overline{\tu}\right)\\
& = & \Im \left(\int \left[i\frac{b}{\lambda}A\nabla\phi\left(\frac{x-\a}{A\l}\right)\cdot\nabla \psi+i\frac{b}{2\l^2}\Delta\phi\left(\frac{x-\a}{A\l}\right)\psi\right]\overline{\ut}\right),
\eea
\bea
\label{bingo}
\nonumber & - & \frac{b}{\l}\Re\left(\int A\nabla\phi\left(\frac{x-\a}{A\l}\right)k(x)(2|w|^2\ut+w^2\overline{\ut})\cdot\overline{\nabla\tu}\right)\\
\nonumber & - &  \half\frac{b}{\l^2}\Re\left(\int\Delta\phi\left(\frac{x-\a}{A\l}\right)k(x)(2|w|^2\ut+w^2\overline{\ut})\overline{\tu}\right)\\
\nonumber & = &  \frac{b}{\l}\Re\left(\int A\nabla\phi\left(\frac{x-\a}{A\l}\right)k(x)(2|\tu|^2w+\tu^2\overline{w})\cdot\overline{\nabla w}\right)\\
 & + & \frac{b}{2\l}\Re\left(\int A\nabla\phi\left(\frac{x-\a}{A\l}\right)\cdot\nabla k(x)\left[2|\ut|^2|w|^2+\overline{\ut}^2w^2\right]\right).
\eea
The last term gains an extra $\l$ smallness:
\bea
\label{lplfuvjeoiog}
\nonumber & & \left|\frac{b}{\l}\Re\left(\int A\nabla\phi\left(\frac{x-\a}{A\l}\right)\cdot\nabla k(x)\left[2|\ut|^2|w|^2+\overline{\ut}^2w^2\right]\right)\right|\\
& \lesssim & \|\ut\|_{L^4}^2\|w\|_{L^4}^2\lesssim \frac{1}{\lambda}\|\nabla \ut\|_{L^2}\|\ut\|_{L^2}=O\left(\|\ut\|_{H^1}^2+\frac{\|\ut\|_{L^2}^2}{\l^2}\right).
\eea
Injecting \fref{lpusfpfwu}, \eqref{bingo}, \fref{lplfuvjeoiog} into \fref{nceohoeoud} yields after a further integration by parts:
\bee
& & \Re\left(\int i\partial_t\ut\left[\overline{\half \Delta\tilde{\phi}\tu+\nabla\tilde{\phi}\cdot\nabla\tilde{u}}\right]\right)\\
\nonumber & = &   \frac{b}{\l^2}\Re\left(\int\nabla^2\phi\left(\frac{x-\a}{A\l}\right)(\n\tu,\overline{\n\tu})\right)-\frac{1}{4}\frac{b}{A^2\l^4}\left(\int\Delta^2\phi\left(\frac{x-\a}{A\l}\right)|\tu|^2\right)\\
\nonumber & + & \frac{b}{\l}\Re\left(\int A\nabla\phi\left(\frac{x-\a}{A\l}\right)k(x)(2|\tu|^2w+\tu^2\overline{w})\cdot\overline{\nabla w}\right)\Bigg)\\
\nonumber & + & \Im\left(\int\left[i\frac{b}{\l}A\nabla\phi\left(\frac{x-\a}{A\l}\right)\cdot\nabla\psi+i\frac{b}{2\l^2}\Delta\phi\left(\frac{x-\a}{A\l}\right)\psi\right]\overline{\tu}\right)+O\left(\|\tu\|^2_{H^1}+\frac{\|\ut\|_{L^2}}{\l^2}\right).
\eee
Together with \fref{cnkoheofh}-\eqref{calculprline}, this yields \fref{crc6}. This concludes the proof of Lemma \ref{lemma:timederivative}.


\subsection{Backwards propagation of smallness}


The first application of the energy estimate \fref{crc6} is a bootstrap control on critical mass solutions to \fref{nlsk:0}. More precisely, let $u$ be a solution to \fref{nlsk:0} defined on $[\tilde{t_0},0)$. Let $\tilde{t_0}<t_1< 0$ and assume that $u$ admits on $[\tilde{t_0},t_1]$ a geometrical decomposition on the form: 
$$u(t,x)=\frac{1}{[k(\alpha(t))]^{\frac{1}{2}}}\frac{1}{\lambda(t)}(Q_{\P(t)}+\e)\left(t,\frac{x-\alpha(t)}{\lambda(t)}\right)e^{i\gamma(t)}
$$ where $\e$ satisfies the  orthogonality conditions \fref{ortho1}- \fref{ortho5} and $\|\e(t)\|_{H^1}+|\mathcal P(t)|\ll 1$. Let 
\be
\label{defcvcvudsbis}
\ut(t,x)=\frac{1}{[k(\alpha(t))]^{\frac{1}{2}}}\frac{1}{\lambda(t)}\e\left(t,\frac{x-\alpha(t)}{\lambda(t)}\right)e^{i\gamma(t)}.
\ee
Assume that the energy $E_0$ satisfies the lower bound \eqref{condblowup}\footnote{which will be proved in \fref{boundinsideezero}}, and define $C_0$ as:
\be
\label{defcvcvuds}
C_0=\sqrt{\frac{\|yQ\|_{L^2}^2}{8E_0+\int\nabla k^2(0)(y,y)Q^4}}.
\ee
We claim the following backwards propagation estimates:

\begin{lemma}[Backwards propagation of smallness]
\label{propbaczkardw}
Assuming that there holds for some $t_1<0$ close enough to 0: 
\be
\label{initltwonorm}
\left|\|u\|_{L^2}-\|Q\|_{L^2}\right|\lesssim \lambda^4(t_1),
\ee
\be
\label{initalaprioriboundut}
\|\nabla \ut(t_1)\|_{L^2}^2+\frac{\|\ut(t_1)\|^2_{L^2}}{\lambda^2(t)} \lesssim  \lambda^2(t_1)
\ee
\be
\label{aprioriinitialparamaters}
\left|\frac{\beta}{\lambda}(t_1)\right|+\left|\frac{\alpha}{\lambda}(t_1)\right|\lesssim \lambda(t_1), \ \ \left|\lambda(t_1)+\frac{t_1}{C_0}\right|\lesssim \lambda^3(t_1), \ \ \left|\frac{b(t_1)}{\lambda(t_1)}-\frac{1}{C_0}\right|\lesssim\lambda^2(t_1).
\ee
Then there exists a backwards time $t_0$ depending only on $C_0$ such that $\forall t\in [t_0,t_1]$,
 \be
 \label{backwardspropagation}
 \|\nabla \ut(t)\|_{L^2}^2+\frac{\|\ut(t)\|^2_{L^2}}{\lambda^2(t)} \lesssim  \|\nabla \ut(t_1)\|_{L^2}^2+\frac{\|\ut(t_1)\|^2_{L^2}}{\lambda^2(t_1)}+\lambda^6(t)
 \ee
 \be
 \label{propparametersbis}
 \left|\frac{b}{\lambda}(t)-\frac{1}{C_0}\right|\lesssim \lambda^2(t),
 \ee
 \be
 \label{propparametersbivmdjvdps} \ \ \left|\lambda(t)+\frac{t}{C_0}\right|\lesssim\lambda^3(t),
 \ee
 \be
 \label{propparameters}
  |\alpha(t)|\lesssim\lambda^2(t), \ \ \left|\frac{\beta}{\lambda}(t)\right|\lesssim \lambda(t).
  \ee
\end{lemma}

{\bf Proof of Lemma \ref{propbaczkardw}}\\

Let us consider from the continuity $u\in \mathcal C([t_0,t_1],H^1)$ a backwards time $t_0$ such that $\forall t\in [t_0,t_1]$,
 \be
 \label{backwardspropagationboot}
 \|\ut(t)\|_{L^2}\lesssim K\lambda^2(t), \ \ \|\ut(t)\|_{H^1}\leq K\lambda(t),
 \ee
 \be
\label{aprioriinitialparamatersboot}
\left|\frac{\beta(t)}{\lambda(t)}\right|+\left|\frac{\alpha(t)}{\lambda(t)}\right|\leq K \lambda(t), \ \ \left|\lambda(t)+\frac{t}{C_0}\right|\leq K \lambda^3(t), \ \ \left|\frac{b(t)}{\lambda(t)}-\frac{1}{C_0}\right|\leq K\lambda^2(t).
\ee
for some large enough universal constant $K>0$. Then we claim that \fref{backwardspropagation},  \fref{propparametersbis}, \fref{propparameters} hold on $[t_0,t_1]$ hence improving the bounds \fref{backwardspropagationboot}, \fref{aprioriinitialparamatersboot} on $[t_0,t_1]$ for $t_0=t_0(C_0)$ small enough independent of $t_1$.\\

{\bf step 1} Monotonicity of the norm.\\

Let us apply Lemma \ref{lemma:timederivative} with:
\be
\label{defwqtilde}
w(t,x)=\tilde{Q}(t,x)=\frac{1}{[k(\alpha(t))]^{\frac{1}{2}}}\frac{1}{\lambda(t)}Q_{\P(t)}\left(\frac{x-\alpha(t)}{\lambda(t)}\right)e^{i\gamma(t)},
\ee
Let $\mathcal I$ be given by \fref{defI}, we claim that \fref{crc6} implies the following coercivity property: 
\be
\label{tobeproved}
\frac{d\mathcal I}{dt}\geq \frac{b}{\lambda^4}\int|\ut|^2+0(K^4\lambda^5+\|\ut\|_{H^1}^2).
\ee
Assume \fref{tobeproved}. By Sobolev embedding and the smallness of $\e$, we have the rough upper bound on $\mathcal I$:  
\be
\label{upepri}
|\mathcal I|\lesssim  \|\nabla \ut\|^2_{L^2}+\frac{1}{\lambda^2}\|\ut\|_{L^2}^2.
\ee Now after renormalization, the proximity of $\qp$ to $Q$ ensures that the quadratic part of $\mathcal I$ is  a small deformation of the linearized energy and hence the coercitivity property \fref{coerclinearenergy} with the choice of orthogonality conditions on $\e$ ensures : 
\bea
\label{coebbeof}
\nonumber \mathcal I& = & \frac{1}{2}\int |\n\tu|^2 +\half\int \frac{|\tu|^2}{\l^2}-\int k(x)\left(F(w+\tu)-F(w)-F'(w)\cdot\tilde{u}\right)\\
& + & \half\frac{b}{\l}\Im\left(\int A\nabla\phi\left(\frac{x-\a}{A\l}\right)\nabla\tu\overline{\tu}\right)\\
\nonumber & = & \frac{1}{2\lambda^2} \left[(L_+\e_1,\e_1)+(L_-\e_2,\e_2)+o\left(\|\e\|_{H^1}^2\right)\right] \geq \frac{c_0}{\lambda^2}\left[\int\|\e\|_{H^1}^2-(\e_1,Q)^2\right].
\eea
We next estimate from the conservation of the $L^2$ norm \fref{degeneracyeq}, \fref{initltwonorm}, \fref{aprioriinitialparamatersboot} and \fref{cm5}:
$$
|\Re(\e,\overline{\qp})|\lesssim \|\e\|_{L^2}^2+\lambda^4(t)+|\a|^2+\left|\int |u|^2-\int Q^2\right|\lesssim \|\e\|_{L^2}^2+K^2\lambda^4(t)
$$
from which
\be
\label{ltwoconsrbis}
(\e_1,Q)^2\lesssim o(\|\e\|_{L^2}^2)+K^4\lambda^8(t).
\ee
We then integrate \fref{tobeproved} in time and inject \fref{upepri}, \fref{coebbeof}, \fref{ltwoconsrbis} to conclude:
\bee
\|\nabla \ut(t)\|_{L^2}^2+\frac{\|\ut(t)\|^2_{L^2}}{\lambda^2(t)} & \lesssim &  \|\nabla \ut(t_1)\|_{L^2}^2+\frac{\|\ut(t_1)\|^2_{L^2}}{\lambda^2(t)}+K^4\lambda^6(t)+\int_t^{t_1}(\|\ut(\tau)\|_{H^1}^2+K^4\lambda^5(\tau))d\tau\\
& \lesssim &  \|\nabla \ut(t_1)\|_{L^2}^2+\frac{\|\ut(t_1)\|^2_{L^2}}{\lambda^2(t)}+K^4\lambda^6(t)+\int_t^{t_1}\left[\|\nabla \ut(\tau)\|_{L^2}^2+\frac{\|\ut(\tau)\|^2_{L^2}}{\lambda^2(\tau)}\right]d\tau
\eee
 for $t_0=t_0(C_0)$ small enough, and \fref{backwardspropagation} follows from Gronwall's lemma. It implies in particular together with \fref{initalaprioriboundut}:
 \be
 \label{estparmaaters}
 \|\nabla \ut(t)\|_{L^2}^2+\frac{\|\ut(t)\|^2_{L^2}}{\lambda^2(t)}\lesssim \lambda^2(t)
 \ee 
 and closes the bootstrap of \fref{backwardspropagationboot}.\\
 
 {\bf step 2} Integration of the law for the parameters.\\
 
 We now integrate the law for the parameters. Indeed, \fref {law2}, \fref{aprioriinitialparamatersboot} and \fref{estparmaaters} imply like for the proof of \fref{ie:rmk2}:
\bea
\label{newetsinvoeo}
\nonumber |b_s+b^2-d_0(\alpha,\alpha)|+\left|\frac{\l_s}{\l}+b\right|+|\b_s+b\b-c_0(\a)\l-\b_3\l^3|+\left|\frac{\a_s}{\l}-2\b\right| & \lesssim  &\lambda^4+K^2\lambda^5\\
& \lesssim &\lambda^4.
\eea
One should be careful when reintegrating the above system which is stable thanks to \fref{newetsinvoeo} {\it and the sign} $d_0(\alpha,\alpha)\leq 0$ provided by \fref{defdzeroalpha}:
$$\left(\frac{b}{\lambda}\right)_s=\frac{b_s+b^2-d_0(\alpha,\alpha)}{\lambda}-\frac{b}{\lambda}\left(\lsl+b\right)+\frac{d_0(\alpha,\alpha)}{\l}\lesssim \l^3$$
and hence: 
\be
\label{cnoeieouodfjjo}
\frac{1}{C_0}-\frac{b}{\l}(s)\leq \frac{1}{C_0}-\frac{b}{\l}(s_1)+\l^2(s)\lesssim \l^2(s)
\ee where we used \fref{aprioriinitialparamaters}. We now write down the  conservation of energy at t and add the conservation of $L^2$ norm. We get using verbatim the same algebra like the one which led to \fref{cme5} and using \fref{estparmaaters}:
\bee
\label{cjeououte}
& & \frac{b^2}{8}\int |y|^2Q^2+\frac{|\b|^2}{2}\int Q^2-\frac{\nabla^2k(0)(\a,\a)}{4}\int Q^2\\
\nonumber &  = & \l^2\left(E_0+\frac{1}{8}\int \nabla^2k(0)(y,y)Q^4\right)+k(\a)\left(\int |u|^2-\int Q^2\right)+O(\l^4),
\eee
and thus from the choice of $C_0$ \fref{defcvcvuds} and \fref{initltwonorm}:
$$\frac{|\b|^2}{\lambda^2}+\frac{|\alpha|^2}{\lambda^2}\lesssim\frac{1}{C_0^2}-\frac{b^2}{\lambda^2}+\l^2 \lesssim \frac{1}{C_0}-\frac{b}{\lambda}+\l^2\lesssim \lambda^2$$ where we used \fref{cnoeieouodfjjo} in the last step. This together with \fref{cnoeieouodfjjo} again yields: 
\be
\label{cnkojheofueou}
\left|\frac{b}{\l}-\frac{1}{C_0}\right|\lesssim \l^2.
\ee
We eventually have from \fref{newetsinvoeo} for the scaling parameter: $$\left|\lambda_t+\frac{b}{\lambda}\right|\lesssim \lambda^3$$ and hence using \eqref{aprioriinitialparamaters}, \fref{cnkojheofueou}:
$$\left|\lambda(t)+\frac{t}{C_0}\right| \lesssim  \left|\lambda(t_1)+\frac{t_1}{C_0}\right| +\int_t^{t_1}\left|\frac{b}{\lambda}(\tau)-\frac{1}{C_0}\right|d\tau+\lambda^3(t)\lesssim \lambda^3(t).$$
This concludes the proof of \fref{propparametersbis}, \fref{propparametersbivmdjvdps}, \fref{propparameters} assuming \fref{tobeproved}.\\

{\bf step 3} Coercivity of the quadratic form in the RHS of \fref{crc6}.\\

We now turn to the proof of the Lyapounov property \fref{tobeproved} and start with computing more explicitly the quadratic terms in the RHS of \fref{crc6} for $w=\tilde{Q}$ given by \fref{defwqtilde}:
\bea
\label{defquu}
\nonumber \mathcal{K}(\ut)& =&-\frac{1}{\l^2}\Im\left(\int k(x)w^2\overline{\tu}^2\right)-\Re\left(\int k(x)w_t\overline{(2|\tu|^2w+\tu^2\overline{w})}\right)\\
\nonumber &+ & \frac{b}{\l^2}\Bigg(\int \frac{|\tu|^2}{\l^2}+\Re\left(\int\nabla^2\phi\left(\frac{x-\a}{A\l}\right)(\n\tu,\overline{\n\tu})\right)-\frac{1}{4A^2}\left(\int\Delta^2\phi\left(\frac{x-\a}{A\l}\right)\frac{|\tu|^2}{\l^2}\right)\\
&+ & \frac{b}{\lambda}\Re\left(\int A\nabla\phi\left(\frac{x-\a}{A\l}\right)k(x)(2|\tu|^2w+\tu^2\overline{w})\cdot\overline{\nabla w}\right)\Bigg),
\eea
we claim: 
\be
\label{firtboundq}
\mathcal{K}(\tu)\geq \frac{c}{\lambda^3}\left(\int|\nabla \e|^2e^{-\frac{|y|}{\sqrt{A}}}+\int|\e|^2\right)+O(K^4\lambda^5)
\ee
for some universal constant $c>0$. Indeed, first observe from \fref{defdzeroalpha}, \fref {law2}, \fref{backwardspropagationboot} and \fref{aprioriinitialparamatersboot}:
\be
\label{newetsinvoeobis}
|Mod(t)|+|d_0(\alpha,\alpha)|+|d_1(\alpha,\alpha)|\lesssim K^2\lambda^4.
\ee
We then first compute from \fref{defwqtilde}:
\bee
\tq_t & = &\ds\bigg(-\frac{(\l\a_t-2\b)\n k(\a)}{2k(\a)\l}-\frac{2\b\n k(\a)}{2k(\a)\l}-\frac{\l\l_t+b}{\l^2}+ \frac{b}{\l^2}+\frac{i}{\l^2}+i\frac{\l^2\widetilde{\gamma}_t-|\b|^2}{\l^2}+i\frac{|\b|^2}{\l^2}\bigg)\tq\\
&+ &\left(-\frac{\l\l_t+b}{\l^2}\left(\frac{x-\a}{\l}\right)+\frac{b}{\l^2}\left(\frac{x-\a}{\l}\right)-\frac{\l\a_t-2\b}{\l^2}-\frac{2\b}{\l^2}\right)\cdot \frac{1}{\lambda[k(\a)]^{\frac{1}{2}}}\n\qp\left(\frac{x-\a(t)}{\l(t)}\right)e^{i\gamma(t)}\\
& +& \frac{1}{k(\a(t))^{1/2}\l(t)}\P_t\frac{\partial\qp}{\partial\P}\left(\frac{x-\a(t)}{\l(t)}\right)e^{i\gamma(t)}\\
& = &\ds\left(\frac{i}{\l^2}+\frac{b}{\l^2}\right)\tq+\frac{b}{\l}\left(\frac{x-\a}{\l}\right)\cdot\n\tq+ O\left(\frac{K}{\lambda}e^{-\frac{|x-\alpha(t)|}{\lambda(t)}}\right)
\eee
where we used \fref{newetsinvoeobis}, \fref{aprioriinitialparamatersboot} in the last step. 
We may thus use the exponential decay of $\qp$ to conclude after renormalization:
\bee
\label{babar}
\nonumber &- & \Re\left(\int k(x)\tq_t\overline{(2|\tu|^2\tq+\tu^2\overline{\tq})}\right)= \ds\frac{1}{\l^2}\Im\left(\int k(x)\tq\overline{(2|\tu|^2\tq+\tu^2\overline{\tq})}\right)\\
\nonumber &- &\frac{b}{\l^2}\Re\left(\int k(x)(2|\tu|^2\tq+\tu^2\overline{\tq})\overline{\tq}\right) -\frac{b}{\l}\Re\left(\int \left(\frac{x-\a}{\l}\right) k(x)(2|\tu|^2\tq+\tu^2\overline{\tq})\cdot\overline{\nabla \tq}\right)\\
& + & \frac{1}{\lambda^3}O\left(K\lambda\|\e\|_{L^2}^2\right).
\eee
We now inject this estimate into \fref{defquu} and rewrite the result in renormalized variables:
\bee
\mathcal{K}(\tu)
&= & \frac{b}{\l^4}\left\{\Re\left(\int\nabla^2\phi\left(\frac{y}{A}\right)(\n\e,\overline{\n\e})\right)+\int|\e|^2 \right .\\
&- &\left . \int k(\alpha+\lambda y)((|\qp|^2+2\S^2)\e_1^2+4\S\T\e_1\e_2+(|\qp|^2+2\T^2)\e_2^2)- \frac{1}{4A^2}\int\Delta^2\phi\left(\frac{y}{A}\right)|\e|^2\right\}\\
& + & \frac{b}{\lambda^3}\Re\left(\int \left(A\nabla\phi\left(\frac{y}{A\l}\right)-y\right)k(\alpha+ \lambda y)(2|\e|^2\qp+\e^2\overline{\qp})\cdot\overline{\nabla \qp}\right)\\
& + &  \frac{1}{\lambda^3}O\left(K\lambda\|\e\|_{L^2}^2\right).
\eee
From the proximity of $\qp$ to $Q$, the above quadratic form is for $A$ large enough a small deformation of the localized in $A$ linearized energy, and hence \fref{coerclinearenergy},  \fref{aprioriinitialparamatersboot} and our choice of orthogonality conditions ensure for $A$ large enough:
$$\mathcal{K}(\tu)\gtrsim\frac{1}{\lambda^3}\left[\int|\nabla \e|^2e^{-\frac{|y|}{\sqrt{A}}}+\int|\e|^2-(\e_1,Q)^2\right]\gtrsim \frac{1}{\lambda^3}\left[\int|\nabla \e|^2e^{-\frac{|y|}{\sqrt{A}}}+\int|\e|^2\right]+O(K^4\lambda^5)
$$
where we used \fref{ltwoconsrbis} in the last step, this is \fref{firtboundq}.\\

{\bf step 4} Control of the remainder terms in the RHS of \fref{crc6}.\\

It remains to control the $\psi$ terms in \fref{crc6}. According to \fref{un2}, \fref{defwqtilde} and the construction of $\qp$, we have: 
\bea
\label{defpsin}
\psi & = & \frac{1}{\l^3k(\a)^{1/2}}\Bigg[i(b_s+b^2)\p_b\qp+i\left(\frac{\l_s}{\l}+b\right)\l\p_\l\qp\\
 \nonumber & + & i\left(\frac{\a_s}{\l}-2\b\right)\l\p_\a\qp+i(\b_s+b\b-\l c_0(\a)-\b_3\l^3-\b_4\l^4)\p_\b\qp- i\left(\frac{\l_s}{\l}+b\right)\Lambda\qp\\
 \nonumber  & - & i\left(\frac{\a_s}{\l}-2\b\right)\cdot\left(\n\qp+\frac{\l}{2}\frac{\n k(\a)}{k(\a)}\qp\right)+(\tgamma_s-|\b|^2)\qp+\psi_{\P}\Bigg]\left(\frac{x-\a(t)}{\l(t)}\right)e^{i\gamma(t)}
\eea
where $\psi_P$ is the remainder in the construction of $\qp$:
$$\psi_P=O(\l^5e^{-c|y|}).
$$
Let us start with a rough bound on $\psi$: from \fref{law2}, \fref{backwardspropagationboot}, \fref{aprioriinitialparamatersboot}, there holds for $i=0,1,2$:
\bea
\label{roughboundpsi}
\nonumber |\nabla^i\psi(x)|& \lesssim & \frac{1}{\lambda^{3+i}}e^{-\frac{|x-\alpha(t)|}{\lambda(t)}}\left[|Mod(t)|+|\alpha|^2+\lambda^5\right]\\
& \lesssim & \frac{1}{\lambda^{3+i}}e^{-\frac{|x-\alpha(t)|}{\lambda(t)}}\left(K\lambda^2\|\e\|_{L^2}+K^2\lambda^4\right)
 \eea
and thus:
\be
\label{eusgidgweitwei}
\|\nabla^i\psi\|_{L^2}\lesssim \frac{1}{\lambda^{2+i}}\left[K\lambda^2\|\e\|_{L^2}+K^2\lambda^4\right].
\ee
This yields in particular the bounds:
\be
\label{dpojepfe}
\lambda^2\|\psi\|_{L^2}^2\lesssim \|\e\|_{L^2}^2+K^4\lambda^6,
\ee
\bea
\label{dpojepfebis}
\nonumber & & \left|\Im\int\left(i\frac{b}{\l}A\nabla\phi\left(\frac{x-\a}{A\l}\right)\cdot\nabla\psi+i\frac{b}{2\l^2}\Delta\phi\left(\frac{x-\a}{A\l}\right)\psi\right)\overline{\tu}\right|\\
\nonumber & \lesssim & \|\nabla\psi\|_{L^2}\|\ut\|_{L^2}+\frac{\|\psi\|_{L^2}}{\lambda}\|\ut\|_{L^2}\lesssim \frac{K^2}{\lambda^3}\left[\lambda\|\e\|_{L^2}+\lambda^4\right]\|\e\|_{L^2}\\
&  \lesssim &  o\left(\frac{\|\e\|_{L^2}^2}{\lambda^3}\right)+K^4\lambda^5.
\eea
The rough bound \fref{roughboundpsi} is not enough to take care of the remainder term in \fref{crc6} for which we need a further cancellation. From \fref{defpsin} and the construction of $\qp$, we have $\psi=\psi_1+\psi_2$ with
\bea
\label{plispsi}
\nonumber \psi_1 =  \frac{1}{\l^3k(\a)^{1/2}}\Big[(b_s+b^2)\frac{|y|^2}{4}Q-(\b_s+b\b-\l c_0(\a)-\b_3\l^3-\b_4\l^4)yQ\\
 -i \left(\frac{\l_s}{\l}+b\right)\Lambda Q-i\left(\frac{\a_s}{\l}-2\b\right)\n Q+ (\tgamma_s-|\beta|^2)Q\Big]\left(\frac{x-\a(t)}{\l(t)}\right)e^{i\gamma(t)}
\eea
and for $i=0,1,2$:
$$
|\nabla^i\psi_2(x)|\lesssim  \frac{1}{\lambda^{3+i}}e^{-\frac{|x-\alpha(t)|}{\lambda(t)}}\left[|\mathcal P|Mod(t)+K^2\lambda^5\right]\lesssim \frac{1}{\lambda^{3+i}}e^{-\frac{|x-\alpha(t)|}{\lambda(t)}}\left[\lambda^2\|\e\|_{L^2}+K^2\lambda^5\right].
 $$
The contribution of this remainder term $\psi_2$ is estimated in brute force:
\bea
\label{cneoneogehoeueoun}
\nonumber & & \left|\Im\int\left(\Delta\psi_2-\frac{\psi_2}{\l^2}+k(x)(2|w|^2\psi_2-w^2\overline{\psi_2})\right)\overline{\tu}\right|\\
\nonumber & \lesssim & \left[\|\Delta \psi_2\|_{L^2}+\frac{\|\psi_2\|_{L^2}}{\l^2}+\|\psi_2\|_{L^{\infty}}\|w\|_{L^4}^2\right]\|\ut\|_{L^2}\\
\nonumber & \lesssim & \frac{1}{\lambda^4}(\lambda^2\|\e\|_{L^2}+K^2\l^5)\|\e\|_{L^2}\lesssim o\left(\frac{\|\e\|_{L^2}^2}{\lambda^3}\right)+K^4\lambda^5.
\eea
We now observe that the $\psi_1$ term lies after renormalization in the generalized null space of $L$ and thus a $O(\mathcal P)$ factor is gained using the choice of orthogonality conditions on $\e$. In other words, after renormalization:
\bea
\label{cneoneogehoeueo}
\nonumber & & \left|\Im\int\left(\Delta\psi_1-\frac{\psi_1}{\l^2}+k(x)(2|w|^2\psi_1-w^2\overline{\psi_1}\right)\overline{\tu}\right|\\
\nonumber & \lesssim &\frac{|Mod(t)|+|\alpha|^2}{\lambda^4}\big[|(\e_2,L_-(|y|^2Q))|+|(\e_2,L_-(yQ))|+|(\e_1,L_+(\n Q))|\\
\nonumber & & +|(\e_2,L_-Q)|+O(|\mathcal P|\|\e\|_{L^2})\big]+\frac{1}{\l^4}\left|\frac{\l_s}{\l}+b\right||(\e_1,L_+(\Lambda Q))|\\
\nonumber & \lesssim & \frac{\lambda\|\e\|_{L^2}+K^2\lambda^4}{\lambda^4}|\mathcal P|\|\e\|_{L^2}+\frac{\lambda\|\e\|_{L^2}+K^2\lambda^5}{\lambda^4}(\l\|\e\|_{L^2}+K^2\l^4)\\
& \lesssim & o\left(\frac{\|\e\|_{L^2}^2}{\l^3}\right)+K^4\l^5,
\eea
where we have used in particular the estimate \eqref{law2:bis} and the conservation of mass \eqref{ltwoconsrbis} to bound the term $\frac{1}{\l^4}|\frac{\l_s}{\l}+b||(\e_1,L_+(\Lambda Q))|$. 

We now inject \fref{dpojepfe}, \fref{dpojepfebis}, \fref{cneoneogehoeueoun} and \fref{firtboundq} into \fref{crc6} and \fref{tobeproved} follows.\\
This concludes the proof of Lemma \ref{propbaczkardw}.


\subsection{Existence of critical mass blow up solutions}


A first consequence of the backwards propagation estimates of Lemma \ref{propbaczkardw} is that one can integrate the flow backwards from the singularity to prove the existence of critical mass finite time blow up solutions. This approach is similar to the one in \cite{Merlemulti}, \cite{Marteluniqueness}, \cite{MMmulti}, \cite{KMR}.

\begin{proposition}[Existence of critical mass blow up solutions]
\label{existenceprop}
Let $$\gamma_{0}\in \RR, \ \ E_0>\frac{1}{8}\int \nabla^2k(0)(y,y)Q^4,$$ and $C_0$ given by \fref{constantcnot}, then there exists $t_0<0$ and a solution $u_c\in \mathcal C([t_0,0),H^{\frac{3}{2}})$ to \fref{nlsk:0} which blows up at $T=0$ with $$E(u_c)=E_0 \ \ \mbox{and} \ \ \|u_c\|_{L^2}=\|Q\|_{L^2}.$$ Moreover, the solution admits on $[t_0,0)$ a geometrical decomposition: 
\be
\label{geomdecompuc}
u_c(t,x)=\frac{1}{[k(\alpha_c(t))]^{\frac{1}{2}}}\frac{1}{\lambda_c(t)}(Q_{\P_c(t)}+\e_c)\left(t,\frac{x-\alpha_c(t)}{\lambda_c(t)}\right)e^{i\gamma_c(t)}=\tilde{Q_c}+\ut_c
\ee where $\e_c$ satisfies the orthogonality  conditions \fref{ortho1}-\fref{ortho5}, and there holds the bounds: 
\be
\label{bounduc}
\|\ut_c\|_{L^2}\lesssim \lambda_c^4, \ \ \|\ut_c\|_{H^1}\lesssim \lambda_c^3, \ \ \|\ut_c\|_{H^{\frac{3}{2}}}\lesssim \lambda^{\frac{3}{2}}_c,
\ee
\be
\label{estcrucailk}
\lambda_c+\frac{t}{C_0}=O(\lambda_c^3), \ \ \frac{b_c}{\lambda_c}-\frac{1}{C_0}=O(\lambda^2_c), 
\ee
\be
\label{estcrucailkbis}
\ \ |\alpha_c|+|\beta_c|\lesssim \l_c^2,\ \ \gamma_c=-\frac{C_0^2}{t}+\gamma_0+O(\lambda_c).
\ee
\end{proposition}

{\bf Proof of Proposition \ref{existenceprop}}\\

{\bf step 1} Backwards uniform bounds.\\

 Let a sequence $t_n\to 0$ and $u_n$ be the solution to \eqref{nlsk:0} with initial data at $t=t_n$ given by:
\be\label{ex2}
u_n(t_n,x)=\frac{1}{\l_n}Q_{\P_n}\left(\frac{x}{\l_n}\right)e^{i\gamma_n(t_n)}
\ee
with $\P_n=(b_n(t_n), \l_n(t_n),\b_n(t_n), \a_n(t_n))$ and:
\be\label{ex9}
b_n(t_n)=-\frac{t_n}{C_0^2},\,\l_n(t_n)=-\frac{t_n}{C_0},\,\a_n(t_n)=\b_n(t_n)=0, \ \ \gamma_n(t_n)=\gamma_0-\frac{C_0^2}{t_n}.
\ee
We have by \eqref{compm4}:
$$\|u_n(t_n)\|^2_{L^2}=\int Q^2+O(t_n^4).
$$
and $\tilde{u}_n(t_n)=0$ by construction. Hence $u_n$ satisfies at $t_1=t_n$ the assumptions of Lemma \ref{propbaczkardw}, and thus we can find a time $t_0$ {\it independent of $n$} such that $\forall t\in [t_0,t_n)$, $u_n$ admits a geometrical decomposition 
$$u_n(t,x)=\frac{1}{[k(\alpha_n(t))]^{\frac{1}{2}}}\frac{1}{\lambda_n(t)}Q_{\P_n(t)}\left(t,\frac{x-\alpha_n(t)}{\lambda_n(t)}\right)e^{i\gamma_n(t)}+\ut_n
$$
with uniform bounds in $n$: 
\be
 \label{backwardspropagationn}
 \|\nabla \ut_n(t)\|_{L^2}^2+\frac{\|\ut_n(t)\|^2_{L^2}}{\lambda_n^2(t)} \lesssim \lambda_n^6(t)
 \ee
 \be
 \label{propparametersbisbis}
\left |\frac{b_n}{\lambda_n}(t)-\frac{1}{C_0}\right|\lesssim \lambda_n^2(t), \ \ \left|\lambda_n(t)+\frac{t}{C_0}\right|\lesssim \lambda_n^3(t),\ \ |\alpha_n(t)|+|\beta_n(t)|\lesssim\lambda_n^2(t).
  \ee
From standard Strichartz bounds, this implies the uniform $H^{\frac{3}{2}}$ bound: 
\be
\label{esthtwo}
\|\ut_n\|_{L^{\infty}([t,t_n],H^{\frac{3}{2}})}\lesssim \lambda_n(t)^{\frac{3}{2}},
\ee
which we prove in step 2.\\
The $H^1$ compactness of $u_n(t_0)$ is now a consequence of a standard localization procedure. Indeed, let a cut off function $\chi(x)=0$ for $|x|\leq 1$ and $\chi(x)=1$ for $|x|\geq 2$ and $\chi_R(x)=\chi(\frac{x}{R})$, then $$\left|\frac{d}{dt}\int \chi_R|u_n|^2\right|=2\left|Im\left(\int \nabla\chi_R\cdot\nabla u_n)\overline{u_n}\right)\right|\lesssim \frac{1}{R},$$ $$ \left|\frac{d}{dt}\int \chi_R\left(\frac{1}{2}|\nabla u_n|^2-\frac{1}{4}\int k(x)|u_n|^4\right)\right|=\left|Im\left(\int \nabla\chi_R\cdot\nabla u_n\overline{(\Delta u_n+k(x)u_n|u_n|^2)}\right)\right|\lesssim \frac{1}{R}$$ where we used \fref{backwardspropagationn}, \fref{propparametersbisbis} and \fref{esthtwo}. Integrating this backwards from $t_1$ to $t_0$ using  \fref{ex2}, \fref{ex9} and \fref{esthtwo} yields up to a subsequence: $$u_n(t_0)\to u_c(t_0) \ \ \mbox{in} \ \ H^1\ \ \mbox{as} \ \ n\to +\infty.$$ Let then $u_c$ be the solution to \fref{nlsk:0} with initial data $u_c(t_0)$, then the $H^1$ continuity of the flow ensures: $\forall t\in[t_0,0)$, $$u_n(t)\to  u_c(t) \ \ \mbox{in} \ \ H^1$$ and $u_c$ admits a geometrical decomposition \fref{geomdecompuc} with $$\P_n(t)\to \P_c(t),$$ see \cite{MR1}, \cite{MM1}, \cite{KMR} for related statements. The $H^1$ bound  \fref{bounduc} and the estimate on the geometrical parameters \fref{estcrucailk} follow by passing to the limit in \fref{backwardspropagationn}, \fref{propparametersbisbis}. This implies in particular that $u_c$ blows up at $t=0$. Eventually, the conservation of the $L^2$ norm ensures: $$\|u_c\|_{L^2}=\lim_{n\to +\infty}\|u_n(t_n)\|_{L^2}=\|Q\|_{L^2}.$$ Similarly,
\bee
E(u_c(t)) & = &  \ds\frac{1}{2}\int |\n u_c(t,x)|^2 -\frac{1}{4}\int k(x)|u_c(t,x)|^4\\
& = & \frac{1}{2\l_c^2}\int \left|\n Q_{\P_c} - i\frac{b_cy}{2}Q_{\P_c}\right|^2 -\frac{1}{4\l_c^2}\int k(\l_c y+\alpha_c)|Q_{\P_c}|^4+O(\lambda_c)\\
& = & \frac{|yQ|_{L^2}^2}{8}\frac{b_c^2}{\lambda_c^2}-\frac{1}{8}\int\nabla^2k(0)(y,y)Q^4+o(1)\to E_0 \ \ \mbox{as} \ \ t\to 0
\eee
and hence the conservation of energy and the choice of $C_0$ \fref{constantcnot} ensure: $$E(u_c)=E_0.$$
Eventually, we derive from \fref{backwardspropagationn}, \fref{propparametersbisbis} the rough bound $$|(\tilde{\gamma_n})_s|\lesssim \lambda_n^2$$ which implies using \fref{propparametersbisbis}:
$$\left|\dt\left(\gamma_n+\frac{C_0^2}{t}\right)\right|=\frac{1}{\lambda_n^2}\left|(\gamma_n)_s-\frac{C_0^2\lambda_n^2}{t^2}\right|=\frac{1}{\lambda_n^2}\left|(\tgamma_n)_s-\left(\frac{C_0^2\lambda_n^2}{t^2}-1\right)\right|\lesssim 1.
$$ We then integrate in time using \fref{ex9} to conclude: $$\gamma_n(t)+\frac{C_0^2}{t}=\gamma_0+O(t)$$ and hence the limit in the phase parameter \fref{estcrucailkbis} follows by taking the limit $n\to+\infty$.\\

{\bf step 2} $H^{\frac{3}{2}}$ bound.\\

It remains to prove the $H^{\frac{3}{2}}$ bound \fref{esthtwo} which implies in particular the bound \fref{bounduc} by passing to the weak  $H^{\frac{3}{2}}$ limit.\\
$\tu_n$ satisfies:
$$
i\p_t\tu_n+\Delta\tu_n =-\psi_n-k(x)\ut_n|\ut_n|^2-F_n
$$
with 
$$
i\p_t\tilde{Q}_n+\Delta\tilde{Q}_n +k(x)\tilde{Q}_n|\tilde{Q}_n|^2=\psi_n, \ \ F_n=k(x)(\tilde{Q}_n+\ut_n)|\tilde{Q}_n+\ut_n|^2-k(x)\tilde{Q}_n|\tilde{Q}_n|^2-k(x)\ut_n|\ut_n|^2.
$$
Hence from standard Strichartz bounds and the smoothing effect of the linear Schr\"odinger flow, there holds using $\tilde{u}(t_n)=0$:
\bea
\label{un38}
\nonumber \|\n^{\frac{3}{2}}\tu_n\|_{L^\infty_{[t,t_n]}L^2} &\lesssim & \|\n^{\frac{3}{2}}\psi_n \|_{L^{4/3}_{[t,t_n]}L^{4/3}}+\|(1+|x|^2)F_n\|_{L^2_{[t,t_n]}H^{1}}\\
& + &\|\n^{\frac{3}{2}}(\ut_n|\ut_n|^2) \|_{L^{4/3}_{[t,t_n]}L^{4/3}}
\eea
The error $\psi_n$ is estimated from \fref{defpsin} which yields a bound: $$\|\n^{\frac{3}{2}}\psi_n\|_{L^{\frac{4}{3}}}\lesssim \frac{1}{\lambda_n ^3}(|Mod_n(t)|+|\alpha_n|^2+\lambda_n^5)\lesssim \l_n$$ with $Mod_n(t)$ corresponding to the modulation equation of $\mathcal P_n(t)$ like for \fref{law2}, and where we used from \fref{law2}, \fref{backwardspropagationn}, \eqref{propparametersbisbis}: $$Mod_n(t)\lesssim \lambda_n^4.$$ Hence: 
\be
\label{cgfgeir}
\|\n^{\frac{3}{2}}\psi_n \|_{L^{4/3}_{[t,t_n]}L^{4/3}}\lesssim \lambda_n^{\frac{7}{4}}.\ee
The $F_n$ term is a local term in $y$ which contains linear and quadratic terms in $\ut$, hence from \fref{backwardspropagationn}:
$$\|(1+|x|^2)F_n\|_{H^{1}}\lesssim \frac{1}{\lambda_n^3}\|\tu_n\|_{L^2}+\frac{1}{\lambda_n^2}\|\tu_n\|_{H^1}+ \frac{1}{\lambda_n^2}\|\tu_n\|^2_{L^4}+\frac{1}{\lambda_n}\|\tu_n\|_{L^4}\|\n\tu_n\|_{L^4}\lesssim \l_n+\l_n^{\frac{5}{2}}\|\n^{\frac{3}{2}}\tu_n\|_{L^2}$$ and thus: 
\be
\label{cnehfeofuye}
\|(1+|x|^2)F_n\|_{L^2_{[t,t_n]}H^{1}}\lesssim \lambda^{\frac{3}{2}}_n+\l_n^3\|\n^{\frac{3}{2}}\tu_n\|_{L^\infty_{[t,t_n]}L^2}.
\ee
Eventually, the nonlinear term is estimated from Sobolev embeddings and standard nonlinear estimates in Besov spaces:
$$\|\n^{\frac{3}{2}}(\ut_n|\ut_n|^2) \|_{L^{4/3}}\lesssim \|\n^{3/2}\ut_n\|_{L^2}\|\ut_n\|_{H^1}^2\lesssim \|\n^{3/2}\ut_n\|_{L^2}\l_n^6.$$

Injecting this together with \fref{cgfgeir}, \fref{cnehfeofuye} into \fref{un38} yields: 
$$
\|\n^{\frac{3}{2}}\tu_n\|_{L^\infty_{[t,t_n]}L^2}\lesssim \l^{\frac{3}{2}}_n+\l^3_n\|\n^{\frac{3}{2}}\tu_n\|_{L^\infty_{[t,t_n]}L^2}$$ and \fref{esthtwo} follows.\\
This concludes the proof of Proposition \ref{existenceprop}.


\section{Critical mass blow up solutions have conformal speed}
\label{sectionspeed}

We now turn to the proof of uniqueness of the critical mass blow up solutions at a nondegenerate point. We let $k$ satisfy {\bf Assumption (H)} and let $u(t)\in H^1$ be a solution to \fref{nlsk:0} with critical mass $$\|u\|_{L^2}=\|Q\|_{L^2}$$ and which blows up at $T=0$. We let $E_0=E(u_0)$.\\
We start with slightly revisiting the variational arguments of \cite{M2} to show that the solution admits near blow up time a decomposition:
$$u(t,x)=\frac{1}{[k(\alpha(t))]^{\frac{1}{2}}}\frac{1}{\lambda(t)}Q\left(t,\frac{x-\alpha(t)}{\lambda(t)}\right)e^{i\gamma(t)}+\tilde{u}$$ with the following a priori bounds:
\be
\label{energyestimate}
 \|\tilde{u}(t)\|_{H^1}\leq C(u_0), 
 \ee
 \be
 \label{nkohvoehe}
 \alpha(t)\to\alpha^* \ \ \mbox{as} \ \ t\to 0 \ \ \mbox{with} \ \ k(\alpha^*)=1,
 \ee
 \be
 \label{speed}
 \lambda(t)\leq C(u_0)|t| \ \ \mbox{ie} \ \ |\nabla u(t)|_{L^2}\geq \frac{C(u_0)}{|t|}.
 \ee
In particular, the solution  will be in the regime described in section \ref{refinedenergybounds}.\\
The main difficulty at this stage is that there holds {\it no a priori upper bound on the blow up rate}, and the remaining subsections are devoted to proving that $u$ must in fact blow up with the conformal speed: 
\be
\label{sharpblowupspeed}
\|\nabla u(t)\|_{L^2}\sim \frac{C(u_0)}{|t|}
\ee 
 which is a key step towards classification. The key will be to exhibit a new Lyapounov type rigidity property adapted to the critical mass setting -Proposition \ref{prop:localvirial}- which will allow us to improve the energy bound \fref{energyestimate} for a {\it dispersive} bound: 
$$ \|\tilde{u}(t)\|_{H^1}\to 0 \ \ \mbox{as} \ \ t\to 0.$$


\subsection{Variational estimates and convergence of the concentration point}


We start with revisiting the variational estimates in \cite{M2} in order to derive the convergence of the concentration point and the lower bound on the blow up speed. 

\begin{lemma}[Variational control of minimal mass blow up solutions]
\label{propinitialvar}
Let $u(t)$ be a critical mass solution to \fref{nlsk:0} which blows up at $T=0$. Then for $t<0$ close enough to 0, $u(t)$ admits a geometrical decomposition 
\be
\label{th1:decompnonline}
u(t,x)=\frac{1}{\lambda(t)}(Q+\e)\left(t,\frac{x-\alpha(t)}{\lambda(t)}\right)e^{i\gamma(t)},
\ee
for some $\mathcal C^1$ parameters $(\lambda(t), \alpha(t),\gamma(t))\in \RR^*_+\times \RR^2\times\RR$ with:\\
{\it 1. Uniform bound on the decomposition}:
\be\label{th1:cme8}
\begin{array}{l}
 |1-k(\a(t))|^{1/2}+\|\e(t)\|_{H^1(\RR^2)}\lesssim\l(t)\to 0 \ \ \mbox{as} \ \ t\to 0.
\end{array}
\ee
{\it 2. Convergence of the concentration point}:
\be
\label{convalpha}
\alpha(t)\to \alpha^* \ \ \mbox{with} \ \ k(\alpha^*)=1.
\ee
{\it 3. Lower bound on the blow up rate}: 
\be
\label{estlambdaprelim}
\lambda(t)\leq C(u_0)|t|.
\ee
\end{lemma}

{\it Proof of Lemma \ref{propinitialvar}}\\
 
The proof is standard and relies on the orbital stability of the ground state solitary wave. We will briefly sketch the argument and the computations which are a simplified version of the proof of Lemma \ref{estprelimestdecomp}.\\

{\bf step 1} Nonlinear decomposition of the flow.\\

Let $$v_0(t,x)=\lambda_0(t)u(t,\lambda_0(t)x) \ \ \mbox{with} \ \ \lambda_0(t)=\frac{|\nabla Q|_{L^2}}{|\nabla u(t)|_{L^2}}, $$ then $$|v_0(t)|_{L^2}=|u(t)|_{L^2}=|Q|_{L^2}, \ \ |\nabla v_0(t)|_{L^2}=|\nabla Q|_{L^2}$$ and the conservation of energy and Assumption (H) imply 
$$ \frac{1}{2}\int |\nabla v_0|^2-\frac{1}{4}\int |v_0|^4\leq \frac{1}{2}\int |\nabla v_0|^2-\frac{1}{4}\int k(\lambda_0(t)y)|v_0|^4=\lambda_0^2(t)E_0\to 0\ \ \mbox{as} \ \ t\to 0.$$ Hence from  a standard concentration compactness argument, see \cite{Lions} and Lemma 1 in \cite{MR1}, there exist $(x_0(t),\gamma_0(t))\in \RR^2\times\RR$ such that $$v_0(t,\cdot+x_0(t))e^{i\gamma_0(t)}\to Q \ \ \mbox{in} \ \  H^1\ \ \mbox{as} \ \ t\to 0.$$ Hence $u(t)$ admits near blow up time a decomposition 
\be
\label{decomporimitive}
u(t,x)=\frac{1}{\lambda(t)}(Q+\e)\left(t,\frac{x-\alpha(t)}{\lambda(t)}\right)e^{i\gamma(t)} 
\ee with 
\be
\label{cnncoenoc}
|1-k(\alpha(t))|+\|\e(t)\|_{H^1}\to 0\ \ \mbox{as} \ \ t\to 0.
\ee
Using the implicit function theorem, the uniqueness of the decomposition \fref{cnncoenoc} can be ensured through a suitable choice of orthogonality conditions\footnote{Note that the profile $Q$ corresponds to a very rough approximate solution to \fref{eqw}, but at this stage, we do not know yet that $\alpha$ is stabilizing around $\alpha^*$ and hence we cannot introduce the refined profiles $Q_{\mathcal P}$, $\mathcal P=(b,\lambda,\beta, \alpha)$ which make sense for $\a$ close to $\a^*$ only.}. We then set the orthogonality conditions on $\e$ to be:
\be\label{th1:ortho0}
(\e_1,|y|^2Q)=0, \ \ (\e_1,yQ)=0, \ \ (\e_2,\rho)=0
\ee
where we wrote \be\label{th1:reim}
\e=\e_1+i\e_2.
\ee
Let $v=Q+\e$, then $v$ satisfies an equation similar to \fref{th1:eqrenormlaizedv}:
$$i\partial_s v+\Delta v-v+k(\lambda(t)y+\alpha(t))v|v|^{2}=i\lsl \Lambda v+i\xsl \cdot \nabla v+\gts v
$$ 
with $\gts=\gamma_s-1.$ We next write down the equation satisfied by $\e$ and compute the geometrical parameters using the orthogonality conditions \fref{th1:ortho0}, see \cite{MR1} for very closely related computations or Appendix A for more details in the setting of Lemma \ref{estprelimestdecomp}. A simple computation left to the reader leads to the $\mathcal C^1$  regularity in time of the geometrical parameters together with the bound:
\be
\label{keyboundsrough}
\left|\lsl\right|+\left|\xsl\right|\lesssim |\e|_{H^1}+|1-k(\alpha(t))|+\lambda(t).
\ee

{\bf step 2} Expansion of the conservation laws.\\

We now expand the conservation laws in the $\e$ variables. For the $L^2$ norm, we have from the critical mass assumption and \fref{decomporimitive}:
$$
\int |u(t,x)|^2dx=\int |v(s,y)|^2dy=\int Q^2,$$ and thus:
\be\label{th1:cm5}
2(\e_1,Q)+\int |\e|^2=0.
\ee
For the energy, we have by rescaling:
\be\label{th1:ce2}
\begin{array}{l}
\lambda^2E_0=\frac{1}{2}\int |\n v(s,y)|^2dy -\frac{1}{4}\int k(\l y+\a)|v(s,y)|^4dy=\tilde{E}(v).
\end{array}
\ee
We compute on $Q$ which has zero energy:
\bee
\tilde{E}(Q)& = & \ds\frac{1}{2}\int |\n Q|^2 -\frac{1}{4}\int k(\l y+\a)Q^4\\
\nonumber & = & \frac{1}{2}\int |\n Q|^2-\frac{1}{4}\int Q^4+\frac{1}{4}\int(1-k(\alpha))Q^4+\frac{1}{4}\int(k(\alpha)-k(\lambda y+\alpha))Q^4\\
\nonumber & = & \frac{1}{4}\int(1-k(\alpha))Q^4 +O(\lambda^2)
\eee
where used in the last step the boundedness of $k's$ derivatives and the radiallity of $Q$.

We then expand  \fref{th1:ce2} and inject the conservation of $L^2$ norm of \fref{th1:cm5} to get after a classical computation:
\bea
\label{nkovoeheoge}
\nonumber 2\lambda^2E_0 & = &\frac{1}{2}\int(1-k(\alpha))Q^4+(L_+\e_1,\e_1)+(L_-\e_2,\e_2)\\
& + & O(\lambda^2+(\lambda+|1-k(\alpha)|)\|\e\|_{H^1}+\|\e\|_{H^1}^3).
\eea
The coercivity property of Lemma \ref{lemmacoerc} together with our choice of orthogonality conditions  \eqref{th1:ortho0} and the degeneracy inherited from \fref{th1:cm5}: $$ |(\e_1,Q)|^2=o(\|\e\|_{H^1}^2)$$ yields: $$(L_+\e_1,\e_1)+(L_-\e_2,\e_2)\geq \frac{c_0}{2}\|\e\|_{H^1}^2.$$ Injecting this into \fref{nkovoeheoge} together with the a priori smallness \fref{cnncoenoc} yields: $$(1-k(\alpha))+\|\e\|_{H^1}^2\lesssim \lambda^2,$$ and \fref{th1:cme8} follows.\\

{\bf step 3} Convergence on the concentration point and upper bound on the blow up rate.\\

We now prove \fref{convalpha}, \fref{estlambdaprelim} which are a straightforward consequence of \fref{th1:cme8} and \fref{keyboundsrough}. Indeed, they first imply $$|\alpha_t|=\frac{1}{\lambda}\left|\xsl\right|\lesssim \frac{\|\e\|_{H^1}+\lambda+|1-k(\alpha)|}{\lambda}\lesssim  1 \ \ \mbox{and thus} \ \ \int_{-1}^0|\alpha_t|dt<+\infty.$$ Hence $\alpha(t)\to \alpha^*$ as $t\to 0$ and $k(\alpha^*)=1$ from \fref{th1:cme8}. This proves \fref{convalpha}. Next, from \fref{th1:cme8} and \fref{keyboundsrough}, $$|\lambda_t|=\frac{1}{\lambda}\left|\lsl\right|\lesssim \frac{\|\e\|_{H^1}+\lambda+|1-k(\alpha)|}{\lambda}\lesssim  1  \ \ \mbox{and thus} \ \ \lambda(t)\lesssim |t|$$ by integration from $t$ to 0 and the blow up assumption at $t=0$, i.e.  $\lambda(0)=0$.\\
This concludes the proof of Lemma \ref{propinitialvar}.


\subsection{Strict lower bound on the energy}


From Lemma \ref{propinitialvar}, we know that the center of mass of $u$ must stabilize around a point $\alpha^*$ where $k$ reaches its maximum. Without loss of generality with respect to the assumptions of Theorem \ref{th2}, we may assume that $$\alpha^*=0, \ \ k(0)=1, \ \ \nabla k(0)=0, \ \ \nabla^2k(0)<0.$$  Hence the solution is from \fref{th1:cme8} on $[t_0,0)$ for $t_0$ close enough to blow up time in the regime described in section \ref{refinedenergybounds}, and we
 may now sharpen the decomposition \fref{decomporimitive} by introducing the modified $Q_{\P}$ profiles given by \fref{solapp2} which by construction are a small deformation of $Q$. We thus decompose the solution for $t$ close to enough to $0$ as 
$$u(t,x)=\frac{1}{[k(\alpha(t))]^{\frac{1}{2}}}\frac{1}{\lambda(t)}(Q_{\P(t)}+\e)\left(t,\frac{x-\alpha(t)}{\lambda(t)}\right)e^{i\gamma(t)}
$$
where $\e$ satisfies the orthogonality conditions \fref{ortho1}-\fref{ortho5}. Moreover, from \fref{estlambdaprelim}, the rescaled time \fref{rescaledtime} is a global time: $$s(t)=\int_{t_0}^t\frac{d\tau}{\lambda^2(\tau)}\to +\infty  \ \ \mbox{as} \ \ t\to0,$$  and thus the solution satisfies the bounds \fref{cme7}, \fref{law2} of Lemma \ref{estprelimestdecomp}: $\forall s\in [s_0,+\infty)$, 
$$
b^2+|\b|^2+|\a|^2+\|\e\|_{H^1(\RR^2)}^2\lesssim \l^2\left(E_0+\frac{1}{8}\int \nabla^2k(0)(y,y)Q^4\right)+O(\P^4),
$$
$$
Mod(t)  \lesssim  \P^4+\P^2\|\e\|_{L^2(\RR^2)}+\|\e\|^2_{L^2}+ \|\e\|_{H^1}^3+\left(b-\frac{\l}{C_0}\right)\P^3+\P(\a^2+\b^2).
$$
We now claim that this implies in particular the strict lower bound on the energy:

\begin{lemma}[Strict lower bound on the energy]
\label{lemmenergystrict}
There holds:
\be
\label{boundinsideezero}
E_0>-\frac{1}{8}\int \nabla^2k(0)(y,y)Q^4.
\ee
\end{lemma}

{\bf Proof of Lemma \ref{lemmenergystrict}}\\

By contradiction, if $E_0+\frac{1}{8}\int \nabla^2k(0)(y,y)Q^4\leq 0$, then \eqref{cme7} yields:
$$b^2+|\b|^2+|\a|^2+\|\e\|_{H^1(\RR^2)}^2\lesssim \l^4,
$$
which together with \fref{law2}, implies:
$$
\frac{\l_s}{\l}=O(\l^2).
$$
Dividing by $\l^2$, we get:
$$\label{cme11}
\frac{\l_s}{\l^3}=O(1) \ \ \mbox{ie} \ \ \frac{\l_t}{\l}=O(1) \ \ \mbox{and hence} \ \ 
\ln(\l(t))=O(1) \ \ \mbox{as} \ \ t\to 0.$$
But this contradicts the fact that $\l(t)\goto 0$ as $t\goto 0$, and \fref{boundinsideezero} is proved. This concludes the proof of Lemma \ref{lemmenergystrict}.


\subsection{The localized virial identity}


We now turn to the heart of the proof which is an improved {\it local bound} for $\e$ locally on the singularity. Our aim is to strictly gain on the global in space energy bound \fref{cme7}. Let us stress that the fact that there holds no a priori upper bound on the blow up rate means that one cannot rule out a priori a regime for which $$b\ll \lambda,$$ where we recall that $b\sim \lambda$ is the expected pseudo-conformal regime. This makes the refined energy identity \fref{crc6}  {\it useless} at this stage because of the presence of the $b$ factor in front of the quadratic term of the RHS of \fref{roughenegryindeitty} which may lead this gain to degenerate.\\
We claim that a Morawetz type computation in the spirit of the local virial estimate initiated by Martel and Merle for the study of the mass critical KdV \cite{MM1} allows one to obtain a space time bound for weighted norms centered on the singularity which improves the energy bound \fref{cme7}. Let us however insist onto the fact that a very specific algebra is at hand here which relies on the critical mass assumption and avoids both the standard troubles when running two dimensional Morawetz estimates for Schr\"odinger, and the delicate study of explicit quadratic forms like in \cite{MM1}, \cite{MR1}, \cite{FMR}, the only required coercivity property here being the energetic lower bound \fref{coerclinearenergy}.
 
\begin{proposition}[Local virial control]
\label{prop:localvirial}
Let $\phi$ be given by \fref{defphi}. There exist  universal constants $c,\underline{c}>0$ and a large enough constant $A>0$ such that for $t$ close enough to 0, we have:
\be\label{vi10}
\begin{array}{l}
\ds\left\{-\left(\frac{b}{\l}\right)\frac{|yQ|_{L^2}^2}{4}+\frac{1}{2\l}\Im\left(\int A\n\phi\left(\frac{y}{A}\right)\cdot\n\e\overline{\e}\right)\right\}_s\\
\ds\geq\frac{\underline{c}}{\l}\left\{|\a|^2+\int |\n\e|^2e^{-\frac{|y|}{\sqrt{A}}}+\int |\e|^2\right\}+O(\|\e\|_{H^1}^2+\P^3+|\b|^2).
\end{array}
\ee
\end{proposition}

The strength of the estimate \fref{vi10} is that the smaller the $\lambda$ -i.e.  the faster the blow up speed-, the better is the estimate. Indeed, the terms involved in the boundary term in time are uniformly bounded from \fref{cme7}:
\be\label{vi11}
\ds\left|-\left(\frac{b}{\l}\right)\frac{|yQ|_{L^2}^2}{4}+\frac{1}{2\l}\Im\left(\int A\n\phi\left(\frac{y}{A}\right)\n\e\overline{\e}\right)\right|\lesssim \frac{|b|}{\l}+\frac{\|\e\|_{H^1}^2}{\l}\lesssim 1,
\ee
while from \fref{cme7} again and the finite time blow assumption: 
\be\label{vi12}
\begin{array}{r}
\ds\int_{s_0}^{+\infty}(\|\e\|_{H^1}^2+\P^3+|\b|^2)ds\lesssim \int_{s_0}^{+\infty}\l^2(s)ds\lesssim \int_{t_0}^0dt\lesssim 1.
\end{array}
\ee
 Hence, integrating \eqref{vi10} between $s_0$ and $+\infty$, and using \eqref{vi11} \eqref{vi12} yields:
\be\label{vi13}
\int_{s_0}^{+\infty}\frac{1}{\l}\left\{|\a|^2+\|\e\|_{L^2}^2+\int|\nabla \e|^2e^{-\frac{|y|}{\sqrt{A}}}\right\}ds\lesssim 1.
\ee
On the other hand,
\be
\label{divergencelawlambda}
\int_{s_0}^{+\infty}\lambda(s)ds=\int_{t_0}^0\frac{d\tau}{\lambda(\tau)}=+\infty
\ee 
from \fref{estlambdaprelim} and hence \fref{vi13} is indeed a strict gain with respect to the energy bound \fref{cme7}.\\

{\bf Proof of Proposition \ref{prop:localvirial}}\\

The proof relies on an algebraic computation and in particular the specific structure of the quadratic terms in $\e$ appearing in the RHS of \fref{vi10}.\\

{\bf step 1} Computation of $\left(\frac{b}{\l}\right)_s$.\\

We have:
$$\left(\frac{b}{\l}\right)_s=\frac{b_s+b^2}{\l}-\frac{b}{\l}\left(\frac{\l_s}{\l}+b\right),
$$
which together with the law \eqref{lawb} of $b$ \footnote{which is a consequence of the critical mass $L^2$ conservation}, and \eqref{cme7}, \eqref{law2} yields:
\bee
& & \left(\frac{b}{\l}\right)_s\left(\frac{|y|^2}{4}Q,\L Q\right)= \frac{1}{\l}\left[-\frac{\nabla^2k(0)(\a,\a)}{2}\int Q^2\right]\\
& + & \frac{1}{\lambda}\left[\int |\e|^2 -(R_1(\e),\L\s)-(R_2(\e),\L\t)\right]\\
& + & O\left(\P^3+\P\|\e\|_{L^2}+\|\e\|_{H^1}^2+|\a|^2+|\b|^2\right).
\eee
Using the definition \eqref{defr1} of $R_1$ and \eqref{defr2} of $R_2$, we obtain:
\bea
\label{compbl3}
& & \left(\frac{b}{\l}\right)_s\left(\frac{|y|^2}{4}Q,\L Q\right)= \frac{1}{\lambda}\left[-\frac{\nabla^2k(0)(\a,\a)}{2}\int Q^2+\int |\e|^2\right .\\
\nonumber & - & \left . (3\s\e_1^2+2\t\e_1\e_2+\s\e_2^2,\L\s) -(3\t\e_2^2+2\s\e_1\e_2+\t\e_1^2,\L\t)\right])\\
\nonumber &+ & O\left(\P^3+\P\|\e\|_{L^2}+\|\e\|_{H^1}^2+|\a|^2+|\b|^2\right).
\eea

{\bf step 2} Computation of the localized virial identity.\\

We now pick a large enough number $A>0$ and compute the localized virial identity:
\be
\label{vi2}
\begin{array}{l}
\ds\left\{\frac{1}{\l}\Im\left(\int A\n\phi\left(\frac{y}{A}\right)\cdot\n\e\overline{\e}\right)\right\}_s=-\frac{\l_s}{\l^2}\Im\left(\int A\n\phi\left(\frac{y}{A}\right)\cdot\n\e\overline{\e}\right)\\
\ds +\frac{2}{\l}\Im\left(\int\left[\frac{1}{2}\Delta \phi\left(\frac{y}{A}\right)\e+A\n\phi\left(\frac{y}{A}\right)\cdot\n\e\right]\overline{\partial_s\e}\right).
\end{array}
\ee
Now from \fref{cme7}, \fref{law2}:
$$\frac{\l_s}{\l^2}=\frac{1}{\l}\left(\frac{\l_s}{\l}+b\right)-\frac{b}{\l}=O(1),
$$
which yields:
$$\ds\frac{\l_s}{\l^2}\Im\left(\int A\n\phi\left(\frac{y}{A}\right)\cdot\n\e\overline{\e}\right)=O(A\|\e\|^2_{H^1}).
$$
We may thus rewrite \eqref{vi2} as:
\bee
\label{vi5}
 \left\{\frac{1}{\l}\Im\left(\int A\n\phi\left(\frac{y}{A}\right)\cdot\n\e\overline{\e}\right)\right\}_s& = &   -\frac{2}{\l}\left(\int A\n\phi\left(\frac{y}{A}\right)\cdot(\p_s\e_1\n\e_2-\p_s\e_2\n\e_1)\right)\\
\nonumber &- & \frac{1}{\l}\left(\int \D\phi\left(\frac{y}{A}\right)(\p_s\e_1\e_2-\p_s\e_2\e_1)\right)+O(A\|\e\|^2_{H^1}).
\eee
We now inject the equation for $\e$  \eqref{eqe1b}-\eqref{eqe2b} and the estimates of Lemma \ref{estprelimestdecomp} to derive:
\bea
\label{vi6}
& & \left\{\frac{1}{2\l}\Im\left(\int A\n\phi\left(\frac{y}{A}\right)\n\e\overline{\e}\right)\right\}_s\\
\nonumber &=&\ds\frac{1}{\l}\left[\int A\n\phi\left(\frac{y}{A}\right)\cdot\left(M_2(\e)\n\e_2 +M_1(\e)\n\e_1\right)\right]-  \frac{1}{2\l}\left[\int \D\phi\left(\frac{y}{A}\right)\left(M_2(\e)\e_2+M_1(\e)\e_1\right)\right]\\
\nonumber&+ & O(A\|\e\|^2_{H^1}+\P^3+\P\|\e\|_{L^2}+|\a|^2+|\b|^2)\\
\nonumber& = &\frac{1}{\l}\left[\int \n^2\phi\left(\frac{y}{A}\right)(\n\e,\overline{\n\e})-\frac{1}{4A^2\l}\int \D^2\phi\left(\frac{y}{A}\right)|\e|^2\right]\\
\nonumber&+ & \frac{1}{2\l}\left[\int A\n\phi\left(\frac{y}{A}\right)(\n(Q^2+2\s^2)\e_1^2+4\n(\s\t)\e_1\e_2+\n(Q^2+2\t^2)\e_2^2)\right]\\
\nonumber&+ & O\left(A\|\e\|^2_{H^1}+\P^3+\P\|\e\|_{L^2}+|\a|^2+|\b|^2\right).
\eea

{\bf step 3} Conclusion.\\

Summing \eqref{compbl3} and \eqref{vi6} and expanding the nonlinear term according to \fref{defr1}, \fref{defr2}, we obtain after a few algebraic manipulations:
\bea
\label{neknvoeheoho}
\nonumber & & \left\{\left(\frac{b}{\l}\right)\left(\frac{|y|^2}{4}Q,\L Q\right)+\frac{1}{2\l}\Im\left(\int A\n\phi\left(\frac{y}{A}\right)\n\e\overline{\e}\right)\right\}_s =\frac{1}{\l}\left(-\frac{\nabla^2k(0)(\a,\a)}{2}\int Q^2\right)\\
\nonumber & + & \frac{1}{\lambda}\left[\int \n^2\phi\left(\frac{y}{A}\right)(\n\e,\overline{\n\e})+\int |\e|^2 -\frac{1}{4A^2}\int \D^2\phi\left(\frac{y}{A}\right)|\e|^2-\int 3Q^2\e_1^2+Q^2\e_2^2\right]\\
\nonumber &+ & \frac{1}{\lambda}\left[\int (A\n\phi\left(\frac{y}{A}\right)-y)\cdot(3Q\n Q \e_1^2+Q\n Q\e_2^2)\right]\\
 & + & O(A\|\e\|^2_{H^1(\RR)}+\P^3+\P\|\e\|_{L^2}+|\a|^2+|\b|^2).
\eea
Now, the choice of the function $\phi$ implies that:
\bea
\label{vi8}
\nonumber & & \int \n^2\phi\left(\frac{y}{A}\right)(\n\e,\overline{\n\e})+\int |\e|^2 -\frac{1}{4A^2}\int \D^2\phi\left(\frac{y}{A}\right)|\e|^2\\
& \geq & \int e^{-\frac{|y|}{2A}}|\nabla \e|^2+\int |\e|^2 +O\left(\frac{1}{A^2}\|\e\|_{L^2}^2\right).
\eea
Hence the quadratic form appearing in the RHS of \fref{neknvoeheoho} is a small deformation of the linearized energy, and a standard localization argument\footnote{see for example Appendix A in \cite{MM1}} together with the coercivity property \fref{coerclinearenergy}, the orthogonality conditions \eqref{ortho1}-\eqref{ortho5} and the conservation of mass \eqref{cm5} ensure the existence of a universal constant $\underline{c}>0$ such that for $A>0$ large enough: 
\be\label{vi9}
\begin{array}{r}
\ds\underline{c}\left[\int |\e|^2+\int |\n\e|^2e^{-\frac{|y|}{\sqrt{A}}}\right]\leq\int |\e|^2+\int |\n\e|^2e^{-\frac{|y|}{2A}}-\int (3Q^2\e_1^2+Q^2\e_2^2)\\
\ds +O(\|\e\|^4_{H^1(\RR^4)}+\P^4).
\end{array}
\ee
Together with \eqref{neknvoeheoho}, \eqref{vi8} and the bound \fref{cme7}, this yields the existence of constants $\underline{c},c>0$ such that for $A>0$ large enough and $t(A)<t<0$ close enough to $0$, we have:
$$
\begin{array}{l}
\ds\left\{-\left(\frac{b}{\l}\right)\frac{|yQ|_{L^2}^2}{4}+\frac{1}{2\l}\Im\left(\int A\n\phi\left(\frac{y}{A}\right)\cdot\n\e\overline{\e}\right)\right\}_s\\
\ds\geq\frac{\underline{c}}{\l}\left\{|\a|^2+\int |\e|^2+\int |\n\e|^2e^{-\frac{|y|}{\sqrt{A}}}\right\}+O(A\|\e\|_{H^1}^2+\P^3+|\b|^2),
\end{array}
$$
where we used the identity $(|y|^2Q,\L Q)=-|yQ|_{L^2}^2$. We may now fix $A$ once and for all and \fref{vi10} follows.\\
This concludes the proof of Proposition \ref{prop:localvirial}.


\subsection{Convergence to 0 of $\tilde{u}$ in $H^1$ away from the concentration point}


We now slighlty change our point of view and consider the decomposition of $u$ in original variables: $$u=\tq +\tu$$ with 
 \be\label{ie2}
\begin{array}{l}
\ds\tq(t,x)=\frac{1}{k(\a(t))^{1/2}\l(t)}\qp\left(\frac{x-\a(t)}{\l(t)}\right)e^{i\gamma(t)}, \\
\ds\tu(t,x)=\frac{1}{k(\a(t))^{1/2}\l(t)}\e\left(s,\frac{x-\a(t)}{\l(t)}\right)e^{i\gamma(t)}.
\end{array}
\ee
Then the bound \fref{cme7} and the definition \eqref{ie2} of $\tu$ yield:
\be\label{ie4}
\|\tu(t)\|_{L^2(\RR^2)}\lesssim\l(t)\textrm{ and }\|\tu(t)\|_{H^1(\RR^2)}\lesssim 1\textrm{ for all }t_0\leq t<0.
\ee
Our aim is to improve the energy bound \fref{ie4} for the dispersive bound
$$\tilde{u}(t)\to 0 \ \ \mbox{in} \ \ H^1 \ \ \mbox{as} \ \ t\to 0.
$$
The first step is dispersion away from the blow up point.

\begin{lemma}[$H^1$ dispersion away from the concentration point]
\label{lemma:convergence0H1}
There holds:
\be\label{ie26bis}
\lim_{t\goto 0}\|\tu(t)\|_{H^1(|x|\geq 3)}=0.
\ee
\end{lemma}

{\bf Proof of Lemma \ref{lemma:convergence0H1}}\\

Let $\tq,\tu$ be given by \fref{ie2}, then $\tu$ satisfies:
\be
\label{coeioeueiu}
i\partial_t\tu+\Delta\tu=-R-k(x)\tu|\tu|^2
\ee with 
\be
\label{defr}
R=k(x)\left[(\tq+\tu)|\tq+\tu|^2-\tq|\tq|^2-\tu|\tu|^2\right]+i\partial_t\tq+\Delta\tq+k(x)\tq|\tq|^2.
\ee

{\bf step 1} $L^2_tH^{3/2}_{loc}$ bound away from the concentration point.\\

We first claim as a consequence of the smoothing effect of the linear Schr\"odinger flow the space time bound:
\be
\label{gainh3halkf}
\int_{t_0}^0|\tu(\tau)|^2_{H^{\frac{3}{2}}(2\leq  r\leq 4)}d\tau<+\infty.
\ee
To wit, let the Fourier multiplier $D=(1-\Delta)^{1/2}$.
Let us start with treating the nonlinear term in \fref{coeioeueiu}. Let $\zt$ be the solution to 
\be
\label{eqztilde}
i\pa_t\zt+\Delta \zt=-k(x)\tilde{u}|\ut|^2, \ \ \zt(0)=0,
\ee 
then from standard Strichartz estimates and smoothing for the linear Schr\"odinger flow, there holds: 
\be
\label{cnekohoeri}
 \|\zt\|_{L^2((t_0,0),H^{\frac{3}{2}}_{loc})}+\|\zt\|_{L^{\infty}((t_0,0),H^1)}\lesssim \|D(k(x)\ut|\ut|^2)\|_{L^{\frac{4}{3}}((t_0,0),L^{\frac 43})}\lesssim 1
 \ee 
 where we used Sobolev embeddings and the $L^{\infty}((t_0,0),H^1)$ bound \fref{ie4}.\\
We then let $\wt=D^{\frac{1}{2}}(\tu-\zt)$ and $\chi$ be a radial smooth function with $$\chi'(r)=\int_0^r\chi''(\rho)d\rho \ \ \mbox{and} \ \ \chi''(r)=\left\{\begin{array}{lll} 
0\ \  \mbox{for} \ \ 0\leq r\leq 1,\\
1 \ \ \mbox{for} \ \ 2\leq r\leq 4,\\
					\frac{1}{r^2} \ \ \mbox{for} \ \ r\geq 5,
					\end{array} \right.
					$$
and such that 
\be
\label{estcutdeivatives}
\forall r\geq 0, \ \ \frac{|\chi'|^2}{r^2}\lesssim  \chi''(r).
\ee					
From \fref{cnekohoeri}, \fref{gainh3halkf} follows from: 
\be
\label{noekeouer}
\int_{t_0}^0\|\nabla \tilde{w}\|_{L^2(2\leq r\leq 4)}^2\lesssim 1.
\ee
Indeed, from \fref{coeioeueiu}, $\wt$ satisfies $$i\pa_t\wt+\Delta \wt=D^{\frac{1}{2}}R,$$ 
where $R$ is given by \eqref{defr}. We compute the associated localized virial identity:
\bea
\label{ie1}
& &  \frac{1}{2}\frac{d}{dt}\left\{\Im\left(\int\nabla\chi\cdot\nabla \wt\overline{\wt}\right)\right\} =    -\Re\left(\Delta\wt-D^{\frac{1}{2}}R,\overline{\frac{\Delta \chi}{2}\wt+\nabla \chi\cdot \nabla \wt}\right)\\
\nonumber & = & \int\chi''(\partial_r \wt)^2+\int\frac{\chi'}{r}(\partial_{\tau}\wt)^2-\frac{1}{4}\int\Delta^2\chi|\wt|^2+  \Re\left(D^{\frac{1}{2}}R,\overline{\frac{\Delta \chi}{2}\wt+\nabla \chi\cdot \nabla \wt}\right)
\eea
where $\partial_\tau\wt=\frac{1}{r}\partial_\theta\wt$. Let us now estimate the various terms in \fref{ie1}. The boundary term in time is bounded using \fref{ie4}, \fref{cnekohoeri}:
\be
\label{ie5}
\left|\Im\left(\int\nabla\chi\cdot\nabla \wt\overline{\wt}\right)\right|\lesssim \|\wt\|_{H^{\frac{1}{2}}}^2\lesssim \|\tu\|^2_{H^1(\RR^2)}+\|\zt\|^2_{H^1(\RR^2)}\lesssim 1.
\ee
Similarily, 
\be\label{ie6}
\left|\int \frac{1}{4}\Delta^2\chi|\wt|^2\right|\lesssim\|\tu\|^2_{H^{1/2}(\RR^2)}+\|\zt\|^2_{H^{1/2}(\RR^2)}\lesssim 1.
\ee
The $R$ terms in \fref{ie1} are treated from Cauchy Schwarz using the space localization of $\chi$:
\bea
\label{cncnoinitenm}
\nonumber \left|\Re\left(D^{\frac{1}{2}}R,\overline{\frac{\Delta \chi}{2}\wt+\nabla \chi\cdot \nabla \wt}\right)\right| & \lesssim & \|xD^{\frac{1}{2}}R\|_{L^2(r\geq 1)}\left[\left\|\frac{\Delta\chi}{r}\wt\right\|_{L^2(r\geq 1)}+\left\|\frac{\chi'}{r}\partial_r\wt\right\|_{L^2}\right]\\
\nonumber & \lesssim & \delta \int\chi''(\partial_r \wt)^2+\frac{1}{\delta}\int_{|x|\geq 1}|x|^2|D^{\frac{1}{2}}R|^2+\|\wt\|^2_{L^2}\\
& \lesssim & \delta \int\chi''(\partial_r \wt)^2+\frac{1}{\delta}\int|xD^{\frac{1}{2}}R|^2+1
\eea
for some small $\delta>0$, where we used \eqref{estcutdeivatives} and the $H^1$ bounds \fref{ie4}, \fref{cnekohoeri}.\\
We now claim 
\be
\label{estcbeiiehr}
\int_{t_0}^0\|xD^{\frac{1}{2}}R\|_{L^2}^2<+\infty
\ee
which together with \fref{ie5}, \fref{ie6}, \fref{cncnoinitenm} injected into \fref{ie1} yields \fref{noekeouer}.\\
{\it Proof of \fref{estcbeiiehr}}:  We use the usual symbolic calculus for the composition of pseudodifferential operators to obtain:
\be\label{fifi2}
D^{\frac{1}{2}}|x|^2D^{\frac{1}{2}}=xDx+a_{-1}(x,D)
\ee
where $a_{-1}(x,\xi)$ is a classical symbol of order $-1$. In particular, $a_{-1}(x,D)$ is bounded from $L^{\frac{4}{3}}$ to $L^4$. Together with \eqref{fifi2}, this yields the bound:
\be\label{fifi3}
\|xD^{\frac{1}{2}}R\|^2_{L^2}  =  (D^{\frac{1}{2}}|x|^2D^{\frac{1}{2}}R,R)\lesssim \|D^{\frac{1}{2}}xR\|^2_{L^2}+\|R\|^2_{L^{\frac{4}{3}}}.
\ee
We change variables: 
\be\label{boubou1}
R(t,x)=\frac{1}{\lambda^3}S\left(\frac{x-\alpha(t)}{\lambda(t)}\right)e^{i\gamma(t)}
\ee 
with 
\bea
\label{defsfinal}
\nonumber S(s,y)& = & i\partial_s \qp+\Delta \qp-\qp+\frac{k(\lambda(t)y+\alpha(t))}{k(\alpha(t))}\qp|\qp|^{2}-i\lsl \Lambda \qp\\
\nonumber & - & i\xsl \cdot\left(\nabla \qp+\frac{\lambda}{2}\frac{\nabla k(\alpha(t))}{k(\alpha(t))}\qp\right)+\gts \qp\\
& + & k(\lambda(t)y+\a(t))\left[(\qp+\e)|\qp+\e|^2-\qp|\qp|^2-\e|\e|^2\right]
\eea
which is well localized in $y$. Together with \eqref{fifi3}, this yields:
\be\label{fifi4}
\|xD^{\frac{1}{2}}R\|_{L^2}\lesssim 
\frac{1}{\l^{\frac{3}{2}}}\left(\|D^{\frac{1}{2}}yS\|_{L^2}+\|D^{\frac{1}{2}}S\|_{L^2}+\|S\|_{L^{\frac{4}{3}}}\right).
\ee
 We now explicitely expand the nonlinear terms in $\e$ in $S$ and use standard commutator estimates together with the good localization in space of $S(s,y)$, the bound on the geometrical parameters \fref{cme7}, \fref{law2}, and the $O(\lambda^3)$ control in the construction of $\qp$ \fref{eqq}, \fref{estpsitildebis} to conclude:
$$
\|xD^{\frac{1}{2}}R\|^2_{L^2}\lesssim 
\frac{1}{\l^{3}}\left(\int|\nabla \e|^2e^{-\frac{|y|}{\sqrt{A}}}+\|\e\|^2_{L^2}+\l^6+|\alpha|^4\right)\lesssim 1+\frac{1}{\l^{3}}\left(\int|\nabla \e|^2e^{-\frac{|y|}{\sqrt{A}}}+\|\e\|^2_{L^2}\right).
$$
\fref{estcbeiiehr} now follows from \fref{vi13}.\footnote{recall that $\frac{ds}{dt}=\frac{1}{\l^2}.$}\\

{\bf step 2} Strong $H^1$ convergence outside the blow up point.\\

The strong convergence \fref{ie26bis} now easily follows. Let a smooth cut off function $\tilde{\psi}$ with $\psit\equiv 1$ on $|x|\geq 3$ and $\psit\equiv 0$ on $|x|\leq 2$, then $w=\psit\tu$ satisfies the following equation:
\be
\label{ie14}i\p_tw+\Delta w=F -k(x)|w|^2w, \ \ w(t)\to 0 \ \ \mbox{in} \ \ L^2 \ \ \mbox{as} \ \ t\to 0,
\ee
with $F$ given by:
$$F=2\nabla\psit\cdot\nabla\tu+\Delta\psit \tu -k(x)\psit (1-\psit^2)|\tu|^2\tu-\psit R
$$ and $R$ given by \eqref{defr}. \eqref{ie26bis} follows from:
\be
\label{coeufeorufe}
\|w\|_{L^{\infty}((t_0,0),H^1)}\to 0 \ \ \mbox{as} \ \ t_0\to 0.
\ee 
Indeed, we write down Duhamel formula for \fref{ie14} and first derive a bound for $\nabla w$ in $L^4((t_0,0), L^4)$. Using the Strichartz estimates for the 2-dimensional Strichartz pair $(4,4)$ and the smoothing effect for the linear Schr\"odinger operator, we obtain:
\bea
\label{ie15bis}
 & & \|\nabla w\|_{L^4((t_0,0),L^4)}+\|\nabla w\|_{L^{\infty}((t_0,0),L^2)}\\
\nonumber & \lesssim & \|\nabla(k(x)|w|^2w)\|_{L^{4/3}((t_0,0)\times\RR^2)}+\ds\||x|F\|_{L^2((t_0,0),H^{\frac 12})}. 
\eea
The nonlinear term is estimated as follows:
\bea
\label{vbjeeoibis}
\nonumber & & \|\nabla(k(x)|w|^2w)\|_{L^{4/3}((t_0,0),L^{\frac{4}{3}})} \lesssim  \|\nabla w\|_{L^4((t_0,0),L^4)}\|\tu\|^2_{L^4((t_0,0),L^4)}\\
\nonumber &\lesssim &\|\nabla w\|_{L^4((t_0,0),L^4)}\|\tu\|_{L^\infty((t_0,0),H^1)}\|\tu\|_{L^\infty((t_0,0),L^2)}\lesssim \|\l\|_{L^\infty(t_0,0)}\|\nabla w\|_{L^4((t_0,0),L^4)}\\
&\lesssim &\d\|\nabla w\|_{L^4((t_0,0),L^4)}
\eea
for some small enough constant $\delta>0$ and where we used \eqref{ie4}. The second term is estimated using the compact support property of $\nabla\psit$, $\Delta\psit$ and $(1-\psit^2)$ and \fref{gainh3halkf}, \fref{estcbeiiehr}:
\bea
\label{ie16}
\nonumber 
\||x|F\|_{L^2_{[t_0,0]}H^{1/2}(\RR^2)} & \lesssim &  \||x|R\|_{L^2_{[t_0,0]}H^{\frac{1}{2}}(|x|\geq 1)} +\|\tu\|_{L^2([t_0,0),H^{\frac{3}{2}}(2\leq r\leq 4))}\\
& \lesssim & o(1) \ \ \mbox{as} \ \ t_0\to 0.
\eea
\fref{ie15bis}, \fref{vbjeeoibis} and \fref{ie16} yield:
\be
\label{cheoeiyoeyt}
\|\nabla w\|_{L^4((t_0,0)\times\RR^2)}+\|\nabla w\|_{L^{\infty}((t_0,0),L^2)}\lesssim o(1) \ \ \mbox{as} \ \ t_0\to 0,
\ee
and \fref{coeufeorufe} is proved.\\
This concludes the proof of Lemma \ref{lemma:convergence0H1}.


\subsection{Convergence to 0 in average of $\tilde{u}$ in $H^1$}


Our aim is now to propagate the $H^1$ convergence of $\tilde{u}$ away from the concentration point \fref{ie26bis} to the blow up up region as well. Here the key is the refined bound \fref{vi13} provided by the refined virial dispersive estimate \fref{vi10} and which implies :$$\liminf_{t\to 0}\frac{\|\tu(t)\|_{L^2}}{\lambda(t)} =0.$$ The first step is to obtain a convergence in average in time.

\begin{lemma}[$H^1$ dispersion in average in time]
\label{lemma:convergence0H1bis}
There holds:
\be\label{ie27}
\lim_{t\goto 0}\frac{1}{|t|}\int_{t}^0\left(\frac{1}{|\tau|}\int_{\tau}^0\|\tu\|^2_{H^1(\RR^2)}d\sigma\right)d\tau=0.
\ee
\end{lemma}

{\bf Proof of Lemma \ref{lemma:convergence0H1bis}}\\

Note that we may restrict the $H^1$ norm in \eqref{ie27} to the region $|x|\leq 3$ in view of \eqref{ie26bis}. \\

{\bf step 1} Morawetz identity.\\

We first claim the following virial type bound:
\be
\label{tobeprovbefhioere}
\int_{t}^0\int_{|x|\leq 3}|\nabla\tu|^2\lesssim o(|t|)+\|\tu(t)\|_{L^2}+\int_t^0\frac{|\alpha|^2+\|\e(\tau)\|_{L^2}^2+\int |\nabla\e|^2e^{-\frac{|y|}{\sqrt{A}}}}{\lambda^2(\tau)}d\tau.
\ee

Let $\chi$ a smooth radial cut off function on $\RR^2$ such that $\nabla^2\chi$ is supported in $|x|\leq 4$ and is positive semidefinite in $\RR^2$, and $\nabla\chi(x)=x$ on $|x|\leq 3$. 
We have:
\bea
\label{ie28}
\nonumber & & \dt\left\{\Im\left(\int\nabla\chi\cdot\nabla \tilde{u}\overline{\tilde{u}}\right)\right\}=-2\Re\left(\int(\Delta\tu+k(x)|\tu|^2\tu+R)\left[\overline{\frac{\Delta \chi}{2}\tu+\nabla\chi\cdot\nabla \tilde{u}}\right]\right)\\
\nonumber & = & 2\int \nabla^2\chi(\nabla \tu,\overline{\nabla \tu})-\frac{1}{2}\int \Delta^2\chi|\tu|^2-\frac{1}{2}\left(\int\Delta\chi k(x)|\tu|^4\right)+\frac{1}{2}\int\nabla\chi \cdot\nabla k(x)|\tu|^4\\
& - & 2\Re\int R\left[\overline{\frac{\Delta \chi}{2}\tu+\nabla\chi\cdot\nabla \tilde{u}}\right]
\eea
where $\tu$, $\tq$ and $R$ are defined as in \eqref{ie2} \eqref{defr}. We estimate the various terms in the right-hand side of \eqref{ie28}. Using \eqref{ie4}, we have:
\bea
\label{ie29}
\nonumber& &  \left|-\frac{1}{2}\int \Delta^2\chi|\tu|^2\ds -\frac{1}{2}\left(\int\Delta\chi k(x)|\tu|^4\right)+\frac{1}{2}\int\nabla\chi \cdot\nabla k(x)|\tu|^4\right|\lesssim \|\tu\|_{L^2(\RR^2)}^2+\|\tu\|^4_{L^4(\RR^2)}\\
& \lesssim & \|\tu\|_{L^2(\RR^2)}^2+\|\tu\|^2_{L^2(\RR^2)}\|\tu\|_{H^1(\RR^2)}^2\lesssim \l^2(t).
\eea
We now claim that we can treat the quadratic $R$ type terms perturbatively thanks to the vanishing of the cut at order 1 at the origin. Indeed, we change variables according to \fref{boubou1} with $S$ given by \fref{defsfinal} well localized in $y$. We thus get from the bound on the geometrical parameters \fref{law2}, $ \nabla\chi(x)\sim x$ near the origin and the $O(\lambda^3)$ control in the construction of $\qp$ \fref{eqq}, \fref{estpsitildebis}:
\bea
\label{vnoneoogjeohjreo}
\nonumber \left|\Re\int R\left[\overline{\frac{\Delta \chi}{2}\tu+\nabla\chi\cdot\nabla \tilde{u}}\right]\right|& = & \frac{1}{\lambda^2}\left|\Re\int S\left[\overline{\frac{\Delta \chi(\lambda(t)y+\alpha(t))}{2}\e+\frac{1}{\lambda}\nabla\chi(\lambda(t)y+\alpha(t))\cdot\nabla \e}\right]\right|\\
& \lesssim & \frac{1}{\lambda^2}\left[|\alpha|^2+\|\e\|_{L^2}^2+\int|\nabla \e|^2e^{-\frac{|y|}{\sqrt{A}}}+\lambda^3\right].
\eea
We inject \fref{ie29}, \fref{vnoneoogjeohjreo} into \fref{ie28} and integrate in time so that:
\bee
& & \int_{t}^0\int_{|x|\leq 3}|\nabla\tu(\tau)|^2d\tau \lesssim  \int_t^0\int \nabla^2\chi(\nabla \tu,\overline{\nabla \tu})\\
& \lesssim & \left[\Im\left(\int\nabla\chi\cdot\nabla \tilde{u}\overline{\tilde{u}}\right)\right]_t^0+\int_t^0\frac{|\alpha|^2+\|\e(\tau)\|^2_{L^2}+\int|\nabla \e(\tau)|^2e^{-\frac{|y|}{\sqrt{A}}}}{\lambda^2(\tau)}d\tau+\int_0^t\lambda(\tau)d\tau\\
& \lesssim & o(|t|)+\|\tu(t)\|_{L^2}+\int_t^0\frac{|\alpha|^2+\|\e(\tau)\|^2_{L^2}+\int|\nabla \e(\tau)|^2e^{-\frac{|y|}{\sqrt{A}}}}{\lambda^2(\tau)}d\tau
\eee
where we used from \fref{ie4}: $$\left|\Im\left(\int\nabla\chi\cdot\nabla \tilde{u}\overline{\tilde{u}}\right)\right|\lesssim \|\nabla \tu\|_{L^2}\|\tu\|_{L^2}\lesssim \|\ut\|_{L^2},$$ this is \fref{tobeprovbefhioere}.\\

{\bf step 2} Averaged in time dispersion.\\

We now divide \fref{tobeprovbefhioere} by $|t|$ and integrate in time. From the pointwise lower bound \fref{estlambdaprelim} and the dispersive bound \fref{vi13}\footnote{Recall that $\frac{ds}{dt}=\frac{1}{\lambda^2}$.}:
\bee
\frac{1}{|t|}\int_{t}^0\frac{|\alpha|^2+\|\e(\tau)\|^2_{L^2}+\int|\nabla \e(\tau)|^2e^{-\frac{|y|}{\sqrt{A}}}}{\lambda^2(\tau)}d\tau&\lesssim &\int_t^0\frac{|\alpha|^2+\|\e(\tau)\|_{L^2}^2+\int|\nabla \e(\tau)|^2e^{-\frac{|y|}{\sqrt{A}}}}{\lambda^3(\tau)}d\tau\\
& \to & 0\ \ \mbox{as} \ \ t\to 0.
\eee Similarily, from Cauchy Schwarz, 
$$\frac{1}{|t|}\int_t^0\frac{\|\tu(\tau)\|_{L^2}}{|\tau|}d\tau\lesssim \left(\int_t^0\frac{\|\tu(\tau)\|^2_{L^2}}{\tau^3}d\tau\right)^{\frac{1}{2}}\lesssim  \left(\int_t^0\frac{\|\tu(\tau)\|^2_{L^2}}{\lambda^3(\tau)}d\tau\right)^{\frac{1}{2}}\to 0 \ \mbox{as} \ \ t\to 0,$$ and \fref{ie27} follows.\\
This concludes the proof of Lemma \ref{lemma:convergence0H1bis}.


\subsection{Control of the modulation parameters}


We now claim that the $H^1$ dispersion \fref{ie27} coupled with the conservation laws implies a sharp control of the modulation parameters and hence the derivation of the blow up speed. Similarly like in \cite{DM1}, the key is to first derive convergence in a weak averaged in time sense, which thanks to the Lyapounov type dispersive control \fref{vi10} may then be turned into  {\it pointwise} controls on the whole sequence in time.

\begin{proposition}[Pointwise dispersive bounds]
\label{prop:improvecontrol}
There holds the pointwise bounds:
\be\label{bis:ex17}
|\b|+|\a|+\|\e\|_{H^1(\RR^2)}\lesssim \l^{2},
\ee 
\be\label{bis:ex19:7}
\begin{array}{ll}
\ds\left|\frac{b(s)}{\l(s)}-\frac{1}{C_0}\right| &\ds \lesssim \l^2,
\end{array}
\ee
\be\label{bis:ex19:9}
\l(t)=-\frac{t}{C_0}+O(|t|^3)\textrm{ for }t_0\leq t\leq 0.
\ee
Moreover, there exists $\gamma_0\in \RR$ such that: 
\be
\label{convergencephaseparameter}
\gamma(t)=-\frac{C_0^2}{t}+\gamma_0+O(|t|).
\ee
\end{proposition}

\begin{remark} Note that \fref{bis:ex17} implies in particular the zero momentum limit \fref{mommentgotozero}, and \fref{bis:ex17} improves the energy bound \fref{ie4} by a whole factor $\lambda$.
\end{remark}

{\bf Proof of Proposition \ref{prop:improvecontrol}}\\

{\bf step 1} Control in average of $\alpha$ and $\beta$.\\

We claim:
\be\label{ie51}
\lim_{t\goto 0}\ds\frac{1}{|t|}\int_t^0\left(\frac{1}{|\tau|}\int_\tau^0\left(\frac{|\a|^2}{\l^2}+\frac{|\b|^2}{\l^2}\right)d\sigma\right)d\tau=0.
\ee
We start with $\a$. Indeed, from \eqref{vi13} and the pointwise bound \fref{estlambdaprelim}:
$$\frac{1}{|t|}\int_t^0\frac{|\a|^2}{\l^2}d\tau\lesssim \int_t^0\frac{|\a|^2}{\l^3}d\tau= \int_s^{+\infty}\frac{|\a|^2}{\l}d\sigma\to 0 \ \ \mbox{as} \ \ t\to 0$$
which in particular implies \fref{ie51} for $\alpha$.\\
We now consider $\b$. We have:
\bea
\label{ie55}
\left(\frac{\a\cdot\b}{\l}\right)_s & = & \ds\a_s\cdot\frac{\b}{\l}+\b_s\cdot\frac{\a}{\l}-\frac{\a\cdot\b}{\l}\frac{\l_s}{\l}\\
\nonumber & = & \ds 2|\b|^2+(-b\b+c_0(\a)\l)\cdot\frac{\a}{\l}+b\frac{\a\cdot\b}{\l}\\
\nonumber & + & \left(\frac{\a_s}{\l}-2\b\right)\cdot\b+(\b_s+b\b-c_0(\a)\l)\cdot\frac{\a}{\l}-\left(\frac{\l_s}{\l}+b\right)\frac{\a\cdot\b}{\l}\\
\nonumber & = & \ds 2|\b|^2+c_0(\a)\cdot\a+O(\l^3+\|\e\|_{L^2}^2),
\eea
where we have used \eqref{law2} in the last equality. Integrating \eqref{ie55} between $s$ and $+\infty$, we obtain:
\bee
 2\int_s^{+\infty}|\b|^2d\sigma  &=& -\frac{\a(s)\cdot\b(s)}{\l(s)}-\int_s^{+\infty}c_0(\a)\cdot\a d\sigma+\int_s^{+\infty}O(\l^3+\|\e\|_{L^2}^2)d\sigma\\
& \lesssim &  |\a(s)|+\int_s^{+\infty}\left(|\a(\sigma)|^2+\|\e\|_{L^2}^2\right) d\sigma+\int_s^{+\infty}\l^3(\sigma) d\sigma
\eee
where we have used \eqref{cme7} which implies in particular that $\a\b/\l\goto 0$ when $s\goto +\infty$. Dividing by $|\tau|$ and using \fref{vi13} yields:
\bee
\nonumber  \frac{1}{|\tau|}\int_\tau^0\frac{|\b|^2}{\l^2}d\sigma & \lesssim & \frac{|\a(\tau)|}{|\tau|}+\frac{1}{|\tau|}\int_\tau^0\frac{d\sigma}{\l(\sigma)^2} \left[|\a(\sigma)|^2+\|\e\|_{L^2}^2\right]+\frac{1}{\tau}\int_\tau^0\l(\sigma) d\sigma\\
\nonumber & \lesssim & \frac{|\a(\tau)|}{|\tau|}+\int_\tau^0\frac{d\sigma}{\l(\sigma)^3} \left[|\a(\sigma)|^2+\|\e\|_{L^2}^2\right]+o(1)\\
& \lesssim & \frac{|\a(\tau)|}{|\tau|}+o(1).
\eee
Integrating once more and using Cauchy-Schwarz, this yields:
$$\frac{1}{|t|}\int_t^0\left(\frac{1}{|\tau|}\int_{\tau}^0\frac{|\b|^2}{\l^2}d\sigma\right)\lesssim  \frac{1}{|t|}\int_t^0\frac{|\a(\tau)|}{|\tau|}d\tau+o(1) \lesssim  \left(\int_t^0\frac{|\a(\tau)|^2}{\lambda^3(\tau)}d\tau\right)^{\frac{1}{2}}+o(1)=o(1)
$$
and \fref{ie51} follows for $\beta$.\\

{\bf step 2} Limit of $\frac{b}{\lambda}$ on a subsequence as $t\goto 0$.\\

\fref{ie27} and \fref{ie51} yield:
$$
\ds \lim_{t\goto 0}\frac{1}{|t|}\int_t^0\left(\frac{1}{|\tau|}\int_\tau^0\left(\frac{|\b|^2}{\l^2}+\frac{|\a|^2}{\l^2}+\frac{\|\e\|^2_{H^1(\RR^2)}}{\l^2}\right)d\sigma\right)d\tau =0.
$$
In particular, this implies the existence of a sequence $t_n\to 0$ such that:
\be\label{ie63}
\ds \lim_{n\goto +\infty}\frac{|\b(t_n)|}{\l(t_n)}+\frac{|\a(t_n)|}{\l(t_n)}+\frac{\|\e(t_n)\|_{H^1(\RR^2)}}{\l(t_n)}=0.
\ee
Injecting this into {\it the conservation of the energy} \eqref{cme5} yields:
\be\label{ie64}
\ds \lim_{n\goto +\infty}\left(\frac{b(t_n)}{\l(t_n)}\right)^2=\frac{E_0+\frac{1}{8}\int \nabla^2k(0)(y,y)Q^4}{\frac{1}{8}\int |y|^2Q^2}=\frac{1}{C_0^2}.
\ee
We now observe from \fref{compbl3} the bound:
$$\left|\left(\frac{b}{\lambda}\right)_s\right|\lesssim \frac{1}{\lambda}\left[\int|\e|^2+|\alpha|^2\right]+O(\lambda^2)$$ and thus \fref{vi13} ensures $$\int_s^{+\infty}\left|\left(\frac{b}{\lambda}\right)_s\right|ds<+\infty>$$ Hence $\frac{b}{\lambda}$ has a limit as $s\to +\infty$, and from \fref{ie64}: $$\frac{b}{\lambda}\to \pm \frac{1}{C_0} \ \ \mbox{as}\ \ t\to T.$$ Now from \fref{cme7}, \fref{law2}, 
\be
\label{cerudestlamnda}
\left|\lambda_t+\frac{b}{\lambda}\right|=\frac{1}{\lambda}\left|\lsl+b\right|\lesssim \lambda
\ee
and hence the finite time blow up assumption together with $\lambda(t)>0$ imply:\be\label{ie65}
\ds \lim_{t\to 0}\frac{b(t)}{\l(t)}=\frac{1}{C_0}>0.
\ee
Injecting this into \fref{cerudestlamnda} and integrating in time yields in particular the pseudo conformal speed: 
\be
\label{pseudoconflambda}
\lambda(t)=\frac{|t|}{C_0}(1+o(1)) \ \ \mbox{as} \ \ t\to 0.
\ee

{\bf step 3} Improved bounds.\\

We now claim that the knowledge of the asymptotic limit \fref{ie65} coupled with the monotonicity \eqref{vi10} and the conservation of energy implies a spectacular improvement on the bounds of $\|\e\|_{H^1}$.\\
Indeed, we integrate the local virial identity \eqref{vi10}  between $t$ and $t_n$ and let $t_n\to 0$. The boundary term in $t_n$ is estimated using \fref{cme7}, \eqref{ie65}:
\be
\label{limboundary}
-\left(\frac{b(t_n)}{\l(t_n)}\right)\int\frac{|y|^2}{4}Q^2+\frac{1}{2\l}\Im\left(\int A\n\phi\left(\frac{y}{A}\right)\n\e(t_n)\overline{\e}(t_n)\right)\to -\frac{1}{C_0}
\ee  as $t_n\to 0$, and we thus get: $\forall s>0$,
$$\begin{array}{r}
\ds \left(\frac{b(s)}{\l(s)}-\frac{1}{C_0}\right)\int\frac{|y|^2}{4}Q^2-\frac{1}{2\l(s)}\Im\left(\int A\n\phi\left(\frac{y}{A}\right)\n\e(s)\overline{\e}(s)\right)\\
\ds\geq\int_s^{+\infty}\left(\frac{\underline{c}}{\l}\left\{|\a|^2+\|\e\|^2_{L^2}\right\}+O(\P^3+|\b|^2+\|\e\|_{H^1}^2)\right)d\sigma,
\end{array}
$$
which together with \eqref{cme7} implies:
\be\label{bis:ex11:1}
\begin{array}{l}
\ds\int_s^{+\infty}\frac{1}{\l}\left\{|\a|^2+\|\e\|_{L^2}^2\right\}d\sigma\\
\ds\lesssim \ds \left(\frac{b(s)}{\l(s)}-\frac{1}{C_0}\right)+\|\e\|_{H^1(\RR^2)}+\int_{s}^{+\infty}(\l^3+|\b|^2+\|\e\|_{H^1}^2)d\sigma.
\end{array}
\ee
We now recall \fref{cme5bis} which implies:
\be\label{bis:ex12}
 \ds |\b|^2+|\a|^2+\|\e\|_{H^1(\RR^2)}^2\lesssim\ds\frac{\l^2}{C_0^2}-b^2+\l^4.
\ee
Dividing by $\l^2$, we obtain:
\be\label{bis:ex13}
\begin{array}{l}
 \ds \frac{|\b|^2}{\l^2}+\frac{|\a|^2}{\l^2}+\frac{\|\e\|_{H^1(\RR^2)}^2}{\l^2}\lesssim\ds\frac{1}{C_0^2}-\frac{b^2}{\l^2}+\l^2.
\end{array}
\ee
Finally, multiplying \eqref{bis:ex11:1} by $\frac{b(s)}{\l(s)}+\frac{1}{C_0}$, adding to \eqref{bis:ex13}, and noticing that the terms $\frac{b^2}{\l^2}-\frac{1}{C_0^2}$ cancel each other yields:
\be\label{bis:ex14}
\begin{array}{l}
\ds \frac{|\b|^2}{\l^2}+\frac{|\a|^2}{\l^2}+\frac{\|\e\|_{H^1(\RR^2)}^2}{\l^2}+\int_s^{+\infty}\frac{1}{\l}\left\{|\a|^2+\|\e\|_{L^2}^2\right\}d\sigma\\
\ds\lesssim \ds \l^2+\int_{s}^{+\infty}(\l^3+|\b|^2+\|\e\|_{H^1}^2)d\sigma\lesssim \l^2+\int_{s}^{+\infty}\l^2d\sigma\lesssim |t|\lesssim \l
\end{array}
\ee
where we used \fref{cme7}, \fref{pseudoconflambda}. This yields the improved pointwise bound $$|\beta|^2+\|\e\|^2_{H^1}\lesssim \l^3$$ which reinjected into \fref{bis:ex14} yields: 
\bea
\label{cnovovhoe}
\nonumber & &  \frac{|\b|^2}{\l^2}+\frac{|\a|^2}{\l^2}+\frac{\|\e\|_{H^1(\RR^2)}^2}{\l^2}+\int_s^{+\infty}\frac{1}{\l}\left\{|\a|^2+\|\e\|_{L^2}^2\right\}d\sigma\\
\nonumber & \lesssim & \ds \l^2+\int_{s}^{+\infty}(\l^3+|\b|^2+\|\e\|_{H^1}^2)d\sigma\lesssim \l^2+\int_{s}^{+\infty}\l^3d\sigma\\
& \lesssim & \lambda^2+\int_t^0\lambda d\tau\lesssim \lambda^2
\eea
where we used \fref{pseudoconflambda} again. This concludes the proof of the improved bound \fref{bis:ex17}.\\ 
We now integrate the localized virial identity \eqref{vi10} between $s$ and $+\infty$ and estimate using \eqref{ie65}, \fref{cnovovhoe}:
\bee
\left|\frac{b(s)}{\l(s)}-\frac{1}{C_0}\right| & \lesssim & \ds\int_s^{+\infty}\frac{1}{\l}\left\{|\a|^2+\|\e\|_{L^2}^2\right\}d\sigma+\|\e\|_{H^1(\RR^2)}+\int_{s}^{+\infty}\l^3d\sigma\\
& \lesssim & \lambda^2
\eee
and \fref{bis:ex19:7} is proved. Next, from \eqref{law2}, \fref{bis:ex17}, \fref{bis:ex19:7}:
\be\label{bis:ex19}
\frac{(\l)_s}{\l}+b=O(\l^3) \ \ \mbox{and thus} \ \ (\l)_t+\frac{1}{C_0}=O(|t|^2)
\ee
which yields \fref{bis:ex19:9} by integration in time. This eventually implies using also \fref{law2}, \fref{bis:ex17}:
$$\left|\dt\left(\gamma+\frac{C_0^2}{t}\right)\right|=\frac{1}{\lambda^2}\left|\gamma_s-\frac{C_0^2\lambda^2}{t^2}\right|=\frac{1}{\lambda^2}\left|\tgamma_s-\left(\frac{C_0^2\lambda^2}{t^2}-1\right)\right|\lesssim 1$$ and \fref{convergencephaseparameter} follows.\\
This concludes the proof of Proposition \ref{prop:improvecontrol}.

\begin{remark}\label{ie:rmk}
In view of \eqref{bis:ex17} and \eqref{bis:ex19:7}, we have:
\be\label{ie:rmk1:bis}
\P(\a^2+\b^2)+\left(b-\frac{\l}{C_0}\right)\P^3\lesssim \P^5.
\ee
In particular, the approximate solution $\qp$ is of order five -see \fref{estpsitilde}- as announced in Remark \ref{rmk:solapp}. Also, in view of \eqref{bis:ex17}, we have:
\be\label{ie:rmk1}
\P^4+\P^2\|\e\|_{L^2}+\|\e\|^2_{L^2}+\|\e\|_{H^1}^3+\P(\a^2+\b^2)+\left(b-\frac{\l}{C_0}\right)\P^3\lesssim \l^4,
\ee
 which together with \eqref{law2} implies:
\be\label{ie:rmk2}
\ds |\alpha|^2+|\beta|^2+|Mod(t)|\lesssim\l^4.
\ee
\end{remark}


\section{Uniqueness}
\label{uniqueness}


We now have obtained the exact blow up speed of critical mass blow up solutions which corresponds to the crucial relation $b\sim \lambda$ and the dispersive behavior $$\ut\to 0\ \ \mbox{in} \ \ H^1 \ \ \mbox{as} \ \ t\to0.$$ Let now $u_c$ be the critical mass blow up solution given by Proposition \ref{existenceprop} which blows up at $T=0$ and $\alpha^*=0$, with energy $E_0$ and phase parameter $\gamma_0$ given by \fref{convergencephaseparameter}. We need to prove that $$u=u_c.$$ The proof proceeds in two steps. We will first show that the refined 
estimates of Proposition \ref{prop:improvecontrol} together with the backwards propagation of smallness of Lemma \ref{propbaczkardw} imply the strong $H^1$ convergence $$u_c-u\to 0\ \ \mbox{in} \ \ H^1\ \ \mbox{as} \ \ t\to0.$$ We then show that the a priori estimates obtained for $u-u_c$ are strong enough to treat {\it perturbatively} the growth induced by the null space of $L$ when linearizing the equation close to $u_c$ and running the energy estimates of Lemma \ref{lemma:timederivative}. Both these steps require having sufficient a priori decay estimates and in particular require the construction of an the approximate solution to at least the order $O(\lambda^5)$ for the error.


\subsection{$H^1$ convergence to the critical element}


We claim the following dispersive property which somehow collects all previous estimates on the solution and is the key to the proof of uniqueness.

\begin{lemma}[$H^1$ convergence to the critical element]
\label{lemmakey}
There holds the strong convergence at blow up time: 
\be
\label{strgoncongevrgence}
u-u_c\to 0\ \ \mbox{in} \ \ H^1.
\ee
More precisely,
\be
\label{keyinformation}
\|\nabla(u-u_c)\|_{L^2}+\frac{\|u-u_c\|_{L^2}}{|t|}\lesssim |t|^3 \ \ \mbox{as} \ \ t\to 0.
\ee
\end{lemma}

{\bf Proof of Lemma \ref{lemmakey}}\\

{\bf step 1} Backwards propagation of smallness and improved bounds on the solution.\\

We first claim that the bounds of Proposition \ref{prop:improvecontrol} coupled with Lemma \ref{propbaczkardw} imply a refined bound\footnote{Note that this bound is the same like for $u_c$, see \eqref{bounduc}}:
\be
 \label{backwarddd}
 \|\ut(t)\|_{L^2}\lesssim \l^4(t),\, \|\nabla \ut(t)\|_{L^2}\lesssim \l^3(t).
 \ee
 Indeed, we decompose $u$ according to the geometrical decomposition \fref{newdecomp}. We then let an increasing sequence of times such that $t_n\to 0$. Observe from the conservation of mass and Proposition \ref{prop:improvecontrol} that the assumptions \eqref{initltwonorm}-\eqref{aprioriinitialparamaters} of Lemma \ref{propbaczkardw} are satisfied at any $t_n$. In particular, there is a time $t_0<0$ such that we have \eqref{backwardspropagation} for any $t_0\leq t\leq t_n$:
$$
 \|\nabla \ut(t)\|_{L^2}^2+\frac{\|\ut(t)\|^2_{L^2}}{\lambda^2(t)} \lesssim  \|\nabla \ut(t_n)\|_{L^2}^2+\frac{\|\ut(t_n)\|^2_{L^2}}{\lambda^2(t_n)}+\lambda^6(t).
$$
We now let $n\to +\infty$ so that \eqref{bis:ex17} yields \fref{backwarddd}.\\

{\bf step 2} Comparison between the modulation parameters of $u$ and $u_c$.\\

Let $b, \l, \a, \b, \gamma$ denote the modulation parameters of $u$, and $b_c, \l_c, \a_c, \b_c, \gamma_c$ denote the modulation parameters of $u_c$. We claim:
\be\label{diff5}
|\gamma-\gamma_c|+\frac{|\a-\a_c|}{|t|}+\frac{|\l-\l_c|}{|t|}+|b-b_c|+|\b-\b_c|\lesssim |t|^4.
\ee
The proof of \eqref{diff5} is a consequence of the modulation equation and the improved a priori bound \fref{backwarddd}, and it is postponed to Appendix C.\\

{\bf step 3} Comparison between $u$ and $u_c$.\\

\fref{keyinformation} is now a simple consequence of \fref{backwarddd}, \fref{diff5}. Indeed, 
\bee
\nonumber  \frac{\|u-u_c\|_{L^2}}{|t|}+\|\nabla(u-u_c)\|_{L^2} & \lesssim&  \frac{\|\tu\|_{L^2}}{|t|}+\|\nabla\tu\|_{L^2}+\frac{\|\tu_c\|_{L^2}}{|t|}+\|\nabla \tu_c\|_{L^2}\\
& + & \frac{\|\tq-\tq_c\|_{L^2}}{|t|}+\|\nabla(\tq-\tq_c)\|_{L^2}.
\eee
Now \eqref{diff5} yields after a simple computation:
$$ \|\tq-\tq_c\|_{L^2}\lesssim |t|^4,\, \|\nabla(\tq-\tq_c)\|_{L^2}\lesssim |t|^3,
$$ and \eqref{bounduc}, \eqref{backwarddd} now yield \eqref{keyinformation}.\\
 This concludes the proof of Lemma \ref{lemmakey}.
 
 
 \subsection{Energy estimates for the flow near $u_c$}
 

Let us now decompose:
\be
\label{un50}
u=u_c+\ttu,  \ \ \ds\ttu(t,x)=\frac{1}{k(\a_c(t))^{1/2}\l_c(t)}\e\left(t,\frac{x-\a_c(t)}{\l_c(t)}\right)e^{i\gamma_c(t)}.
\ee
Here we do not impose modulation theory and orthogonality conditions on $\e$. Indeed, one should have in mind that we do not have uniform well localized bounds on the $\tilde{u}_{c}$ part of $u_c$ and that the $\pa_t\tilde{u}_c$ part is in particular  controlled only through the $H^{\frac 32}$ bound \fref{bounduc}. We however claim that the a priori estimate from \fref{bounduc}, \eqref{keyinformation}:
\be
\label{backwardttu}
 \|\ttu\|_{L^2}\lesssim \l^4,\, \|\ttu\|_{H^1}\lesssim \l^3
 \ee
 is enough to treat the instability generated by the null space of $L$ perturbatively.\\
 Let 
 \be\label{un70}
N(t):=\sup_{t<\tau<0}\left(\|\ttu(\tau)\|^2_{H^1}+\frac{\|\ttu(\tau)\|^2_{L^2}}{\l^2_c(\tau)}\right),
\ee 
 and 
 \bea
 \label{defst}
 \mbox{Scal}(t) & = & (\e_1,Q)^2+(\e_2,\Lambda Q)^2+(\e_1,|y|^2Q)^2+(\e_2,\rho)^2\\
 \nonumber & + & (\e_1,yQ)^2+(\e_2,\nabla Q)^2.
 \eea
 
 We first claim the following energy bound:
 
 \begin{lemma}
 \label{lemma:basicestimateun}
 There holds for $t$ close enough to 0: 
 \be
 \label{controlnscal}
 N(t)\lesssim \sup_{t\leq \tau<0}\frac{\mbox{Scal}(\tau)}{\lambda_c^2(\tau)}+\int_t^0\frac{\mbox{Scal}(\tau)}{\lambda_c^3(\tau)}d\tau.
 \ee
 \end{lemma}
 
 {\bf Proof of Lemma \ref{lemma:basicestimateun}}\\
 
 It is a consequence of the energy estimate \fref{crc6} together with the a priori bound \fref{bounduc}.\\
 
 {\bf step 1} Application of Lemma \ref{lemma:timederivative}.\\
 
 Let $\mathcal I(\ttu)$ be given by \fref{defI}, we claim that:
 \be
 \label{iegigi}
\frac{\|\ttu\|_{L^2}^2}{\lambda_c^3}+O\left(N(t)+\frac{\mbox{Scal}(t)}{\lambda_c^3}\right)\lesssim  \frac{d\mathcal I}{dt}.
\ee 
 Indeed, let us apply Lemma \ref{lemma:timederivative} with $w=u_c=(u_c)_1+i(u_c)_2$, then the bound \fref{aprioriboundwsw} holds from \fref{bounduc} and $\psi$ given by \fref{un2} is identically zero. Hence \fref{crc6} becomes:
 \bea
 \label{un55}
\nonumber & & \frac{d\mathcal I}{dt}= -\frac{1}{\l_c^2}\Im\left(\int k(x)u_c^2\overline{\ttu}^2\right)-\Re\left(\int k(x)\p_tu_c\overline{(2|\ttu|^2u_c+\ttu^2\overline{u_c})}\right)\\
 \nonumber & + & \frac{b_c}{\l_c^2}\Bigg\{\int \frac{|\ttu|^2}{\l_c^2}+\Re\left(\int\nabla^2\phi\left(\frac{x-\a_c}{A\l_c}\right)(\n\ttu,\overline{\n\ttu})\right) -\frac{1}{4A^2}\left(\int\Delta^2\phi\left(\frac{x-\a_c}{\l_c}\right)\frac{|\ttu|^2}{\l_c^2}\right)\\
\nonumber & + & \l_c\Re\left(\int A\nabla\phi\left(\frac{x-\a_c}{A\l_c}\right)k(x)(2|\ttu|^2u_c+\ttu^2\overline{u_c})\cdot\overline{\nabla u_c}\right)\Bigg\}\\
& + & O\left(\frac{\|\ttu\|^{2}_{L^2}}{\l_c^2}+\|\ttu\|^{2}_{H^1}\right).
\eea
We consider the first two terms in the right-hand side of \eqref{un55} and expand $u_c=\tq_c+\tu_c$:\footnote{Keep in mind that we do not have satisfactory well localized bounds in $\tilde{u}_c$, and the corresponding terms will be treated using the smallness \fref{bounduc}}
\bea
\label{un56}
& & -\frac{1}{\l_c^2}\Im\left(\int k(x)u_c^2\overline{\ttu}^2\right)-\Re\left(\int k(x)\p_tu_c\overline{(2|\ttu|^2u_c+\ttu^2\overline{u_c})}\right)\\
\nonumber & = & -\frac{1}{\l_c^2}\Im\left(\int k(x)\tq_c^2\overline{\ttu}^2\right)-\Re\left(\int k(x)\p_t\tq_c\overline{(2|\ttu|^2\tq_c+\ttu^2\overline{\tq_c})}\right)\\
\nonumber & & -\frac{1}{\l_c^2}\Im\left(\int k(x)(2\tu_c\tq_c+\tu_c^2)\overline{\ttu}^2\right)-\Re\left(\int k(x)\p_t\tq_c\overline{(2|\ttu|^2\tu_c+\ttu^2\overline{\tu_c})}\right)\\
\nonumber & - & \Re\left(\int k(x)\p_t\tu_c\overline{(2|\ttu|^2\tq_c+\ttu^2\overline{\tq_c})}\right)-\Re\left(\int k(x)\p_t\tu_c\overline{(2|\ttu|^2\tu_c+\ttu^2\overline{\tu_c})}\right).
\eea
Arguing like for the proof of \eqref{firtboundq}, we may rewrite the first two terms in the right-hand side of \eqref{un56} as:
\bea
\label{un57}
& & -\frac{1}{\l_c^2}\Im\left(\int k(x)\tq_c^2\overline{\ttu}^2\right)-\Re\left(\int k(x)\p_t\tq_c\overline{(2|\ttu|^2\tq_c+\ttu^2\overline{\tq_c})}\right)\\
\nonumber & = &  -\frac{b_c}{\l^2_c}\int k(x)((|\tq_c|^2+2\ts_c^2)\ttu_1^2+4\ts_c\tt_c\ttu_1\ttu_2+(|\tq_c|^2+2\tt_c^2)\ttu_2^2)\\
\nonumber &- & \frac{b_c}{\l_c}\Re\left(\int \left(\frac{x-\a_c}{\l_c}\right) k(x)(2|\ttu|^2\tq_c+\ttu^2\overline{\tq_c})\cdot\overline{\nabla \tq_c}\right)+  O\left(\frac{\|\ttu\|^{2}_{L^2(\RR^2)}}{\l^2_c}+\|\ttu\|^{2}_{H^1(\RR^2)}\right)
\eea
where $\tq_c=\ts_c+i\tt_c$. For the next two terms in the right-hand side of \eqref{un56}, we use Sobolev embeddings and \fref{bounduc}  to obtain:
\bea
\label{un58}
\nonumber & & \bigg|-\frac{1}{\l_c^2}\Im\left(\int k(x)(2\tu_c\tq_c+\tu_c^2)\overline{\ttu}^2\right)-\Re\left(\int k(x)\p_t\tq_c\overline{(2|\ttu|^2\tu_c+\ttu^2\overline{\tu_c})}\right)\bigg|\\
\nonumber & \lesssim & \frac{1}{\l_c^2}\|\tq_c\|_{L^\infty(\RR^2)}\|\tu_c\|_{L^2(\RR^2)}\|\ttu\|^2_{L^4(\RR^2)}
+\frac{1}{\l_c^2}\|\tu_c\|^2_{L^4(\RR^2)}\|\ttu\|^2_{L^4(\RR^2)}\\
& + & \|\p_t\tq_c\|_{L^\infty(\RR^2)}\|\tu_c\|_{L^2(\RR^2)}\|\ttu\|^2_{L^4(\RR^2)}\lesssim  \frac{\|\ttu\|^{2}_{L^2(\RR^2)}}{\l^2_c}+\|\ttu\|^{2}_{H^1}
\eea
where we used the bound\footnote{The worst term is generated by the phase $|(\gamma_c)_t|\lesssim \frac{1}{\lambda_c^2}$.} $$\|\pa_t\tq_c\|_{L^{\infty}}\lesssim \frac{1}{\lambda_c^3}.$$
The last two terms in the right-hand side of \eqref{un56} require using the equation satisfied by $\tu_c$:
\be\label{tintin}
i\p_t\tu_c=-\Delta\tu_c-k(x)(|u_c|^2u_c-|\tq_c|^2\tq_c)-\tilde{\psi}_c
\ee
where $\tilde{\psi}_c$ is defined by:
$$\tilde{\psi}_c=i\p_t\tq_c+\Delta\tq_c+k(x)|\tq_c|^2\tq_c.
$$
Expanding this term like for \fref{defpsin}, we have using \fref{law2}, \fref{bounduc}, \fref{estcrucailk}, \fref{estcrucailkbis}: 
$$|\tilde{\psi_c}|_{L^2}\lesssim \frac{Mod_c(t)+|\alpha_c|^2+\lambda_c^5}{\lambda_c^2}\lesssim \lambda_c^2.$$ Using this together with \eqref{tintin}, integration by parts, Sobolev embeddings and the $H^{\frac{3}{2}}$ bound \fref{bounduc} now yields:
\bea
\label{un58bis}
\nonumber & & \bigg|-\Re\left(\int k(x)\p_t\tu_c\overline{(2|\ttu|^2\tq_c+\ttu^2\overline{\tq_c})}\right)-\Re\left(\int k(x)\p_t\tu_c\overline{(2|\ttu|^2\tu_c+\ttu^2\overline{\tu_c})}\right)\bigg|\\
\nonumber & \lesssim & \|\tu_c\|_{H^{\frac{3}{2}}}\left[\|2|\ttu|^2\tq_c+\ttu^2\overline{\tq_c}\|_{H^{\frac{1}{2}}}+\|2|\ttu|^2\tu_c+\ttu^2\overline{\tu_c}\|_{H^{\frac{1}{2}}}\right]\\
\nonumber & + & \|k(x)(|u_c|^2u_c-|\tq_c|^2\tq_c)+\tilde{\psi}_c\|_{L^2}\left[\|\tq_c\|_{L^\infty(\RR^2)}\|\ttu\|^2_{L^4}+\|\tu_c\|_{L^6}\|\ttu\|^2_{L^6}\right]\\
&\lesssim & \frac{\|\ttu\|^{2}_{L^2(\RR^2)}}{\l^2_c}+\|\ttu\|^{2}_{H^1}.
\eea
We now consider the last term in the right-hand side of \eqref{un55} and compute:
\bea
\label{un60}
& &  \Re\left(\int A\nabla\phi\left(\frac{x-\a_c}{A\l_c}\right)k(x)(2|\ttu|^2u_c+\ttu^2\overline{u_c})\overline{\nabla u_c}\right)\\
 \nonumber &= & \Re\left(\int A\nabla\phi\left(\frac{x-\a_c}{A\l_c}\right)k(x)(2|\ttu|^2\tq_c+\ttu^2\overline{\tq_c})\overline{\nabla \tq_c}\right)+\mbox{Error}
 \eea
 with from Sobolev embeddings and \fref{bounduc}:
\bea
\label{un61}
\nonumber &  & |\mbox{Error}|\lesssim  \left|\Re\int A\nabla\phi\left(\frac{x-\a_c}{A\l_c}\right)k(x)(2|\ttu|^2\tu_c+\ttu^2\overline{\tu_c})\overline{\nabla \tq_c}\right|\\
\nonumber & +&\left|\Re\int A\nabla\phi\left(\frac{x-\a_c}{A\l_c}\right)k(x)(2|\ttu|^2\tq_c+\ttu^2\overline{\tq_c})\overline{\nabla \tu_c}\right|\\
\nonumber & +& \left|\Re\int A\nabla\phi\left(\frac{x-\a_c}{A\l_c}\right)k(x)(2|\ttu|^2\tu_c+\ttu^2\overline{\tu_c})\overline{\nabla\tu_c}\right|\\
\nonumber & \lesssim & (\|\n\tq_c\|_{L^\infty}\|\tu_c\|_{L^2}+\|\tq_c\|_{L^\infty}\|\n\tu_c\|_{L^2})\|\ttu\|^2_{L^4}\\
&+ & \|\tu_c\|_{L^6}\|\n\tu_c\|_{L^2}\|\ttu\|^2_{L^6}
 \lesssim  \frac{\|\ttu\|^{2}_{L^2}}{\l^2_c}+\|\ttu\|^{2}_{H^1}.
\eea
Collecting the estimated \eqref{un55}-\eqref{un61} yields:
\bee
\nonumber & & \frac{d\mathcal I}{dt}  =  \frac{b_c}{\l_c^2}\left[\int \frac{|\ttu|^2}{\l_c^2}+\Re\left(\int\nabla^2\phi(\frac{x-\a_c}{A\l_c})(\n\ttu,\overline{\n\ttu})\right)\right .\\
\nonumber& -& \left . \int k(x)((|\tq_c|^2+2\ts_c^2)\ttu_1^2+4\ts_c\tt_c\ttu_1\ttu_2+(|\tq_c|^2+2\tt_c^2)\ttu_2^2)-\frac{1}{4A^2}\int\Delta^2\phi\left(\frac{x-\a_c}{A\l_c}\right)\frac{|\ttu|^2}{\l_c^2}\right .\\
\nonumber&+ &\left . \l_c\Re\left(\int \left(A\nabla\phi\left(\frac{x-\a_c}{A\l_c}\right)-\left(\frac{x-\a_c}{\l_c}\right)\right)k(x)(2|\ttu|^2\tq_c+\ttu^2\overline{\tq_c})\cdot\overline{\nabla \tq_c}\right)\right]\\
& + & O\left(\frac{\|\ttu\|^{2}_{L^2}}{\l_c^2}+\|\ttu\|^{2}_{H^1}\right).
\eee
We now use the uniform proximity of $Q_{\mathcal P_c}$ to $Q$ and the coercitivity property \fref{coerclinearenergy} to conclude like for the proof of \fref{firtboundq}:
\bee
 \frac{b_c}{\lambda_c^4}\left[\int|\nabla \e|^2e^{-\frac{|y|}{\sqrt{A}}}+\int|\e|^2+O(\mbox{Scal}(t))\right]+O
\left(\frac{\|\ttu\|^{2}_{L^2}}{\l_c^2}+\|\ttu\|^{2}_{H^1(\RR^2)}\right)\lesssim \frac{d\mathcal I}{dt}
\eee
which implies \fref{iegigi}.\\

{\bf step 2} Coercivity of $\mathcal I$.\\

We now recall from \fref{defI} the formula:
\bee
\mathcal I(t) & = & \frac{1}{2}\int |\n\ttu|^2 +\half\int \frac{|\ttu|^2}{\l_c^2}-\frac{1}{4}\int k(x)|u_c+\ttu|^4+\frac{1}{4}\int k(x)|u_c|^4\\
\nonumber & + & \int k(x)|u_c|^2(u_c)_1\ttu_1+\int k(x)|u_c|^2(u_c)_2\ttu_2+\half\frac{b_c}{\l_c}\Im\left(\int A\nabla\phi\left(\frac{x-\a_c}{A\l_c}\right)\cdot\nabla\ttu\overline{\ttu}\right)
\eee
Expanding $u_c=\tq_c+\tu_c$ and arguing like for the proof of \fref{iegigi}, we get using \fref{backwardttu} the rough upper bound: 
\be
\label{boundarytem}
|\mathcal I|\lesssim \|\ttu(t)\|^2_{H^1}+\frac{\|\ttu(t)\|^2_{L^2}}{\l^2_c(t)}\to 0 \ \ \mbox{as} \ \ t\to 0.
\ee
 We now claim the lower bound:
\bea
\label{coervifityagain}
\nonumber \nonumber \mathcal I(t) & \geq & \frac{1}{2\lambda_c^2} \left[(L_+\e_1,\e_1)+(L_-\e_2,\e_2)+o\left(\|\e\|_{H^1}^2\right)\right] \\
& \geq & \frac{\underline{c}}{\lambda_c^2}\left[\int\|\e\|_{H^1}^2-\mbox{Scal}(t)\right].
\eea
The proof is very similar to the one of \fref{coebbeof} using also the control of interaction terms similar to \fref{un58}, \fref{un58bis}, \fref{un61}. The details are left to the reader.\\
Integrating \fref{iegigi} from $t$ to $0$ using the boundary condition \fref{boundarytem} and the lower bound \fref{coervifityagain} now yields \fref{controlnscal}.\\
This concludes the proof of Lemma \ref{lemma:basicestimateun}.


\subsection{Control of the scalar products and proof of Theorem \ref{th2}}


It now remains to control the possible growth of the scalar product terms in \fref{controlnscal}. We claim:

\begin{lemma}[A priori control of the null space]
\label{scalarproductterms}
There holds for $t$ close enough to $0$:
\be
\label{estscalarprouctmain}
\mbox{Scal}(t)\lesssim |t|^{\frac{1}{2}}|t|^2N(t).
\ee
\end{lemma}

Let us assume Lemma \ref{scalarproductterms} and conclude the proof of Theorem \ref{th2}.\\

{\bf Proof of Theorem \ref{th2}}\\

From \fref{controlnscal}, \fref{estscalarprouctmain} and the law $\lambda_c\sim |t|$, we have for $t$ close enough to $0$:
$$N(t)\lesssim |t|^{\frac{1}{2}}N(t)+\int_t^0  \frac{N(\tau)}{\sqrt{\tau}}d\tau\lesssim|t|^{\frac 1 2}N(t)$$ and hence $N(t)=0$ for $t$ small enough. From the definition \fref{un70} of $N$, this yields $u=u_c$  and concludes the proof of Theorem \ref{th2}.\\

{\bf Proof of Lemma \ref{scalarproductterms}}\\

The proof follows by deriving the null space close to $Q_{\mathcal P_c}$ to sufficiently high order and reintegrating the corresponding modulation equations from blow up time. The worst behavior is on the even terms where the modulation equations are a deformation of the ones for $L_+,L_-$, and roughly correspond to the system of ODE's: 
\be
\label{modequations}
(\e_2,\Lambda Q)_s=2(\e_1,Q), \ \ (\e_1,|y|^2Q)_s=-4(\e_2,\Lambda Q), \ \ (\e_2,\rho)_s=-(\e_1,|y|^2Q)
\ee 
with initial degeneracy provided by the $L^2$ norm conservation law and the a priori bound \fref{backwardttu}: $$|(\e_1,Q)|\lesssim \int|\e|^2\lesssim \|\e\|_{L^2}\lambda(t)\sqrt{N(t)}\lesssim \lambda^4(t) \lambda(t)\sqrt{N(t)}.$$ The control of the worst paramater (related to the phase) requires 
\bee
& & |(\e_2,\rho)(s)|\lesssim \int_s^{+\infty}ds_1\int_{s_1}^{+\infty}ds_2\int_{s_2}^{+\infty}|(\e_1,Q)(s_3)|ds_3\\
& \lesssim & \lambda(s)\sqrt{N(s)}\int_s^{+\infty}ds_1\int_{s_1}^{+\infty}ds_2\int_{s_2}^{+\infty}\lambda_c^4ds_3\lesssim \frac{\lambda(s)\sqrt{N(s)}}{s}.
\eee
This implies \fref{estscalarprouctmain} and explains why we needed a small enough estimate $\|\e\|_{L^2}\lesssim \l_c^4$ in \fref{keyinformation}\footnote{ In fact, any bound $\|\e\|_{L^2}\lesssim \l_c^{3+\eta}$ for some $\eta>0$ is enough}, \fref{backwardttu}.\\
Let us now implement  the above strategy which requires being careful with respect to polynomial losses, and in particular demands a high order approximation of the null space close to $\qp$ .\\

{\bf step 1} Approximate equation in conformal variables to the order $O(\l_c^4)$.\\

Let $v,w$ be defined by:
\be
\label{ps3}
u(t,x)=\frac{1}{[k(\alpha_c(t))]^{\frac{1}{2}}}\frac{1}{\lambda_c(t)}v\left(t,\frac{x-\alpha_c(t)}{\lambda_c(t)}\right)e^{i\gamma_c(t)}, \ \ w(s,y)=v(s,y)e^{i\frac{b_c|y|^2}{4}-i\beta_c\cdot y}
\ee
then $w$ satisfies the equation:
\bee
& &  i\partial_s w+\Delta w-w+((b_c)_s+b_c^2)\frac{|y|^2w}{4}\\
& &  -\left\{((\beta_c)_s+b_c\beta_c)\cdot y+i\lambda_c\beta_c\cdot \frac{\nabla k(\alpha_c(t))}{k(\alpha_c(t))}\right\}w+\frac{k(\lambda_c(t)y+\alpha_c(t))}{k(\alpha_c(t))}w|w|^{2}\\
& = & i\left(\frac{(\l_c)_s}{\l_c}+b_c\right) \left(\Lambda w +\left(-ib_c\frac{|y|^2}{4}+i\b_cy\right)w\right)+((\gamma_c)_s-|\b_c|^2)w\\
&+ & i\left(\frac{(\a_c)_s}{\l_c}-2\b_c\right) \cdot\left(\nabla w+\left(-ib_c\frac{y}{2}+i\b_c\right)w+\frac{\lambda_c}{2}\frac{\nabla k(\alpha_c(t))}{k(\alpha_c(t))}w\right).
\eee
Starting with $u_c$, we also define $v_c$ and $w_c(s,y)=v_c(s,y)e^{i\frac{b_c|y|^2}{4}-i\beta_c\cdot y}$. We let $u=u_c+\ttu$ and define: 
\be
\label{defetilde}
v=v_c+\e, \ \ w=w_c+\tilde{\e}, \ \ \mbox{i.e.} \ \ \et=\e e^{i\frac{b_c|y|^2}{4}-i\beta_c\cdot y}.
\ee
Since $u_c$ satisfies \eqref{nlsk:0}, $\et$ satisfies:
\bea
\label{ps7}
\nonumber & &  i\partial_s\tw+\Delta \tw-\tw+((b_c)_s+b_c^2)\frac{|y|^2\tw}{4}-\left\{((\beta_c)_s+b_c\beta_c)\cdot y+i\lambda_c\beta_c\cdot \frac{\nabla k(\alpha_c(t))}{k(\alpha_c(t))}\right\}\tw\\
\nonumber&+ & \frac{k(\lambda_c(t)y+\alpha_c(t))}{k(\alpha_c(t))}(w|w|^{2}-w_c|w_c|^2)=((\gamma_c)_s-|\b_c|^2)\tw\\
\nonumber&= &\ds i\left(\frac{(\l_c)_s}{\l_c}+b_c\right) \left(\Lambda \tw +\left(-ib_c\frac{|y|^2}{4}+i\b_cy\right)\tw\right)\\
& + &i\left(\frac{(\a_c)_s}{\l_c}-2\b_c\right) \cdot\left(\nabla \tw+\left(-ib_c\frac{y}{2}+i\b_c\right)\tw+\frac{\lambda_c}{2}\frac{\nabla k(\alpha_c(t))}{k(\alpha_c(t))}\tw\right).
\eea
\eqref{bis:ex17}-\eqref{bis:ex19:9}, \eqref{ie:rmk2} and \eqref{ps7} yield:
\bea
\label{pscl8}
\nonumber & &  i\partial_s\tw+\Delta \tw-\tw -(c_0(\a_c)\l_c+\b_3\l_c^3)\cdot y\tw+\frac{k(\lambda_c(t)y+\alpha_c(t))}{k(\alpha_c(t))}(w|w|^{2}-w_c|w_c|^2)\\
&  = &  O\left(\l_c^4(1+|y|^2)\tw+\l_c^4(1+|y|)\n \tw\right).
\eea
Let $\tw=\tw_1+i\tw_2$ and $w_c=(w_c)_1+i(w_c)_2$. We now expand the nonlinear term in \eqref{pscl8} as well as $w_c=P_{\mathcal P_c}+\tilde{\e}_c$ and the expansion of the approximate solution \fref{solapp} to derive the equation at order $O(\l_c^4)$:
\be
\label{eqioheo}
-i\pa_s\te+\mq(\te)+(c_0(\a_c)\l_c+\b_3\l_c^3)\cdot y\te=-\psi
\ee
where $\mq$ is the expansion of $M$ at order 4:
\be
\label{defnmsnkfo}
\mq(\te)=\mq_1(\te_1)+i\mq_2(\te_2),
\ee
\be\label{pscl11}
\begin{array}{ll}
\ds \mq_1(\tw)= & \ds -\D\tw_1+\tw_1-\bigg(3Q^2+6QT_2+\frac{3}{2}\nabla^2k(0)(y,y)\l_c^2Q^2+3\nabla^2k(0)(y,\a_c)\l_c Q^2\\
&\ds +6QT_3+\frac{1}{2}\nabla^3k(0)(y,y,y)\l_c^3Q^2\bigg)\tw_1-2QS_3\tw_2,
\end{array}
\ee
\be\label{pscl12}
\begin{array}{ll}
\ds \mq_2(\tw)= &\ds -\D\tw_2+\tw_2-\bigg(Q^2+2QT_2+\frac{1}{2}\nabla^2k(0)(y,y)\l_c^2Q^2+\nabla^2k(0)(y,\a_c)\l_c Q^2\\
&\ds +2QT_3+\frac{1}{6}\nabla^3k(0)(y,y,y)\l_c^3Q^2\bigg)\tw_2-2QS_3\tw_1,
\end{array}
\ee
and where the remainder $\psi$ satisfies:
\be\label{pscl13}
\psi=O\left(\l_c^4(1+|y|^2)\tw+\l_c^4(1+|y|)\n \tw+\tw_c\tw+\tw_c^2\tw+w_c\tw^2+\tw^3\right).
\ee

\vspace{0.2cm}

{\bf step 2} Approximate null space.\\ 

Let $f(s,y)=O(e^{-c|y|})$ be a smooth well localized slowly time dependent function, then \fref{eqioheo} yields:
\be
\label{feiuyieyfei}
\frac{d}{ds}\left\{\Im(\e,\overline{f})\right\}=-\Re(\e,\overline{\mq(f)-i\pa_sf+(c_0(\a_c)\l_c+\b_3\l_c^3)\cdot yf})+O((\psi,f))
\ee
with 
\bea
\label{cofeongorgh}
\nonumber |(\psi,f)| &\lesssim &\l_c^4\|\tw\|_{L^2(\RR^2)}+\|\tw_c\|_{L^2(\RR^2)}\|\tw\|_{L^2(\RR^2)} + \|\tw_c\|^2_{L^4(\RR^2)}\|\tw\|_{L^2(\RR^2)}\\
\nonumber &+ &\|\tw_c\|_{L^4(\RR^2)}\|\tw\|_{L^4(\RR^2)}\|\tw\|_{L^2(\RR^2)}+ \|\tw\|^2_{L^4(\RR^2)}\|\tw\|_{L^2(\RR^2)}\\
&\lesssim & \l_c^4\|\tw\|_{L^2(\RR^2)}.
\eea
We now claim that we can find some real valued smooth well localized functions\\
 $A_2(y), B_2(y), D_3(s,y), E_2(s,y), E_3(s,y), F_2(s,y)$ with 
 $$A_2,B_2=O(\l_c^2e^{-c|y|}), \ \ E_2,F_2=O(\l_c^2e^{-c|y|}), \ \ E_3,D_3=O(\l_c^3e^{-c|y|})$$
 such that we have the following approximate null space relations:
\be
\label{estunfin}
\mq(\nabla Q+\l_c^2A_2)-i\pa_s(\nabla Q+\l_c^2A_2)=-a_1\l_c^2yQ+O(\l_c^3e^{-|y|}),
\ee
\be
\label{estdeuxfin}
\mq(i(yQ+\l_c^2B_2))-i\pa_s\left[i(yQ+\l_c^2B_2))\right]=-2i(\nabla Q+\l_c^2A_2)+O(\l_c^3e^{-|y|}),
\ee
for some universal constant 
\be
\label{tobecomputed}
a_1>0,
\ee
and for the even part:
\bea
\label{esttroisfin}
\nonumber & & \mq(\Lambda P_{\mathcal P_c}+E_2+E_3+iD_3)-i\pa_s\left[\Lambda P_{\mathcal P_c}+E_2+E_3+iD_3\right]\\
\nonumber &+ &(c_0(\a_c)\l_c+\b_3\l_c^3)\cdot y\left[\Lambda P_{\mathcal P_c}+E_2+E_3+iD_3\right]\\
& = &  -2P_{\mathcal P_c}+\l_ca_2(\a_c)\cdot yQ+\l_c^3a_3\cdot yQ+O(\l_c^4e^{-|y|}),
\eea
for a linear map $a_2$ on $\RR^2$ and a vector $a_3$ in $\RR^2$, 
\bea
\label{estquatrefin}
\nonumber & & \mq(i(|y|^2 P_{\mathcal P_c}+F_2))-i\pa_s\left[i(|y|^2 P_{\mathcal P_c}+F_2)\right]\\
& = &  -4i(\Lambda P_{\mathcal P_c}+E_2)+O(\l_c^3e^{-|y|}),
\eea
and
\be\label{estcinqfin}
 \mq(\rho)= |y|^2Q+O(\l_c^2e^{-|y|}).
\ee
The algebraic derivation of \fref{estunfin}, \fref{estdeuxfin}, \fref{tobecomputed}, \fref{esttroisfin}, \fref{estquatrefin},  \fref{estcinqfin} is done in Appendix D.\\

{\bf step 3} Control of $(\e_1,yQ)$.\\

From \fref{cofeongorgh}, \fref{estunfin}, \fref{estdeuxfin}, let $A=\nabla Q+\l_c^2A_2$, $B=yQ+\l_c^2B_2$ and rewrite  
\fref{feiuyieyfei} with $f=A$ and then $f=iB$ to get:
\be\label{ps46}
\left\{\begin{array}{l}
(\tw_1,B)_s=-2(\tw_2,A)+O(\l_c^3\|\tw\|_{L^2(\RR^2)}),\\
(\tw_2,A)_s=a_1\l_c^2(\tw_1,yQ)+O(\l_c^3\|\tw\|_{L^2(\RR^2)}).
\end{array}\right.
\ee
Now, using the fact that $\frac{ds}{dt}=\frac{1}{\l^2}$ together with \eqref{bis:ex19:9} yields:
\be\label{mickey}
\l_c(s)=\frac{C_0}{s}+O\left(\frac{1}{s^2}\right).
\ee
\eqref{ps46}, \eqref{mickey} and the fact that $B=yQ+O(\l_c^2)$ imply:
$$\left\{\begin{array}{l}
\ds (\tw_1,B)_s=-2(\tw_2,A)+F_1,\\
\ds (\tw_2,A)_s=\frac{a_1 C_0^2}{s^2}(\tw_1,B)+F_2,
\end{array}\right. \ \ \mbox{with} \ \ F_1,F_2=O\left(\frac{1}{s^3}\|\tw\|_{L^2}\right).
$$
We conclude from Lemma \ref{ps50:0} with $\varsigma=a_1 C_0^2>0$ from \fref{tobecomputed} that 
$$
\begin{array}{l}
\ds |(\tw_1,B)|+|(\tw_2,A)|\ds\lesssim \ds\int_{s}^{+\infty}\frac{\ln(\sigma)}{\sigma^2}\|\tw\|_{L^2(\RR^2)}d\sigma. 
\end{array}
$$
Using the fact that $A=\n Q+O(s^{-2})$ and $B=yQ+O(s^{-2})$ yields:
$$\begin{array}{l}
\ds |(\tw_1,yQ)|+|(\tw_2,\n Q)|\ds\lesssim \ds\frac{\|\tw\|_{L^2(\RR^2)}}{s^2}+\int_{s}^{+\infty}\frac{\ln(\sigma)}{\sigma^2}\|\tw\|_{L^2(\RR^2)}d\sigma. 
\end{array}
$$
This yields together with \fref{defetilde}, \fref{estcrucailk}, \fref{estcrucailkbis} and the definition \fref{un70}:
\bee
 |(\e_1,yQ)(t)|+|(\e_2,\n Q)(t)| & \lesssim & |t|\|\ttu\|_{L^2(\RR^2)}+\int_{t}^{0}\ln(|\sigma|)\|\ttu\|_{L^2(\RR^2)}d\sigma\nonumber\\
 & \lesssim & |t|^{\frac{1}{2}}|t|\sqrt{N(t)}\label{poum0}
\eee
and \fref{estscalarprouctmain} is proved for the odd directions.\\

{\bf step 4} Control of $(\e_1,Q)$, $(\e_2,\Lambda Q)$, $(\e_1,|y|^2Q)$, $(\e_2,\rho)$.\\

We first have from the $L^2$ conservation law, the critical mass assumption and \fref{defetilde}:
$$\Re(\tw,\overline{w_c})=\Re(\e,\overline{v_c})=-\frac{1}{2}\int |\e|^2=O(\l_c^4|t|\sqrt{N(t)}).
$$
Hence from \fref{bounduc}:
\be
\label{cneoheog}
\Re(\tw,\overline{P_{\mathcal P_c}})=2\Re(\tw,\overline{w_c})+O(\l_c^4\|\e\|_{L^2})=O(\l_c^4|t|\sqrt{N(t)}).
\ee

Using \fref{cofeongorgh}, \fref{esttroisfin} and \fref{feiuyieyfei} with $f=\L P_{\mathcal P_c}+E_2+E_3+iD_3$, we have:
\be\label{poum1}
\begin{array}{l}
\left(\Im(\te,\overline{\L P_{\mathcal P_c}+E_2+E_3+iD_3})\right)_s\\
=2\Re(\te,\overline{P_{\mathcal P_c}})-\l_c a_2(\a_c)\cdot (\te_1,yQ)-\l_c^3 a_3\cdot (\te_1,yQ)+O(\l_c^4\|\tw\|_{L^2(\RR^2)}). 
\end{array}
\ee
\eqref{mickey} and \eqref{poum1} yield:
$$
\begin{array}{l}
\left|\Im(\te,\overline{\L P_{\mathcal P_c}+E_2+E_3+iD_3})\right|\\
\lesssim \ds\int_s^{+\infty}\left(|\Re(\te,\overline{P_{\mathcal P_c}})|+\frac{|(\te_1,yQ)|}{\sigma^3}+\frac{\|\tw\|_{L^2(\RR^2)}}{\sigma^4}\right)ds. 
\end{array}
$$
Together with \eqref{poum0}, \eqref{cneoheog} and the definition \eqref{un70} of $N(t)$, we obtain:
\be
\label{poum2bis}\left|\Im(\te,\overline{\L P_{\mathcal P_c}+E_2+E_3+iD_3})\right|\lesssim \ds |t|^{\frac{5}{2}}|t|\sqrt{N(t)}.
\ee
Using \fref{cofeongorgh}, \fref{estquatrefin} and \fref{feiuyieyfei} with $f=i|y|^2P_{\mathcal P_c}+iF_2$, we have:
\be\label{poum3}
-\left(\Re(\te,\overline{|y|^2P_{\mathcal P_c}+F_2})\right)_s
=4\Im(\te,\overline{\L P_{\mathcal P_c}+E_2})+O(\l_c^3\|\tw\|_{L^2(\RR^2)}). 
\ee
\eqref{mickey} and \eqref{poum3} yield:
\be\label{poum4}
\left|\Re(\te,\overline{|y|^2P_{\mathcal P_c}+F_2})\right|\\
\lesssim \ds\int_s^{+\infty}\left(|\Im(\te,\overline{\L P_{\mathcal P_c}+E_2})|+\frac{\|\tw\|_{L^2(\RR^2)}}{\sigma^3}\right)ds. 
\ee
Together with \eqref{poum2bis} and the definition \eqref{un70} of $N(t)$, we obtain:
\be\label{poum5}
\left|\Re(\te,\overline{|y|^2P_{\mathcal P_c}+F_2})\right|\lesssim \ds |t|^{\frac{3}{2}}|t|\sqrt{N(t)}.
\ee

Using \fref{cofeongorgh}, \fref{estcinqfin} and \fref{feiuyieyfei} with $f=\rho$, we have:
\be\label{poum6}
\left(\te_2,\rho\right)_s
=-(\te_1,|y|^2Q)+O(\l_c^2\|\tw\|_{L^2(\RR^2)}). 
\ee
\eqref{mickey} and \eqref{poum6} yield:
\be\label{poum7}
\left|(\te_2,\rho)\right|
\lesssim \ds\int_s^{+\infty}\left(|(\te_1,|y|^2Q)|+\frac{\|\tw\|_{L^2(\RR^2)}}{\sigma^2}\right)ds. 
\ee
Together with \eqref{poum5}, the fact that $P_{\mathcal P_c}=Q+O(t^2)$ and the definition \eqref{un70} of $N(t)$, we obtain:
\be\label{poum8}
\left|(\te_2,\rho)\right|\lesssim \ds |t|^{\frac{1}{2}}|t|\sqrt{N(t)}.
\ee

Finally, \eqref{cneoheog}, \eqref{poum2bis}, \eqref{poum5}, \eqref{poum8} together with the fact that $P_{\mathcal P_c}=Q+O(t^2)$ imply:
$$\ds |(\e_1,Q)|+|(\e_2,\L Q)|+|(\e_1,|y|^2Q)|+|(\e_2,\rho)|\lesssim |t|^{\frac{1}{2}}|t|\sqrt{N(t)}
$$
and \fref{estscalarprouctmain} is proved for the even directions.\\


\section*{Appendix A: Computation of the modulation parameters}


This Appendix is devoted to the algebraic computation of the modulation parameters in the setting of a decomposition \fref{newdecomp} satisfying the orthogonality conditions \fref{ortho1}-\fref{ortho5}. We in particular aim at deriving the estimates \eqref{law2}, \fref{law2:bis}. Similar computations can be found in \cite{MR4}, \cite{PR}. We will in fact prove a slightly more general result which will be used in Appendix C.\\

Let $Mod(t)$ be the vector of modulation equations given by \fref{defmod}, we claim that we can find a polynomial vector $\mathcal R$ with 
\be
\label{esttnsplu}
|R(\mathcal P)|\lesssim \P(|\alpha|^2+|\beta|^2)+|\mathcal P|^4+\left(b-\frac{\lambda}{C_0}\right)|\mathcal P|^3
\ee
such that:
\bea
\label{computationparamatersappenidxgeneralized}
& & Mod(t) =  \mathcal R(\mathcal P)\\
\nonumber &+ & O\left(\P^7+[\P^2+|Mod(t)|]\|\e\|_{L^2}+\|\e\|_{L^2}^2+\|\e\|_{H^1}^3+\left|\|u\|_{L^2}-\|Q\|_{L^2}\right|\right)
\eea
and this implies in particular \fref{law2}.\\

{\bf step 1} Inner products.\\

We compute the inner products needed to compute the law of the parameters from the $\qp$ equation \eqref{eqq} where $M_1,M_2$ are given by \fref{defm1}, \fref{defm2}:
\bea
\label{pscal1}
\nonumber & & -\left(-M_2(\e)+b\L\e_1-2\b\cdot\n\e_1,\n\t\right)+\left(M_1(\e)+b\L\e_2-2\b\cdot\n\e_2,\n\s\right)\\
& = & O(\P^2\|\e\|_{L^2}),
\eea
\bea
\label{pscal2}
\nonumber & & \left(-M_2(\e)+b\L\e_1-2\b\cdot\n\e_1,y\s\right)+\left(M_1(\e)+b\L\e_2-2\b\cdot\n\e_2,y\t\right)\\
& = & O(\P^2\|\e\|_{L^2}),
\eea
\bea
\label{pscal3}
\nonumber & & -\left(-M_2(\e)+b\L\e_1-2\b\cdot\n\e_1,\L\t\right)+\left(M_1(\e)+b\L\e_2-2\b\cdot\n\e_2,\L\s\right)\\
& = & -2\Re(\e,\overline{\qp})+O(\P^2\|\e\|_{L^2}),
\eea
\bea
\label{pscal4}
\nonumber & &  \left(-M_2(\e)+b\L\e_1-2\b\cdot\n\e_1,|y|^2\s\right)+\left(M_1(\e)+b\L\e_2-2\b\cdot\n\e_2,|y|^2\t\right)\\
& = & O(\P^2\|\e\|_{L^2}),
\eea
\bea
\label{pscal5}
\nonumber & &  -\left(-M_2(\e)+b\L\e_1-2\b\cdot\n\e_1,\rho_2\right)+\left(M_1(\e)+b\L\e_2-2\b\cdot\n\e_2,\rho_1\right)\\
& = & O(\P^2\|\e\|_{L^2}).
\eea
To prove \eqref{pscal1}-\eqref{pscal5}, remark using the definition \eqref{solapp2} of $\qp$ that:
\be\label{solapp3}
\qp=Qe^{-ib\frac{|y|^2}{4}+i\b\c y}+O(\P^2).
\ee
Thus, it suffices to prove \eqref{pscal1}-\eqref{pscal5} with $\qp$ replaced by $Qe^{-ib\frac{|y|^2}{4}+i\b\c y}$. This is done in \cite{PR} for \eqref{pscal1}-\eqref{pscal4}. Let us briefly sketch the proof for \eqref{pscal5}. From:
$$ M_1(\e)= -\D\e_1+\e_1-((Q^2+2\s^2)\e_1+2\s\t\e_2)+O(\P^2\e),
$$
$$\ds M_2(\e)= -\D\e_2+\e_2-((Q^2+2\t^2)\e_2+2\s\t\e_1)+O(\P^2\e),
$$
we have:
$$\begin{array}{ll}
&\ds -\left(-M_2(\e)+b\L\e_1-2\b\cdot\n\e_1,\rho_2\right)+\left(M_1(\e)+b\L\e_2-2\b\cdot\n\e_2,\rho_1\right)\\
\ds = &\ds\left(\e_1,-\Delta\rho_1+\rho_1-(Q^2+2\s^2)\rho_1-2\s\t\rho_2+b\L\rho_2-2\b\cdot\nabla\rho_2\right)\\
&\ds +\left(\e_2,-\Delta\rho_2+\rho_2-(Q^2+2\t^2)\rho_2-2\s\t\rho_1-b\L\rho_1 +2\b\cdot\n\rho_1\right)+O(\P^2\e).
\end{array}
$$Next:
$$\s=Q+O(\P^2),\,\t=\left(-b\frac{|y|^2}{4}+\b \cdot y\right)Q+O(\P^2),
$$$$
\rho_1=\rho+O(\P^2),\,\rho_2=\left(-b\frac{|y|^2}{4}+\b\cdot  y\right)\rho+O(\P^2),
$$from which:
$$\begin{array}{ll}
&\ds -\left(-M_2(\e)+b\L\e_1-2\b\cdot\n\e_1,\rho_2\right)+\left(M_1(\e)+b\L\e_2-2\b\cdot\n\e_2,\rho_1\right)\\
\ds = &\ds\left(\e_1,L_+(\rho)\right)+\Bigg(\e_2,\left(-b\frac{|y|^2}{4}+\b y\right)L_-(\rho)-2\left(-b\frac{y}{2}+\b \right)\cdot\n\rho-\rho(-b)\\
&\ds -2\left(-b\frac{|y|^2}{4}+\b y\right)Q^2\rho-b\L\rho+2\b\cdot\n\rho\Bigg)+O(\P^2\|\e\|_{L^2})\\
\ds = &\ds\left(\e_1,L_+\rho\right)+\left(\e_2,\left(-b\frac{|y|^2}{4}+\b \cdot y\right)L_+\rho\right)+O(\P^2\|\e\|_{L^2})\\
\ds = &\ds\left(\e_1,|y|^2Q\right)+\left(\e_2,\left(-b\frac{|y|^2}{4}+\b\cdot y\right)|y|^2Q\right)+O(\P^2\|\e\|_{L^2})\\
\ds = &\ds\left(\e_1,|y|^2\s\right)+\left(\e_2,|y|^2\t\right)+O(\P^2\|\e\|_{L^2}).
\end{array}
$$
The orthogonality condition \fref{ortho4} now yields \fref{pscal5}.\\

{\bf step 2} Simplification of the equations.\\

To compute the modulation equations driving the geometrical parameters, we first simplify the equations \eqref{eqe1} and \eqref{eqe2} using the explicit construction of $\qp$ as follows:
\be\label{eqe1b}
\begin{array}{l}
\ds (b_s+b^2)\p_b\s+(\b_s+b\b-c_0(\a)\l-\b_3\l^3-\b_4\l^4)\p_\b\s\\
\ds +\p_s\e_1-M_2(\e)+b\L\e_1-2\b\n\e_1\\
\ds =\left(\frac{\l_s}{\l}+b\right)\L\s+(\tgamma_s-|\b|^2)\t+\left(\frac{\a_s}{\l}-2\b\right)\n\s+\Im(\tpsi)-R_2(\e),
\end{array}
\ee
and 
\be\label{eqe2b}
\begin{array}{l}
\ds (b_s+b^2)\p_b\t+(\b_s+b\b-c_0(\a)\l-\b_3\l^3-\b_4\l^4)\p_\b\t\\
\ds +\p_s\e_2+M_1(\e)+b\L\e_2-2\b\n\e_2\\
\ds =\left(\frac{\l_s}{\l}+b\right)\L\t-(\tgamma_s-|\b|^2)\s+\left(\frac{\a_s}{\l}-2\b\right)\n\s -\Re(\tpsi)+R_1(\e),
\end{array}
\ee
where the remainder $\tpsi$ has he following form from  \fref{noiveauerror}, \fref{estpsitilde}: there exists $F_{\mathcal P}(y)$ and a matrix $S_{\mathcal \P}(y)$ polynomial in $\mathcal P$ and smooth and well localized in $y$ such that$$
\tpsi(y)  = F_{\mathcal P}(y)+\mathcal S_{\mathcal P}(y) Mod(t)+ O\left(\P^2+Mod(t)\right)\left[|\e(y)|+|(1+|y|)\nabla \e(y)|\right]+O(\mathcal P^7)e^{-c|y|}
$$
with
\be
\label{ests}
|\pa_y^{k}\mathcal S_{\mathcal P}(y)|\lesssim |\mathcal P|e^{-C_{k}|y|},
\ee 
\be\label{jepeteuncable}
\left|\pa_y^{k}F_{\mathcal P}(y)\right|\lesssim \left(\P(|\alpha|^2+|\beta|^2)+|\mathcal P|^5+\left(b-\frac{\lambda}{C_0}\right)|\mathcal P|^3\right)e^{-C_{k}|y|}.
\ee

{\bf step 3} The law of $b$.\\

We take the inner product of the equation \eqref{eqe1b} of $\e_1$ with $-\L\t$ and we sum it with the inner product of equation \eqref{eqe2b} of $\e_2$ with $\L\s$.   We obtain after integrating by parts:
\bee
& &  -(b_s+b^2)\left\{\left(\frac{|y|^2}{4}Q,\L Q\right)\right\}-2\Re(\e,\overline{\qp})\\
& = &(R_1(\e),\L\s)+(R_2(\e),\L\t)+ \mathcal R_1(\mathcal P)+\mathcal S_1(\mathcal P)Mod (t)\\
&+ & O\left(\P^2+|Mod(t)|\right)\|\e\|_{L^2}+O(\mathcal P^7)
\eee
with $\mathcal R_1$ polynomial satisfying \fref{esttnsplu}, $\mathcal S_1$ polynomial in $\mathcal P$ of order 1, and where we have used the orthogonality condition \eqref{ortho3} and the computation of the inner product \eqref{pscal3}. We now rewrite the conservation of mass \fref{degeneracyeq}:
\bea
\label{coefeou}
\nonumber 2\Re(\e,\overline{\qp}) = -\int |\e|^2+k(\alpha)\left[\int|u|^2-\int Q^2\right]-\left[\int|\qp|^2-k(\alpha)\int Q^2\right]\\
=  -\int |\e|^2+\frac{\nabla^2k(0)(\a,\a)}{2}\int Q^2+k(\alpha)\left[\int|u|^2-\int Q^2\right]+\mathcal R_2(\P)
\eea
where from \fref{compm4} $\mathcal R_2(\mathcal P)$ is polynomial and satisfies \fref{esttnsplu}. This yields using the definition \fref{defdzeroalpha}:
\bea
\label{lawb}
 \nonumber& - & (b_s+b^2)\left\{\left(\frac{|y|^2}{4}Q,\L Q\right)\right\}-\frac{\nabla^2k(0)(\a,\a)}{2}\int Q^2 =  \frac{\|yQ\|_{L^2}^2}{4}\left(b_s+b^2-d_0(\alpha,\alpha)\right)\\
 & = & -\int|\e|^2 + \mathcal R_3(\mathcal P) +\mathcal S_1(\mathcal P)Mod (t)+(R_1(\e),\L\s)+(R_2(\e),\L\t)\\
 &+ & O\left((\P^2+|Mod(t)|)\|\e\|_{L^2}+\left|\|u\|_{L^2}-\|Q\|_{L^2}\right|+\P^7\right)
\eea
with $\mathcal R_3$ polynomial satisfying \fref{esttnsplu}.\\

{\bf step 4} The law of $\l$.\\

We take the inner product of the equation \eqref{eqe1b} of $\e_1$ with $|y|^2\s$ and we sum it with the inner product of equation \eqref{eqe2b} of $\e_2$ with $|y|^2\t$.   We obtain after integrating by parts:
\bea
\label{lawlambda}
0 & = &\ds \left(\frac{\l_s}{\l}+b\right)\left\{(\L Q,|y|^2Q)\right\}+\mathcal R_4(\P)+\mathcal S_2(\mathcal P)Mod (t)\\
\nonumber &+ &O\left((\P^2+|Mod(t)|)\|\e\|_{L^2}+\mathcal P^7+\|\e\|_{L^2}^2+\|\e\|_{H^1}^3\right)
\eea
with $\mathcal R_4$ polynomial satisfying from \fref{jepeteuncable}:
\be
\label{chiofyeiyo}
|\mathcal{R}_4(\P)|\lesssim \P(|\a|^2+|\b|^2)+|\P|^5+\left(b-\frac{\l}{C_0}\right)|\P|^3,
\ee
$\mathcal S_2$ polynomial in $\mathcal P$ of order 1, and where we have used the orthogonality condition \eqref{ortho4} and the computation of the inner product \eqref{pscal4}.\\

{\bf step 5} The law of $\a$.\\

We take the inner product of the equation \eqref{eqe1b} of $\e_1$ with $y\s$ and we sum it with the inner product of equation \eqref{eqe2b} of $\e_2$ with $y\t$.   We obtain after integrating by parts:
\bea
\label{lawalpha}
 0 & =&\ds\left(\frac{\a_s}{\l}-2\b\right)\left\{(yQ,\n Q)\right\}+\mathcal R_5(\P)+\mathcal S_3(\mathcal P)Mod(t)\\
\nonumber &+ & O\left((\P^2+|Mod(t)|)\|\e\|_{L^2}+\mathcal P^7+\|\e\|_{L^2}^2+\|\e\|_{H^1}^3\right)
\eea
with $\mathcal R_5$ polynomial satisfying \fref{esttnsplu}, $\mathcal S_3$ polynomial in $\mathcal P$ of order 1, and where we have used the orthogonality condition \eqref{ortho2} and the computation of the inner product \eqref{pscal2}.\\

{\bf step 6} The law of $\b$.\\

We take the inner product of the equation \eqref{eqe1b} of $\e_1$ with $-\n\t$ and we sum it with the inner product of equation \eqref{eqe2b} of $\e_2$ with $\n\s$.   We obtain after integrating by parts:
\bea
\label{lawbeta}
\nonumber & &  -(\b_s+b\b-c_0(\a)\l-\b_3\l^3)\left\{(yQ,\n Q)\right\}=\mathcal R_6(\P)++\mathcal S_4(\mathcal P)Mod(t)\\
 &+ & O\left((\P^2+|Mod(t)|)\|\e\|_{L^2}+\mathcal P^7+\|\e\|_{L^2}^2+\|\e\|_{H^1}^3\right)
\eea
with $\mathcal R_6$ polynomial satisfying \fref{esttnsplu}, $\mathcal S_4$ polynomial in $\mathcal P$ of order 1, and where we have used the orthogonality condition \eqref{ortho1} and the computation of the inner product \eqref{pscal1}.\\

{\bf step 7} The law of $\tgamma$.\\

We take the inner product of the equation \eqref{eqe2b} of $\e_2$ with $\rho$. We obtain after integrating by parts:
\bea
\label{lawtgamma}
\nonumber && -(b_s+b^2)\left\{\left(\frac{|y|^2}{4}Q,\rho\right)\right\}=  (\tgamma_s-|\b|^2)\left\{(Q,\rho)\right\}+\mathcal R_7(\P)+\mathcal S_5(\mathcal P)Mod(t)\\
  & + & O\left((\P^2+|Mod(t)|)\|\e\|_{L^2}+\mathcal P^7+\|\e\|_{L^2}^2+\|\e\|_{H^1}^3\right)
\eea
with $\mathcal R_7 $ polynomial satisfying \fref{esttnsplu}, $\mathcal S_5$ polynomial in $\mathcal P$ of order 1, and where we have used the orthogonality condition \eqref{ortho5} and the computation of the inner product \eqref{pscal5}. We now inject \fref{lawb} and the definition \fref{defdzeroalpha} and get:
\bee
& & (\tgamma_s-|\b|^2+d_1(\alpha,\alpha))\left\{(Q,\rho)\right\}=\mathcal R_8(\P)+\mathcal S_6(\mathcal P)Mod(t)\\
  & + & O\left((\P^2+|Mod(t)|)\|\e\|_{L^2}+\mathcal P^7+\|\e\|_{L^2}^2+\|\e\|_{H^1}^3+\left|\|u\|_{L^2}-\|Q\|_{L^2}\right|\right)
\eee
with $\mathcal R_8 $ polynomial satisfying \fref{esttnsplu}, $\mathcal S_6$ polynomial in $\mathcal P$ of order 1.\\

{\bf step 8} Conclusion.\\

Putting together \fref{lawb}, \fref{lawlambda}, \fref{lawalpha}, \fref{lawbeta}, \fref{lawtgamma} and estimating the nonlinear interaction terms using Sobolev embeddings yields: \bee
& & \left\{A+\mathcal S_7(\mathcal P)\right\}Mod(t)  =  \mathcal R_9(\P)\\
& + &  O\left([\P^2+|Mod(t)|]\|\e\|_{L^2}+\|\e\|_{L^2}^2+\|\e\|_{H^1}^3+\left|\|u\|_{L^2}-\|Q\|_{L^2}\right|\right)
\eee for some polynomial vector $\mathcal R_9$ satisfying \fref{esttnsplu}, some $O(1)$ invertible matrix $A$ and some polynomial matrix $\mathcal S_7(\mathcal P)$ of order 1 in $\mathcal P$. Inverting this relation to compute $Mod(t)$ and computing the Taylor expansion of $(A+\mathcal S_7)^{-1}$ to sufficiently high order yields \fref{computationparamatersappenidxgeneralized} which implies \eqref{law2}.\\
We now now inject \eqref{law2}, \fref{chiofyeiyo} and the almost critical mass assumption \fref{almostcriticalmass} into \fref{lawlambda} and \fref{law2:bis} follows.


\section*{Appendix B: A priori bounds on an ODE}


In this appendix, we study an ordinary differential system which will intervene in Appendix C and Appendix D.

\begin{lemma}
\label{ps50:0}
Let $Z=(Z_1,Z_2):\RR\goto\RR^2$ satisfying the following ordinary differential system:
\be\label{ps50}
Z_s=\left(\begin{array}{cc}
\ds 0 & -2\\
\ds \frac{\varsigma}{s^2} & 0
\end{array}\right)Z+F
\ee
where $\varsigma>0$ is a constant and where $F=(F_1,F_2):\RR\goto\RR^2$. Assume also that
\be\label{ps51}
\lim_{s\goto +\infty}Z(s)=0\textrm{ and }|F(s)|\lesssim \frac{1}{s^3}.
\ee
Then, we have the following estimate for $Z$ for $s\geq 2$:
\be
\label{ps52}
|Z_1(s)|+s|Z_2(s)|\lesssim \int_s^{+\infty}\left(|F_1(\sigma)|+|F_2(\sigma)|\sigma\right)\log \sigma d\sigma.
\ee
\end{lemma}

\noindent{\bf Proof of Lemma \ref{ps50:0}} We start with the homogeneous differential equation. We look for $Z_h$ solution of:
\be
\label{ps53}
(Z_h)_s=\left(\begin{array}{cc}
\ds 0 & -2\\
\ds \frac{\varsigma}{s^2} & 0
\end{array}\right)Z_h.
\ee
We have:
$$Z_h=\left(\begin{array}{c}
\ds z\\
\ds -\frac{1}{2}z_s
\end{array}\right)
$$
where $z:\RR\goto\RR^2$ satisfies the following second order differential equation:
$$z_{ss}+\frac{2\varsigma}{s^2}z=0.
$$
We look for solutions of the form $z=s^\vartheta$, which yields the equation:
$$\vartheta^2-\vartheta+2\varsigma=0.
$$
We find a basis $z_+, z_-$  with wronskian$$W(z_+,z_-)=-\frac{1}{2}z_+(z_-)_s+\frac{1}{2}(z_+)_sz_-$$ as follows:
$$\textrm{if }0< \varsigma<\frac{1}{8},\textrm{ then }\vartheta_\pm=\frac{1\pm\sqrt{1-8\varsigma}}{2},\,z_\pm=s^{\vartheta_\pm}\textrm{ and }W=\sqrt{1-8\varsigma},
$$
$$
\textrm{if }\varsigma=\frac{1}{8},\textrm{ then },\,z_+=\sqrt{s}\log s,\,z_-=\sqrt{s}\textrm{ and }W=\frac{1}{2},
$$
and
$$
\begin{array}{l}
\ds\textrm{if }\varsigma>\frac{1}{8}, \ds\textrm{ then },\,z_+=\sqrt{s}\cos\left(\frac{\sqrt{8\varsigma-1}}{2}\log(s)\right),\\
\ds z_-=\sqrt{s}\sin\left(\frac{\sqrt{8\varsigma-1}}{2}\log(s)\right)\textrm{ and }W=-\frac{\sqrt{8\varsigma-1}}{2}.
\end{array}
$$
We define:
$$
Z_\pm=\left(\begin{array}{c}
\ds z_\pm\\
\ds -\frac{1}{2}(z_\pm)_s
\end{array}\right).
$$
The assumption $\varsigma>0$ yields the bound 
\be
\label{estzh}
|Z_{\pm}(s)|\lesssim s \ \ \mbox{for} \ \ s\geq 1
\ee
and for all $2\leq s\leq \sigma$: 
\be
\label{boundcrossed}
\begin{array}{l}
\displaystyle |(Z_{\pm})_1(s)||(Z_{\mp})_1(\sigma)|\lesssim \sigma\log \sigma,\,|(Z_{\pm})_1(s)||(Z_{\mp})_2(\sigma)|\lesssim \log \sigma,\\[1mm]
\displaystyle |(Z_{\pm})_2(s)||(Z_{\mp})_1(\sigma)|\lesssim \frac{\sigma}{s}\log \sigma,\,|(Z_{\pm})_2(s)||(Z_{\mp})_2(\sigma)|\lesssim \frac{\log \sigma}{s}.
\end{array}
\ee
Moreover, if $Z_h$ solves the homogeneous equation \fref{ps53} with 
\be\label{ps62}
\lim_{s\goto +\infty}Z_h=0,\textrm{ then }Z_h\equiv 0.
\ee
We now turn to the inhomogeneous ordinary differential equation \eqref{ps50} which we solve using the variation of constants method:
\be\label{ps63}
Z(s)=a_+(s)Z_+(s)+a_-(s)Z_-(s)
\ee
where
\be\label{ps65}
\begin{array}{l}
\ds a_+=a_+^0-\int_s^{+\infty}\frac{F_1(Z_-)_2-F_2(Z_-)_1}{W}d\sigma,\\
\ds a_-=a_-^0-\int_s^{+\infty}\frac{F_2(Z_+)_1-F_1(Z_+)_2}{W}d\sigma,
\end{array}
\ee
where $a_\pm^0$ are real constants and where $W$ is the wronskian of $z_\pm$. Observe from \eqref{ps51}, \fref{estzh}, \fref{boundcrossed} that the integrals in \eqref{ps65} are absolutely convergent, and there holds the bounds:
\bea
\label{boundZ1}
\nonumber
& & \left|(Z_+)_1(s)\int_s^{+\infty}\frac{F_1(Z_-)_2-F_2(Z_-)_1}{W}d\sigma+(Z_-)_1(s)\int_s^{+\infty}\frac{F_2(Z_+)_1-F_1(Z_+)_2}{W}d\sigma\right|\\
& \lesssim & \int_s^{+\infty} \left(|F_1(\sigma)|+|F_2(\sigma)|\sigma\right)\log \sigma d\sigma\to 0 \ \ \mbox{as} \ \ s\to +\infty,
\eea
and similarly:
\bea
\label{boundZ2}
\nonumber
& & \left|(Z_+)_2(s)\int_s^{+\infty}\frac{F_1(Z_-)_2-F_2(Z_-)_1}{W}d\sigma+(Z_-)_2(s)\int_s^{+\infty}\frac{F_2(Z_+)_1-F_1(Z_+)_2}{W}d\sigma\right|\\
& \lesssim & \frac{1}{s}\int_s^{+\infty} \left(|F_1(\sigma)|+|F_2(\sigma)|\sigma\right)\log \sigma d\sigma\to 0 \ \ \mbox{as} \ \ s\to +\infty.
\eea
We conclude from \fref{ps62} that $a_+^0=a_-^0=0$ and \fref{ps63}, \fref{ps65}, \fref{boundZ1}, \fref{boundZ2} now yield \fref{ps52}. This concludes the proof of Lemma \ref{ps50:0}. 


\section*{Appendix C: Proof of the estimate \eqref{diff5}}


This appendix is devoted to the proof of the estimate \eqref{diff5}. We in fact claim the stronger estimates:
\be
\label{eststronger}
|\lambda-\lambda_c|+|b-b_c|\lesssim |t|^6, \ \ |\alpha-\alpha_c|+|\beta-\beta_c|\lesssim |t|^5, \ \ |\gamma-\gamma_c|\lesssim t^4.
\ee

{\bf step 1} Improved bound on the modulation equations.\\

Let $Mod(t)$ be the vector of modulation equations given by \fref{defmod},  and rewrite the modulation equations \fref{computationparamatersappenidxgeneralized} using the bounds \fref{backwarddd}:
\bea
\label{computationparamatersappenidxgeneralizedbis}
\nonumber Mod(t) & =  & \mathcal R(\mathcal P)+  O\left(\P^7+\P^2\|\e\|_{L^2}\right)\\
& = & \mathcal R(\mathcal P)+O(|t|^6)
\eea
with $R(\mathcal P)$ polynomial in $\mathcal P$ satisfying \fref{esttnsplu}. Unfortunately, this estimate is just not enough to conclude and we need to slightly improve the inner products \eqref{pscal1} and \eqref{pscal3} to gain a cancellation on the null space for the $(b,\b)$ laws from \fref{computationparamatersappenidxgeneralizedbis} to 
\bea
\label{boundbetter}
\nonumber \left(b_s+b^2-d_0(\alpha,\alpha),\beta_s+b\beta-c_0(\alpha)\lambda-\b_3\l^3\right) & =  &  \tilde{\mathcal R}(\mathcal P)+  O\left(\P^7+\P^3\|\e\|_{L^2}\right)\\
& = & \tilde{\mathcal R}(\mathcal P)+  O(|t|^7)
\eea
with $\tilde{\mathcal R}$ satisfying \fref{esttnsplu}.\\
Let us assume the improved bound \fref{boundbetter} and conclude the proof of \eqref{eststronger}.\\
We rewrite \fref{computationparamatersappenidxgeneralizedbis} for $Mod(t),Mod_c(t)$ using the fact that $\mathcal R$ vanishes at least at order 2 at the origin from \fref{esttnsplu}:
\be
\label{estdifferencekey}
\left|Mod(t)-Mod_c(t)\right|\lesssim |t|^2|\mathcal P-\mathcal P_c|+O(|t|^6).
\ee
Using similarly the improved bound \fref{boundbetter}, the definition of $d_0$ \fref{defdzeroalpha} and the degeneracy $\alpha,\alpha_c=O(t^2)$, we get:
\bea
\label{imprbehffP}
\nonumber & & \left|(b_s+b^2)-((b_c)_s+b_c^2)\right|\\
\nonumber & + & \left|(\beta_s+b\beta-c_0(\alpha)\lambda-\b_3\l^3)-((\beta_c)_s+b_c\beta_c-c_0(\alpha_c)\lambda_c-\b_3\l_c^3)\right|\\
& \lesssim &  |t|^2|\mathcal P-\mathcal P_c|+O(|t|^7).
\eea

{\bf step 2} Estimates for $\l-\l_c$ and $b-b_c$.\\
 
Let us define: 
\be
\label{cnkheoheo}
\underline{\P}=\P-\P_c=O(|t|^2)
\ee from \eqref{bis:ex17}, \fref{bis:ex19:7}, \fref{bis:ex19:9}. We claim the bound:
\be
\label{diff36}
|b-b_c|+|\l-\l_c|\lesssim |t|^6+\int_0^t |\underline{\P}(\tau)|d\tau+|t|\left(\int_0^t|\underline{\P}(\tau)|\frac{d\tau}{\tau}\right).
\ee
Indeed, we first have:
\be
\label{diff27}
\l_t+\frac{b}{\l}=\frac{1}{\l}\left(\frac{\l_s}{\l}+b\right), \ \ b_t+\frac{b^2}{\l^2}=\frac{1}{\l^2}\left(b_s+b^2\right).
\ee
Observe now from \fref{bis:ex19:7}, \fref{bis:ex19:9} that $$\lambda(t)=-\frac{t}{C_0}+O(|t|^3), \ \ b(t)=-\frac{t}{C^2_0}+O(|t|^3)$$ from which:
$$\frac{b}{\l}-\frac{b_c}{\l_c}=-\frac{C_0}{t}\underline{b}+\frac{1}{t}\underline{\l}+O(t\underline{\mathcal P}).$$ Injecting this into \fref{diff27} together with \fref{estdifferencekey} for $\lambda$ and the improved bound \fref{imprbehffP} for $b$ yields:
\be\label{diff28}
\underline{\l}_t-\frac{C_0}{t}\underline{b}+\frac{1}{t}\underline{\l}=F_1,\ \  \underline{b}_t-\frac{2}{t}\underline{b}+\frac{2}{C_0t}\underline{\l}=F_2,
\ee
where $F_1$ and $F_2$ satisfy:
\be\label{diff29}
|F_1|+|F_2|\lesssim |t|^5+|\underline{\P}|.
\ee
We rewrite \eqref{diff28}
\be\label{diff30}
Z_t=\frac{1}{t}MZ+F,
\ee
where $Z$, $M$ and $F$ are given by:
\be\label{diff31}
Z=\left(\begin{array}{c}
\underline{\l}\\
\underline{b}
\end{array}\right),\,M=\left(\begin{array}{cc}
-1 & C_0\\
-\frac{2}{C_0} & 2
\end{array}\right),\,F=\left(\begin{array}{c}
F_1\\
F_2
\end{array}\right).
\ee
The eigenvalues of the matrix $M$ are 0 and 1, and hence the system can be rewritten in an eigenbasis: 
\be
\label{mpjpgjerpo}
\tilde{Z}_t=\frac{D}{t}\tilde{Z}+\tilde{F} \ \ \mbox{with} \ \ D=\left(\begin{array}{cc}
0 & 0\\
0 & 1\\
\end{array}\right),
\ee  with from \fref{bis:ex19:7}, \fref{bis:ex19:9}, \fref{diff29}:$$|\tilde{F}|\lesssim |t|^5+|\underline{\P}|, \ \ \ \frac{\tilde{Z}}{t}\to 0 \ \ \mbox{as} \ \ t\to 0.$$ The explicit integration of \fref{mpjpgjerpo} yields: $$|Z(t)|\lesssim |\tilde{Z}(t)|\lesssim \int_t^0|\tilde{F}(\tau)|d\tau+|t|\int_t^0\left|\frac{\tilde{F}(\tau)}{|\tau|}\right|d\tau \lesssim |t|^6+\int_0^t |\underline{\P}(\tau)|d\tau+|t|\left(\int_0^t|\underline{\P}(\tau)|\frac{d\tau}{\tau}\right)$$ which implies \fref{diff36}.\\

{\bf step 3} Estimates for $\a-\a_c$ and $\b-\b_c$.\\

We claim the bound: 
\be
\label{diff48}
|\alpha-\alpha_c|+|\beta-\beta_c|\lesssim |t^6\ln(|t|)|+|t|\left(\int_0^t\frac{|\log(\tau)|}{|\tau|}|\underline{\P}(\tau)|d\tau\right).
\ee
Indeed, we first have:
\be\label{diff37}
\begin{array}{l}
\ds\a_t-2\frac{\b}{\l}=\frac{1}{\l}\left(\frac{\a_s}{\l}-2\b\right),\\ 
\ds \b_t+\frac{b\b}{\l^2}-\frac{c_0(\a)\l+\b_3\l^3}{\l^2}=\frac{1}{\l^2}\left(\b_s+b\b-c_0(\a)\l-\b_3\l^3\right).
\end{array}
\ee
We now derive from \eqref{bis:ex17}-\eqref{bis:ex19:9} the bounds:
$$\frac{\b}{\l}-\frac{\b_c}{\l_c}=-\frac{C_0}{t}\underline{\b}+O(\underline{\P}),$$
$$\frac{b\b}{\l^2}-\frac{c_0(\a)\l+\b_3\l^3}{\l^2}-\frac{b_c\b_c}{\l_c^2}+\frac{c_0(\a_c)\l_c+\b_3\l_c^3}{\l_c^2}=-\frac{1}{t}\underline{\beta}+\frac{C_0}{t}c_0(\underline{\alpha})+O(|\underline{\mathcal P}|).$$ Injecting this together with \fref{estdifferencekey} and the improved bound \fref{imprbehffP} into \fref{diff37} yields:
\be\label{diff38}
\underline{\a}_t+\frac{2C_0}{t}\underline{\b}=G_1,\ \  \underline{\b}_t-\frac{1}{t}\underline{\b}+\frac{C_0}{t}c_0(\underline{\a})=G_2,
\ee
with:
\be\label{diff39}
|G_1|+|G_2|\lesssim t^5+|\underline{\P}|.
\ee
We set:
\be\label{diff40}
\underline{\eta}=\frac{\underline{\b}}{t},
\ee
and rewrite \eqref{diff38} as:
$$\begin{array}{l}
\underline{\a}_t+2C_0\underline{\eta}=G_1, \ \ \underline{\eta}_t+\frac{C_0}{t^2}c_0(\underline{\a})=\frac{G_2}{t}.
\end{array}
$$
In order to solve this system, we first diagonalize the matrix $c_0$ which is negative definite from \eqref{choicec0} and hence diagonalizable with eigenvalues $r_1<0, r_2<0$. We obtain in an eigenbasis: for  $j=1, 2$,
\be\label{diff43}
\left\{\begin{array}{l}
\ds(\underline{\a}_j)_t=-2C_0\underline{\eta}_j+O(t^5+\underline{\P}),\\
\ds(\underline{\eta}_j)_t=-\frac{C_0r_j\underline{\a}_j}{t^2}+O\left(t^4+\frac{\underline{\P}}{t}\right).
\end{array}\right.
\ee
We then perform the change of variables $s=\frac{1}{|t|}$ and rewrite \fref{diff43}: for $j=1,2$,
$$Z_j=\left(\begin{array}{c}
\underline{\eta}_j\\
\frac{r_j}{2C_0}\underline{\alpha}_j
\end{array}\right), \ \ (Z_j)_s=\left(\begin{array}{cc}
\ds 0 & -2\\
\ds \frac{\varsigma}{s^2} & 0
\end{array}\right)Z_j+F_j
$$
with $\varsigma=-r_jC_0^2>0$ and 
$$F_j=((F_j)_1,(F_j)_2), \ \ (F_j)_1=O\left(\frac{1}{s^6}+\frac{\underline{\P}}{s}\right), \ \ (F_j)_2=O
\left(\frac{1}{s^7}+\frac{\underline{\P}}{s^2}\right).$$
We have from \fref{bis:ex19:7}, \fref{bis:ex19:9}, \fref{cnkheoheo} the bounds: $$\lim_{s\to +\infty}Z_j= 0\ \ \mbox{and} \ \ F_j(s)=O\left( \frac{1}{s^3}\right),$$ and we thus conclude from Lemma \ref{ps50:0}: 
\bee
|\underline{\eta}_j|+s|\underline{\alpha}_j|\lesssim \int_s^{+\infty}\left(\frac{1}{\sigma^6}+\frac{|\underline{\mathcal P}|}{\sigma}\right)\log \sigma d\sigma\lesssim |t^5\log |t||+\int_t^0\frac{|\log(\tau)|}{|\tau|}|\underline{\P}(\tau)|d\tau
\eee
which from \fref{diff40} implies \fref{diff48}.\\

{\bf step 4} Bound on $\underline{\mathcal P}$.\\

We conclude from \fref{diff36}, \fref{diff48} that 
$$|\underline{\P}|\lesssim |t^6\log |t||+\int_t^0|\underline{\mathcal P}(\tau)|d\tau+|t|\int_t^0\frac{|\log(\tau)|}{|\tau|}|\underline{\P}(\tau)|d\tau.
$$
Using the initial a priori bound \fref{cnkheoheo} and iterating several times yields: 
$$ |\underline{\P}|\lesssim t^6|\log |t||^4.
 $$
 Injecting this into \fref{diff36}, \fref{diff48} yields in return:
 \be
 \label{eijvdpodpeu}
 |\l-\l_c|+|b-b_c|\lesssim t^6, \ \ |\alpha-\alpha_c|+|\b-\b_c|\lesssim |t|^5.
 \ee
 
 {\bf step 5} Bound on the phase parameter.\\
 
 We are now in position to derive the key estimate on the phase parameter:
 \be\label{diff57}
|\gamma-\gamma_c|\lesssim t^4.
\ee 
 Indeed, we have:
\be\label{diff53}
\gamma_t-\frac{1+|\b|^2+d_1(\alpha,\alpha)}{\l^2}=\frac{1}{\l^2}\left(\tgamma_s-|\b|^2+d_1(\alpha,\alpha)\right),
\ee
and from the improved bound \fref{eijvdpodpeu}\footnote{which is critical for the scaling parameter}:
\bea
\label{cjoeoujeop}
\nonumber \frac{1+|\b|^2+d_1(\alpha,\alpha)}{\l^2}-\frac{1+|\b_c|^2+d_1(\alpha_c,\alpha_c)}{\l_c^2}& =& O\left(\frac{|\lambda-\lambda_c|}{|t|^3}+|\b-\b_c|+|\alpha-\alpha_c|\right)\\
& = & O(|t|^3).
\eea 
\fref{estdifferencekey}, \fref{diff53} and \fref{cjoeoujeop} now yield:
$$\underline{\gamma}_t=O\left(|t|^3+|\underline{\mathcal P}|\right),
$$
which after integration in time yields \eqref{diff57}. This concludes the proof of \fref{eststronger} assuming \fref{boundbetter}.\\

{\bf step 6} Proof of the improved bound \fref{boundbetter}.\\

We now turn to the proof of \fref{boundbetter}. We change the orthogonality conditions \eqref{ortho1} and \eqref{ortho3} which respectively govern the law of $\b$ and $b$ for 
\be\label{diff8}
(\e_2,\n\S+A\l^2)-(\e_1,\n\T)=0
\ee
and
\be\label{diff9}
(\e_2,\L\S+B\l^2)-(\e_1,\L\T)=0
\ee
where $A$ and $B$ are well localized functions to be chosen. We claim that for a suitable choice of $A,B$, the computation of the modulation equations like for \fref{lawb}, \fref{lawbeta} directly leads to the improved bound \fref{boundbetter}\footnote{We could have chosen this improved set of orthogonality conditions from the very beginning. We chose not to do so because this refinement is at first hand not natural and is really needed only here, see Remark \ref{cneoheoheoei}.}.\\
Indeed, let $M_1, M_2$ be given by \fref{defm1}, \fref{defm2}, and let the complex operator for $\e=\e_1+i\e_2$
; $$\tilde{M}(\e)=\tilde{M}_1(\e_1,\e_2)+\tilde{M_2}(\e_1,\e_2)=M_1(\e_2,\e_2)+iM_2(\e_1,\e_2)-ib\Lambda \e+2i\beta\cdot\nabla \e.$$ An integration by parts yields the adjunction formula: 
$$\Re(\tilde{M}(\e),\overline{f})=\Re(\e,\overline{\tilde{M}(f)}).$$
Given an orthogonality condition $$\Im(\e,\overline{f})=0,$$ the linear term in the computation of the modulation equations is from \eqref{eqe1} \fref{eqe2} up to terms $O(t^3\|\e\|_{L^2})$:$$-(\tilde{M}_1(\e_1,\e_2),f_1)+(\e_2,\pa_sf_1)-(\tilde{M}_2(\e_1,\e_2),f_2)-(\e_1,\pa_sf_2)=-\Re(\e,\overline{\tilde{M}(f)-i\pa_sf}).$$ Following closely the proof of \fref{lawb}, \fref{lawbeta}, the improvement \fref{boundbetter} follows from\footnote{using also the $L^2$ conservation law \fref{coefeou}}:
\be
\label{tobeprovedlvs;vs}
Re(\e,\overline{\tilde{M}(\nabla \qp+\l^2A)-i\pa_s(\nabla \qp+\l^2A)})=O(|t|^3\|\e\|_{L^2}),
\ee
\be
\label{tobeprovedbislm;lkss}
Re(\e,\overline{\tilde{M}(\Lambda \qp+\lambda^2B)-i\pa_s(\Lambda \qp+\lambda^2B)})=-2\Re(\e,\overline{\qp})+O(|t|^3\|\e\|_{L^2}).
\ee
{\it Proof of \fref{tobeprovedlvs;vs}}: We rewrite the equation \eqref{eqq} of $\qp$:
\bea
\label{diff10}
\nonumber & &-ib^2\partial_b\qp +\Delta \qp-\qp+\frac{k(\lambda(t)y+\alpha(t))}{k(\alpha(t))}\qp|\qp|^{2}+ib\L\qp-2i\b\cdot\n\qp\\
&= & O(t^3e^{-c|y|}).
\eea
By differentiating \eqref{diff10}, we obtain:
\bee
& &  -ib^2\partial_b\n\qp -M_1(\n\qp)-iM_2(\n\qp)+ib\L\n\qp-2i\b\cdot\n^2\qp\\
& = & -\frac{\l\n k(\l y+\a)}{k(\a)}|\qp|^2\qp-ib\n\qp+O(t^3e^{-|y|})
\eee
which yields:
\bea
\label{neoheoheointerm}
\nonumber \tilde{M}(\nabla \qp+\l^2A)-i\pa_s(\nabla \qp+\l^2A)& = & \lambda^2\left(\nabla^2k(0)(y,\cdot)Q^3+L_+(A)\right)+ib\nabla \qp\\
& + & O(|t|^3e^{-|y|}).
\eea
From \fref{ellitpicestinamtes}, we then chose $A$ solution to 
$$L_+(A)=-\nabla^2k(0)(y,.)Q^3+a yQ \ \ \mbox{with} \ \ a=\frac{(\nabla^2k(0)(y,.)Q^3,\n Q)}{(\n Q,yQ)},
$$
and thus from \fref{neoheoheointerm} with $f=\n\qp+\l^2A$: 
\bee
Re(\e,\overline{\tilde{M}(f)-i\pa_sf})& = & b\Im(\e,\nabla \overline{\qp})+a\lambda^2\Re(\e,\overline{yQ})+O(|t|^3\|\e\|_{L^2})\\
& = & O(|t|^3\|\e\|_{L^2})
\eee where we used \fref{ortho2}, \fref{diff8} in the last step. This concludes the proof of \fref{tobeprovedlvs;vs}.\\
{\it Proof of \fref{tobeprovedbislm;lkss}}: We similarly derive from \fref{diff10} the relation:
\bee
& & -ib^2\partial_b\L\qp -M_1(\L\qp)-iM_2(\L\qp)+ib\L^2\qp-2i\b\cdot\n\L\qp\\
& = & 2ib^2\partial_b\qp+2\qp-\frac{\l y\cdot\n k(\l y+\a)}{k(\a)}|\qp|^2\qp-2ib\L\qp+2i\b\cdot\n\qp+O(t^3e^{-|y|})
\eee
which implies:
\bea
\label{neoheoheointermbis}
\nonumber & & \tilde{M}(\Lambda\qp+\l^2B)-i\pa_s(\Lambda \qp+\l^2B) =  \lambda^2\left(\nabla^2k(0)(y,y)Q^3+L_+(B)\right)\\
& - & 2ib^2\pa_b\qp-2\qp+2ib\Lambda \qp-2i\beta\cdot\nabla \qp+O(|t|^3e^{-|y|}).
\eea
We hence choose from \fref{ellitpicestinamtes} $B$ solution to: $$L_+(B)=-\nabla^2k(0)(y,y)Q^3.$$ We now observe that $$\pa_b\qp=-i\frac{|y|^2}{4}Q+O(|t|e^{-|y|})$$ so that \fref{neoheoheointermbis} yields:
\bee
& & \left|Re(\e,\overline{\tilde{M}(\Lambda\qp+\l^2B)-i\pa_s(\Lambda \qp+\l^2B)})+2\Re(\e,\overline{\qp})\right|\\
& \lesssim & b^2|\Re(\e,\overline{|y|^2Q})|+b|\Im(\e,\overline{\Lambda \qp})|+\beta|\Im(\e,\overline{\nabla \qp})|+O(|t|^3\|\e\|_{L^2})\\
& = & O(|t|^3\|\e\|_{L^2})
\eee
where we used the orthogonality conditions \fref{diff8}, \fref{ortho4}, \fref{diff9}. This concludes the proof of \fref{tobeprovedbislm;lkss}.\\
This concludes the proof of \fref{boundbetter}.

\begin{remark}
\label{cneoheoheoei}
The estimate \fref{eststronger} is stronger than \eqref{diff5}. We need these estimates because we use $|\l-\l_c|\lesssim t^6$ in the proof of \fref{cjoeoujeop}, and thus it is only to control the phase paramater in \fref{eststronger} that we need the improved bound \fref{boundbetter}, or equivalently the refinement of the orthogonality conditions \eqref{ortho1}, \eqref{ortho3}.
\end{remark}


\section*{Appendix D: Algebraic derivation of \fref{estunfin}, \fref{estdeuxfin}, \fref{esttroisfin}, \fref{estquatrefin} }


Let us recall the definitions \fref{defnmsnkfo}, \fref{pscl11}, \fref{pscl12} of $\mq , \mq_1, \mq_2$.\\

{\bf Proof of \fref{estunfin}, \fref{tobecomputed}}: We compute from \fref{defnmsnkfo}, \fref{pscl11}, \fref{pscl12} and the control of the modulation parameters on $u_c$:
\bee
\mq(\nabla Q+\l_c^2A_2)-i\pa_s(\nabla Q+\l_c^2A_2) & = & L_+(\nabla Q+\l_c^2A_2)-6\l_c^2 QT_2^0\nabla Q\\
& - & \frac{3}{2}\l_c^2\nabla^2k(0)(y,y)Q^2\nabla Q + O(\l_c^3e^{-c|y|})
\eee 
where from \fref{eqT2S2} and the degeneracy $|\alpha_c|\lesssim \l_c^2$ it suffices to let $T_2^0$ be the radial solution to:
\be
\label{esttdexunot}
 L_+(T^0_2)=\frac{1}{2}\nabla^2k(0)(y,y)Q^3.
 \ee
To fulfill \fref{estunfin}, we recall $L_+(\nabla Q)=0$ and chose $A_2$ solution to\be\label{ps38}
L_+(A_2)=6QT_2^0\n Q+\frac{3}{2}\nabla^2k(0)(y,y)Q^2\n Q -a_1 yQ
\ee
with from \fref{ellitpicestinamtes} and the fact that $Q, T_2^0$ are radial:
\be\label{ps40}
a_1= -\frac{(6QT_2^0\pa_j Q+\frac{3}{2}\nabla^2k(0)(y,y)Q^2\pa_j Q ,\pa_j Q)}{\int Q^2}.
\ee
We now observe that $a_1$ can be computed in an explicit form. Indeed, we first differentiate twice the equation of $Q$ \eqref{eqQ} to derive:
\be\label{ps41}
L_+(\Delta Q)=6|\n Q|^2Q,
\ee
and hence using the fact that $Q, T_2^0$ are radial and \fref{esttdexunot}: for $j=1,2$,
\bee
2(6QT_2^0\pa_j Q,\pa_j Q) & = & 6(T_2^0,Q|\n Q|^2)=(T_2^0,L_+(\Delta Q))=(L_+(T_2^0),\Delta Q)\\
&  = & \left(\frac{1}{2}\nabla^2k(0)(y,y)Q^3,\Delta Q\right)\\
& = & -\left(\frac{3}{2}\nabla^2k(0)(y,y)Q^2\n Q,\n Q\right)-\left(\frac{1}{2}\n(\nabla^2k(0)(y,y))Q^3,\n Q\right)\\
&  = & -\left(\frac{3}{2}\nabla^2k(0)(y,y)Q^2\n Q,\n Q\right)+\frac{1}{8}\int \Delta(\nabla^2k(0)(y,y))Q^4\\
&  = & -\left(\frac{3}{2}\nabla^2k(0)(y,y)Q^2\n Q,\n Q\right)+\frac{\textrm{tr}(\nabla^2k(0))}{4}\int Q^4\\
& = & -3\left(\nabla^2k(0)(y,y)Q^2\pa_j Q,\pa_j Q\right)+\frac{\textrm{tr}(\nabla^2k(0))}{4}\int Q^4
\eee
which together with \eqref{ps40} yields:
\be\label{ps43}
a_1= -\frac{\textrm{tr}(\nabla^2k(0))}{8}\frac{\int Q^4}{\int Q^2}>0
\ee
since $\nabla^2k(0)$ is negative definite, and \fref{tobecomputed} is proved.\\

{\bf Proof of \fref{estdeuxfin}}: We compute from \fref{defnmsnkfo}, \fref{pscl11}, \fref{pscl12} and the control of the modulation parameters on $u_c$:
\bee
M^{(4)}(i(yQ+\l_c^2B_2))-i\pa_s\left[i(y Q+\l_c^2B_2)\right] & = &i L_-(y Q+\l_c^2B_2)-2i\l_c^2 QT_2^0y Q\\
& - & \frac{i}{2}\l_c^2\nabla^2k(0)(y,y)Q^2yQ + O(\l_c^3e^{-c|y|}).
\eee 
From $L-(yQ)=-2\nabla Q$, \fref{estdeuxfin} follows by choosing $B_2$ solution to:
$$L_-(B_2)=-2A_2+2QT_2^0yQ+\frac{\nabla^2k(0)(y,y)}{2}Q^2yQ.
$$
which is solvable from \fref{ellitpicestinamtesbis} since its RHS is orthogonal to $Q$ by radial symmetry and \fref{ps38}.\\

{\bf Proof of \fref{esttroisfin}}: 

We compute from \fref{defnmsnkfo}, \fref{pscl11}, \fref{pscl12},  \eqref{solapp} and the control of the modulation parameters on $u_c$:
\bea
\label{xixi1}
\nonumber & & \mq(\Lambda P_{\mathcal P_c}+E_2+E_3+iD_3)-i\pa_s\left[\Lambda P_{\mathcal P_c}+E_2+E_3+iD_3\right]\\
\nonumber &+ &(c_0(\a_c)\l_c+\b_3\l_c^3)\cdot y\left[\Lambda P_{\mathcal P_c}+E_2+E_3+iD_3\right]\\
\nonumber & = &  \mq(\Lambda P_{\mathcal P_c})-i\p_s\Lambda T_2+(c_0(\a_c)\l_c+\b_3\l_c^3)\cdot y\Lambda Q\\
& & +L_+(E_2)+L_+(E_3)+iL_-(D_3)-i\p_s(E_2)+O(\l_c^4e^{-|y|}).
\eea

We now evaluate $\mq(\Lambda P_{\mathcal P_c})$. From \eqref{eqpbis}, \eqref{bis:ex17}-\eqref{bis:ex19:9} and the fact that $\p_\l\pp=O(\P)$, $\p_\a\pp=O(\P)$, $\p_b\pp=O(\P^2)$ and $\p_\b\pp=O(\P^2)$, $\pp$ satisfies the following equation:
\bea
\label{pscl4} 
 \nonumber & & -i\l b\p_\l\pp+\Delta \pp-\pp+\frac{k(\lambda(t)y+\alpha(t))}{k(\alpha(t))}\pp|\pp|^{2}-(\l c_0(\a)+\b_3\l^3)\cdot y\pp\\
 & = & O(\l^4e^{-c|y|}).
\eea
Differentiating this relation yields:
\bea
\label{pscl5}
& & ib\l\p_\l(\L\pp)+M(\L\pp)+(\l c_0(\a)+\b_3\l^3)\cdot y\Lambda\pp=-2\pp \\
\nonumber & + &  \l y\cdot\frac{\n k(\l y+\a)}{k(\a)}|\pp|^2\pp -3(\l c_0(\a)+\b_3\l^3)\cdot y\pp -2ib\l\p_\l(\pp)
+O(\l^4e^{-|y|}).
\eea
Finally, \eqref{xixi1}, \eqref{pscl5} and the control of the modulation parameters on $u_c$ yield:
\bea
\label{xixi2}
\nonumber & & \mq(\Lambda P_{\mathcal P_c}+E_2+E_3+iD_3)-i\pa_s\left[\Lambda P_{\mathcal P_c}+E_2+E_3+iD_3\right]\\
\nonumber &+ &(c_0(\a_c)\l_c+\b_3\l_c^3)\cdot y\left[\Lambda P_{\mathcal P_c}+E_2+E_3+iD_3\right]\\
\nonumber & = &  -i\p_s\Lambda T_2-ib\l\p_\l\Lambda T_2 -2P_{\mathcal P_c}+\l y\cdot\frac{\n k(\l y+\a)}{k(\a)}Q^3\\
\nonumber & & -3(c_0(\a_c)\l_c+\b_3\l_c^3)\cdot y Q -2ib_c\l_c\p_\l(T_2)\\
\nonumber & & +L_+(E_2)+L_+(E_3)+iL_-(D_3)-i\p_s(E_2)+O(\l_c^4e^{-|y|})\\
\nonumber & = &  -2P_{\mathcal P_c}+\nabla^2k(0)(y,y)\l_c^2Q^3+\nabla^2k(0)(y,\a_c)\l_c Q^3+\frac{1}{2}\nabla^3k(0)(y,y,y)\l_c^3Q^3\\
\nonumber & & -3(c_0(\a_c)\l_c+\b_3\l_c^3)\cdot y Q -2ib_c\l_c\p_\l(T_2)\\
& & +L_+(E_2)+L_+(E_3)+iL_-(D_3)-i\p_s(E_2)+O(\l_c^4e^{-|y|}).
\eea
To fulfill  \fref{esttroisfin}, we choose $E_2$, $E_3$ and $D_3$ solutions to:
\be\label{pscl24}
L_+(E_2)=-\nabla^2k(0)(y,y)\l_c^2Q^3-\nabla^2k(0)(y,\a_c)\l_cQ^3 +3c_0(\a_c)\l_c\cdot y Q 
+\l_c a_2(\a_c)\cdot yQ,
\ee
where $a_2$ is given by:
\be\label{pscl25}
(a_2(.))_j=\frac{2(\nabla^2k(0)(y,.)Q^3,\partial_j Q)}{(yQ,\partial_j Q)}-3(c_0(.))_j, j=1,2.
\ee
\be\label{pscl26}
L_+(E_3)=-\frac{1}{2}\nabla^3k(0)(y,y,y)\l_c^3Q^3+3\b_3\l_c^3\cdot y Q +\l_c^3a_3\cdot yQ,
\ee
where $a_3$ is given by:
\be\label{pscl27}
(a_3)_j=\frac{(\nabla^3k(0)(y,y,y)Q^3,\partial_j Q)}{2(yQ,\partial_j Q)}-3(\b_3)_j,
\ee
and
\be\label{pscl28}
L_-(D_3)=\frac{2}{C_0}\l^2_c\p_\l T_2-\frac{1}{C_0}\l^2_c\p_\l(E_2).
\ee
From \fref{ellitpicestinamtes} and the definition \eqref{pscl25} and \eqref{pscl27} of $a_2$ and $a_3$, we may solve \eqref{pscl24} and \eqref{pscl26}. To conclude the proof of \eqref{esttroisfin}, it remains to prove that we may solve \eqref{pscl28}. In fact, we have:
\be\label{pscl29}
(E_2,Q)=-\frac{1}{2}(L_+(E_2),\L Q)=\frac{1}{2}(\nabla^2k(0)(y,y)Q^3,\L Q)\l_c^2=0,
\ee
where we have used \eqref{pscl24} and \eqref{explicitcomp}. Moreover, we have from \eqref{calcultdeuxq}:
\be\label{pscl30}
(T_2,Q)=0.
\ee
We deduce from \eqref{ellitpicestinamtesbis}, \eqref{pscl29} and \eqref{pscl30} that we may indeed solve \eqref{pscl28}.\\

{\bf Proof of \fref{estquatrefin}}: We compute from \fref{defnmsnkfo}, \fref{pscl11}, \fref{pscl12} and the control of the modulation parameters on $u_c$:
\bea
\label{xixi3}
\nonumber & & \mq(i(|y|^2 P_{\mathcal P_c}+F_2))-i\pa_s\left[i(|y|^2 P_{\mathcal P_c}+F_2)\right]\\
& = &  \mq(i|y|^2 P_{\mathcal P_c})+ iL_-(F_2)+O(\l_c^3e^{-|y|}).
\eea

We now evaluate $\mq(i|y|^2 P_{\mathcal P_c})$. We deduce from \eqref{pscl4}:
\begin{equation}\label{pscl6}
\ds \mq(i|y|^2 P_{\mathcal P_c})=-4i\L P_{\mathcal P_c}+O(\l_c^3e^{-|y|}).
\end{equation} 
\eqref{xixi3} and \eqref{pscl6} yield:
\bee\nonumber & & \mq(i(|y|^2 P_{\mathcal P_c}+F_2))-i\pa_s\left[i(|y|^2 P_{\mathcal P_c}+F_2)\right]\\
& = &  -4i\L P_{\mathcal P_c}+ iL_-(F_2)+O(\l_c^3e^{-|y|}).
\eee
To fulfill  \fref{estquatrefin}, we chose $F_2$ solution to:
$$L_-(F_2)=-4E_2,
$$
which is possible from \fref{ellitpicestinamtesbis} and \eqref{pscl29}.\\

{\bf Proof of \fref{estcinqfin}}: \fref{estcinqfin} follows from the fact that $\mq(\rho)=L_+(\rho)+O(\l_c^2e^{-|y|})$ and $L_+(\rho)=|y|^2Q$.


\begin{thebibliography}{10}


\bibitem{Ba} Banica, V., Remarks on the blow-up for the Schr\"odinger equation with critical mass on a plane domain. Ann. Sc. Norm. Super. Pisa Cl. Sci. (5) 3 (2004), no. 1, 139--170.

%
\bibitem{BCD} Banica, V.; Carles, R.; Duyckaerts, T., Minimal blow-up solutions to the mass-critical inhomogeneous focusing NLS equation. Preprint.

\bibitem{Mas} Blanchet, A.; Carrillo, J. A.; Masmoudi, N., Infinite time aggregation for the critical Patlak-Keller-Segel model in $\Bbb R^2$. Comm. Pure Appl. Math. 61 (2008), no. 10, 1449--1481

\bibitem{BW} Bourgain, J.; Wang, W., Construction of blowup solutions for the nonlinear Schr\"odinger equation with critical nonlinearity, Ann. Scuola Norm. Sup. Pisa Cl. Sci. (4) \textbf{25} (1997), no. 1-2, 197--215.

\bibitem{BGR}  Burq, N.; G\'erard, P.; Rapha\"el, P., Log-log blow up solutions for the mass critical NLS on a manifold, in preparation.

\bibitem{BGT} Burq, N.; G\'erard, P.; Tzvetkov, N. Two singular dynamics of the nonlinear Schr\"odinger equation on a plane domain. Geom. Funct. Anal. 13 (2003), no. 1, 1--19.

\bibitem{nakspectral} Chang, S.M.; Gustafson, S.; Nakanishi, K.; Tsai, T.P., Spectra of linearized operators for NLS solitary waves. SIAM J. Math. Anal. 39 (2007/08), no. 4, 1070--1111.

\bibitem{DM1} Duyckaerts, T.; Merle, F., Dynamic of threshold solutions for energy-critical NLS, Geom. Funct. Anal. 18 (2009), no. 6, 1787--1840.

\bibitem{FMR} Fibich, G.; Merle, F.; Rapha\"el, P., Proof of a spectral property related to the singularity formation for the $L\sp 2$ critical nonlinear Schr\"odinger equation. Phys. D 220 (2006), no. 1, 1--13. 
%
\bibitem{GNN} Gidas, B.; Ni, W.M.; Nirenberg, L.,
Symmetry and related properties via the maximum principle,
Comm. Math. Phys. {\bf 68} (1979), 209---243.

%
\bibitem{GV} Ginibre, J.; Velo, G., On a class of nonlinear Schr\"odinger equations. I. The Cauchy problem, general case. J.
Funct. Anal. 32 (1979), no. 1, 1--32. 

\bibitem{KM} Kenig, C. E.; Merle, F., Global well-posedness, scattering and blow-up for the energy-critical, focusing, non-linear Schr\"odinger equation in the radial case, Invent. Math. 166 (2006), no. 3, 645--675.

\bibitem{visanetal} Killip, R.; Li, D.; Visan, M.; Zhang, X., Characterization of minimal-mass blowup solutions to the focusing mass-critical NLS, SIAM J. Math. Anal. 41 (2009), no. 1, 219--236.

\bibitem{KLR} Krieger, J.; Lenzmann, E.; Rapha\"el, P, On the stability of pseudo conformal for $L^2$ critical Hartree-NLS, to appear in Ann. Henri Poincare.

\bibitem{KMR} Krieger, J.; Martel, Y.; Rapha\"el, P., Two soliton solutions to the gravitational Hartree equation, to appear in Comm. Pure App. Math.

\bibitem{KS} Krieger, J.; Schlag, W. Non-generic blow-up solutions for the critical focusing NLS in 1-D. J. Eur. Math. Soc. (JEMS) 11 (2009), no. 1, 1--125.

%
\bibitem{Kw} Kwong, M. K., Uniqueness of positive solutions of $\Delta u-u+u\sp p=0$ in ${R}\sp n$. Arch. Rational Mech. Anal. 105 (1989), no. 3, 243--266.

\bibitem{Lions} Lions, P.-L.; The concentration-compactness principle in the calculus of variations. The locally compact case. I. Ann. Inst. H. Poincar\'e Anal. Non Linéaire 1 (1984), no. 2, 109--145. 


\bibitem{Marteluniqueness} Martel, Y., Asymptotic $N$-soliton-like solutions of the subcritical and critical generalized Korteweg-de Vries equations. Amer. J. Math. 127 (2005), no. 5, 1103--1140.

\bibitem{MM1} Martel, Y.; Merle, F., Stability of blow-up profile and lower bounds for blow-up rate for the critical generalized KdV equation. Ann. of Math. (2) 155 (2002), no. 1, 235--280.

\bibitem{MMduke} Martel, Y.; Merle, F., Nonexistence of blow-up solution with minimal $L\sp 2$-mass for the critical gKdV equation. Duke Math. J. 115 (2002), no. 2, 385--408.

\bibitem{MMmulti} Martel, Y.; Merle, F., Multi solitary waves for nonlinear Schr\"odinger equations. Ann. Inst. H. Poincar\'e Anal. Non Lin\'eaire 23 (2006), no. 6, 849--864.

\bibitem{M1} Merle, F., Determination of minimal blow-up solutions with minimal mass for nonlinear Schr\"odinger equations with critical power. Duke Math. J. \textbf{69} (1993), no. 2, 427--454.

\bibitem{Merlemulti} Merle, F., Construction of solutions with exactly $k$ blow-up points for the Schr\"odinger equation with critical nonlinearity. Comm. Math. Phys. 129 (1990), no. 2, 223--240. 

\bibitem{M2} Merle, F., Nonexistence of minimal blow-up solutions of equations $iu_t=-\Delta u-k(x)|u|^{4/N}u$ in $\RR^N$. Ann. Inst. H. Poincar\'e Phys. Th\'eor. \textbf{64} (1996), no. 1, 33--85.

\bibitem{MR1} Merle, F.; Rapha\"el, P., The blow-up dynamic and upper bound on the blow-up rate for critical nonlinear Schr\"odinger equation. Ann. of Math. (2) 161 (2005), no. 1, 157--222.

\bibitem{MR2} Merle, F.; Rapha\"el, P., Sharp upper bound on the blow-up rate for the critical nonlinear Schr\"odinger equation. Geom. Funct. Anal. 13 (2003), no. 3, 591--642.

\bibitem{MR3} Merle, F.; Rapha\"el, P., On universality of blow-up profile for $L\sp 2$ critical nonlinear Schr\"odinger equation. Invent. Math. 156 (2004), no. 3, 565--672

\bibitem{MR4} Merle, F.; Rapha\"el, P., On a sharp lower bound on the blow-up rate for the $L^2$ critical nonlinear Schr\"odinger equation, J. Amer. Math. Soc. 19 (2006), no. 1, 37--90.

\bibitem{PR} Planchon, F.; Rapha\"el, P., Existence and stability of the log-log blow-up dynamics for the $L^ 2$-critical nonlinear Schr\"odinger equation in a domain, Ann. Henri Poincar\'e 8 (2007), no. 6, 1177--1219.


\bibitem{R1} Rapha\"el, P., Stability of the log-log bound for blow up solutions to the critical non linear Schr\"odinger equation. Math. Ann. 331 (2005), no. 3, 577--609.

\bibitem{RodRaph} Rapha\"el, P.; Rodnianski, I., Stable blow up dynamics for the critical corotational wave map and equivariant Yang Mills, submitted (2009).

%
\bibitem{W1}  Weinstein, M.I., Nonlinear Schr\"odinger equations and sharp interpolation estimates, Comm. Math. Phys.
{\bf 87} (1983), 567---576.

\bibitem{W2} Weinstein, M. I., Lyapunov stability of ground states of nonlinear dispersive evolution equations. Comm. Pure Appl. Math. 39 (1986), no. 1, 51--67.


\end{thebibliography}
\end{document}